\documentclass{amsart}
\usepackage{graphicx}
\usepackage{amsfonts}
\usepackage{amscd}
\usepackage{amssymb}
\usepackage{xypic}
\usepackage{hyperref}
\usepackage{url}
\usepackage{tikz}
\newtheorem{thm}{Theorem}[section]
\newtheorem{cor}[thm]{Corollary}
\newtheorem{lem}[thm]{Lemma}
\newtheorem{prop}[thm]{Proposition}
\newtheorem{conj}[thm]{Conjecture}
\theoremstyle{definition}
\newtheorem{defn}[thm]{Definition}
\theoremstyle{remark}

\numberwithin{equation}{section}
\DeclareMathOperator{\Gal}{Gal}
\DeclareMathOperator{\unr}{unr}
\DeclareMathOperator{\Tr}{Tr}
\DeclareMathOperator{\Sim}{sim}
\DeclareMathOperator{\Hom}{Hom}
\DeclareMathOperator{\Ind}{Ind}
\DeclareMathOperator{\cInd}{c-Ind}
\DeclareMathOperator{\Supp}{Supp}
\DeclareMathOperator{\Span}{Span}
\DeclareMathOperator{\Irr}{Irr}
\DeclareMathOperator{\Par}{Par}
\DeclareMathOperator{\Ad}{Ad}
\DeclareMathOperator{\Contra}{Contra}
\DeclareMathOperator{\Wh}{Wh}
\DeclareMathOperator{\Sh}{Sh}
\DeclareMathOperator{\Id}{Id}

\newcommand{\rvec}[1]{\overrightarrow{#1}}
\newcommand{\lvec}[1]{\overleftarrow{#1}}

\newcommand{\Matrix}[4]{ \left( \begin{array}{cc}  #1 & #2 \\  #3 & #4 \\ \end{array} \right) }

\newcommand{\TMatrix}[9]{ \left( \begin{array}{c|c|c} #1 & #2 & #3 \\ \hline #4 & #5 & #6 \\ \hline #7 & #8 & #9  \\ \end{array} \right) }
\newcommand{\lset}{ \left\{ }
\newcommand{\rset}{ \right\} }
\newcommand{\Norm}{\mathsf N}

\newcommand{\Type}[1]{\mathbf{\mathsf{#1}}}

\newcommand{\Lie}[1]{ {\mathfrak{#1}} }

\newcommand{\alg}[1]{\mathbf{#1}}
\newcommand{\mc}[1]{\mathcal{#1}}

\newcommand{\ZZ}{\mathbb Z}
\newcommand{\rats}{\mathbb Q}

\newcommand{\complex}{\mathbb C}
\newcommand{\CC}{\mathbb C}

\newcommand{\octs}{\mathbb O}

\newcommand{\adeles}{\mathbb A}

\newcommand{\isom}{\cong}

\newcommand{\trans}{\mathsf{T}}

\begin{document}

\bibliographystyle{plain}

\title{Dichotomy for generic supercuspidal representations of $G_2$}%
\author{Gordan Savin and Martin H. Weissman}%
\date{\today}

\address{Dept. of Mathematics, University of California, Santa Cruz, CA 95064}
\email{weissman@ucsc.edu}%


\begin{abstract}
The local Langlands conjectures imply that to every generic supercuspidal irreducible representation of $G_2$ over a $p$-adic field, one can associate a generic supercuspidal irreducible representation of either $PGSp_6$ or$PGL_3$.  We prove this conjectural dichotomy, demonstrating a precise correspondence between certain representations of $G_2$ and other representations of $PGSp_6$ and $PGL_3$.  This correspondence arises from theta correspondences in $E_6$ and $E_7$, analysis of Shalika functionals, and spin L-functions.  Our main result reduces the conjectural Langlands parameterization  of generic supercuspidal irreducible representations of $G_2$ to a single conjecture about the parameterization for $PGSp_6$.

\end{abstract}

\maketitle

\tableofcontents

\section*{Introduction}
Let $k$ be a finite extension of $\rats_p$, $p$ a prime number; we work here with the $k$-points of algebraic groups.  In this paper, we prove a precise correspondence between the generic supercuspidal irreducible representations (abbreviated to ``irreps'') of the exceptional group $G_2$ and certain generic supercuspidal irreps of the classical groups $PGL_3$ and $PGSp_6$.  This correspondence is phrased as a \textbf{dichotomy}, in which to every generic supercuspidal irrep $\tau$ of $G_2$, we associate \textbf{either} a generic supercuspidal irrep $\sigma$ of $PGSp_6$ whose spin L-function has a pole at $s=0$, \textbf{or} a contragredient pair (or self-contragredient singleton) of generic supercuspidal irreps $\{ \rho, \tilde \rho \}$ of $PGL_3$.  Symbolically, we write this dichotomy as a function $\Delta$:
$$\Delta \colon \Irr_g^\circ(G_2) \rightarrow \Irr_{g,Spin}^\circ(PGSp_6) \sqcup \frac{\Irr_g^\circ(PGL_3)}{\Contra}.$$
After constructing this function $\Delta$, we prove that it is bijective when $p \neq 2$.  When $p = 2$, we can prove that $\Delta$ is injective, but there is a subtlety involving self-dual supercuspidal irreps of $PGL_3$ which prevents a proof of bijectivity for now.

This main result is suggested by Langlands' conjectural parameterization of the generic supercuspidal irreps of these groups $G_2$, $PGL_3$, and $PGSp_6$.  For this reason, we demonstrate the precise dichotomy at the level of Langlands parameters in the first section.  The results on Langlands parameters depend essentially on the structure theory of the complex simple groups $G_2(\CC)$, $SL_3(\CC)$, and $Spin_7(\CC)$:  embeddings of $SL_3(\CC)$ into $G_2(\CC)$, embeddings of $G_2(\CC)$ into $Spin_7(\CC)$, and classification of parabolic and other subgroups.

The second section is devoted to the structure theory of certain algebraic groups over the $p$-adic field $k$, including constructions and embeddings of exceptional groups and their parabolic subgroups.  At different times in this paper, we require different embeddings of groups.  As Jacquet modules play a crucial role, we describe in detail two kinds of parabolic subgroups:  minuscule parabolics (arising from Jordan algebras), and two-step parabolic subgroups arising from the structurable algebras of Allison \cite{All1}, \cite{All2}.

The third section provides the definition of the dichotomy map $\Delta$.  Specifically, the dichotomy is realized via theta correspondences using dual pairs $G_2 \times PGL_3 \subset E_6$ and $G_2 \times PGSp_6 \subset E_7$, and the minimal representations (see \cite{GS4}) of $E_6$ and $E_7$.  Such theta correspondences have been studied in the literature -- we mention the results of Ginzburg-Jiang \cite{GJ1}, Gan-Savin \cite{GS1} \cite{GS2}, Savin \cite{Sav}, Magaard-Savin \cite{MS1}, Loke-Savin \cite{LS1}, Gross-Savin \cite{GrS1} \cite{GrS2} .  We refine results of Ginzburg-Rallis-Soudry \cite{GRS}, who first considered the ``tower of theta correspondences'' for $G_2$.  Using extensive analysis of Jacquet modules for the minimal representations of $E_6$ and $E_7$, we are able to demonstrate that this pair of theta correspondences determines a dichotomy function $\Delta$, taking a generic supercuspidal irrep of $G_2$ either to a (unique, up to isomorphism) generic supercuspidal irrep of $PGSp_6$ or to a (unique, up to isomorphism and contragredient) generic supercupsidal irrep of $PGL_3$.

The fourth section is devoted to proving the injectivity of the dichotomy map $\Delta$, through a study of Whittaker and Shalika functionals.  When considering a generic supercuspidal irrep $\rho$ of $PGL_3$, the fibre $\Delta^{-1}(\{\rho, \tilde \rho \})$ has cardinality at most the dimension of a space of Whittaker functionals on $\rho$.  The uniqueness of Whittaker functionals immediately yields injectivity of the dichotomy map in this case.  However, when considering a generic supercuspidal irrep $\sigma$ of $PGSp_6$, the fibre $\Delta^{-1}(\sigma)$ has cardinality equal to the dimension of a space of Shalika functionals on $\sigma$.  Here, the ``Shalika subgroup'' is nearly isomorphic to $GL_2(k[\epsilon] / \epsilon^3)$, embedded appropriately in $GSp_6$.  This subgroup is a cubic analogue of the Shalika subgroup $GL_n(k[\epsilon] / \epsilon^2)$ studied by Jacquet-Rallis \cite{J-R} and others (see \cite{JNQ} for a recent example).  In this fourth section, we prove a result of some independent interest -- the uniqueness of such Shalika functionals for arbitrary supercuspidal irreps of $GSp_6$.  It almost immediately follows that the dichotomy map is injective.

The fifth section is devoted to characterizing the image of the dichotomy map $\Delta$, finishing the proof of a bijection when $p \neq 2$.  The dichotomy map surjects onto the set of generic non-self-contragredient (an automatic condition when $p \neq 2$) supercuspidal irreps of $PGL_3$.  When $p = 2$, we cannot yet exclude the possibility that a generic self-contragredient supercuspidal irrep of $PGL_3$ occurs in the theta correspondence with a generic supercuspidal irrep $\tau$ of $G_2$, and also a generic supercuspidal irrep of $PGSp_6$ occurs in the theta correspondence with the same $\tau$.  In other words, we cannot yet prove that a ``second occurrence'' in a tower of theta lifts is not supercuspidal in residue characteristic two.  From the way we define our dichotomy map $\Delta$, we cannot therefore prove that the image of $\Delta$ includes all self-contragredient supercuspidal irreps of $PGL_3$, though all such irreps of $PGL_3$ occur in a theta correspondence with a generic supercuspidal irrep of $G_2$.

The fifth section focuses on the set of generic supercuspidal irreps of $PGSp_6$ in the image of $\Delta$.  Precisely those generic supercuspidal irreps of $PGSp_6$ with nonvanishing Shalika functional occur in this image.  However, Langlands' conjectures predict another characterization of the image of dichotomy:  a generic supercuspidal irrep $\sigma$ of $PGSp_6$ should occur in the image of dichotomy if and only if its degree 8 spin L-function has a pole at $s=0$.  Thus to characterize the image of dichotomy, we prove that $\sigma$ has a nonvanishing Shalika functional if and only if $L(\sigma, Spin, s)$ has a pole at $s=0$.  This is a local version of the main result of Ginzburg-Jiang \cite{GJ1}.  One direction -- that a nonvanishing Shalika functional implies that the L-function has a pole -- requires an analysis of the minimal representation of $E_8$ (!), the construction of Shahidi \cite{Sha} of the spin L-function and connections to reducibility points for representations of $F_4$ parabolically induced from $GSp_6$.  The other direction -- that if $L(\sigma, Spin, s)$ has a pole at $s=0$ then $\sigma$ has a nonvanishing Shalika functional -- requires the Bump-Ginzburg \cite{BG} integral representation of the spin L-function, results of Vo \cite{Vo} on this L-function, and global methods to demonstrate that the Bump-Ginzburg construction agrees (in its poles) with Shahidi's for the spin L-function.

The dichotomy proven in this paper comes close to proving Langlands' conjectural parameterization of generic supercuspidal irreps of $G_2$ by parameters (representations of the Weil group) with values in $G_2(\CC)$.  Indeed, the dichotomy reduces this parameterization (when $p \neq 2$) to a conjecture related to the Langlands parameterization for $PGSp_6$.  While Langlands parameters for generic irreps of $PGSp_6$ are now known (by functoriality for classical groups, due to Cogdell, Kim, Piatetski-Shapiro, and Shahidi \cite{CKPSS} and the local Langlands correspondence for $GL_7$ by Henniart \cite{Hen2} \cite{Hen3}, Kutzko-Moy \cite{KM}, Harris-Taylor \cite{HT}), it remains to be proven that the currently understood parameterization for $PGSp_6$ is compatible with spin L-functions.  Thus the local Langlands parameterization of generic supercuspidal irreps of $G_2$ is reduced to a single question about the classical group $PGSp_6$ when $p \neq 2$.

Of course, a complete parameterization of supercuspidal irreps of $G_2$ satisfying Langlands' conjectures would require also an analysis of the nongeneric supercuspidal irreps, and the partition of all supercuspidal irreps into L-packets.  For example, many nongeneric representations arise from inner forms $PD^\times$ of $PGL_3$ (see \cite{GS3}),  but we do not address such phenomena in this paper.

\subsection{Conventions}

The letter $k$ will always denote a finite extension of $\rats_p$, where $p$ is a prime number.  A $k$-algebra will always mean a {\em unital} (except for Lie algebras, of course), {\em finite-dimensional} $k$-algebra.  An involution on a $k$-algebra will always mean an anti-automorphism of order $2$, which {\em fixes every element} of $k$.  We do not assume $k$-algebras to be commutative or associative; in fact, non-associative algebras play a central role.  For a $k$-vector space $A$, we write $\Lie{End}_k(A)$ for the Lie algebra of $k$-linear endomorphisms of $A$.

We fix a split Cayley algebra $\octs$ over $k$, in what comes later.  We also fix a smooth, nontrivial, additive character $\psi_k$ of $k$.  From $\psi_k$, we may define a smooth additive character $\psi_\octs$ by:
$$\psi_\octs(\omega) = \psi_k(\Tr(\omega)) \mbox{ for all } \omega \in \octs.$$

We use a boldface letter, such as $\alg{G}$ to denote an algebraic group over $k$.
We use an ordinary letter, such as $G$, to denote the $k$-points of $\alg{G}$, viewed naturally as a topological group.  All representations of such groups $G$ will be assumed to be smooth representations on complex vector spaces.  An \textbf{irrep} of $G$ will mean a smooth irreducible representation of $G$ on a complex vector space.  If $G \rightarrow G'$ is a surjective group homomorphism, and $\pi$ is a representation of $G'$, we often also write $\pi$ for the representation of $G$ arising by pullback.

If $H \subset G$ is a closed subgroup, and $\pi$ is a representation of $H$ on a complex vector space $V$, then we write $\Ind_H^G$ for the represenation of $G$ obtained by smooth (unnormalized) induction:
$$\Ind_H^G \pi = \{ f \in C^\infty(G,V) : f(hg) = \pi(h) f(g) \mbox{ for all } h \in H \}.$$
Here, $C^\infty(G,V)$ denotes the space of uniformly locally constant functions from $G$ to $V$.  Induction is adjoint to restriction, by the appropriate version of Frobenius reciprocity:
$$\Hom_G(\tau, \Ind_H^G \pi) \isom \Hom_H(\tau, \pi)$$
for every smooth representation $\tau$ of $G$ and every smooth representation $\pi$ of $H$.

When $H \backslash G$ is noncompact, it is often more useful to consider the smooth compact induction:
$$\cInd_H^G \pi = \{ f \in \Ind_H^G \pi : \Supp(f) \subset H \cdot K \mbox{ for some compact subset } K \subset G \}.$$
Then $\cInd_H^G \pi$ is again a smooth representation of $G$, and is a subrepresentation of $\Ind_H^G \pi$.

If $\pi$ is a representation of $G$, and $\rho$ is an irrep of $G$, then we say that $\rho$ is a \textbf{constituent} of $\pi$ if $\rho$ is isomorphic to a quotient of a subrepresentation of $\pi$.  However, we almost exclusively work with {\em supercuspidal} constituents in this paper; the injectivity and projectivity of supercuspidal irreps, in the category of smooth representations, implies that when supercuspidal irreps occur as constituents, they also occur as subrepresentations and as quotients.

\subsection{Acknowledgments}
The authors wish to thank the American Institute of Mathematics, where collaboration on this paper began, and the IAS Park City Mathematics Institute for their hospitality and support while this paper was finished.  

The first author was supported by the National Science Foundation grant DMS-0852429 during the preparation of this paper.  The second author wishes to thank the University of Michigan, where parts of this paper were completed.  He also thanks Daniel Bump and Wee Teck Gan for some useful conversations.

\section{Dichotomy of parameters}
\subsection{The local Langlands conjectures}
Recall that $k$ is a finite extension of $\rats_p$, fix an algebraic closure $\bar k$ of $k$, and let $\Gamma = \Gal(\bar k / k)$.  Let $k^{\unr}$ denote the maximal unramified extension of $k$ in $\bar k$.  There is a unique continuous isomorphism from $\Gal(k^{\unr} / k)$ to the profinite group $\hat \ZZ$ which sends the geometric Frobenius to $1$.  This isomorphism yields a surjective homomorphism from $\Gamma$ to $\hat \ZZ$.  The preimage of $\ZZ$ is the subgroup $W_k \subset \Gamma$, called the Weil group of $k$.

The Weil group contains $\Gal(\bar k / k^{\unr})$, and $W_k$ is given the coarsest topology for which $\Gal(\bar k / k^{\unr})$ is an open subgroup endowed with the subspace topology from $\Gal(\bar k / k)$.  Thus there is a short exact sequence of topological groups and continuous homomorphisms:
$$1 \rightarrow \Gal(\bar k / k^{\unr}) \rightarrow W_k \rightarrow \ZZ \rightarrow 1.$$

Let $\alg{G}$ be a semisimple, split, {\em adjoint} algebraic group over $k$, and let $G = \alg{G}(k)$.  Let $\Irr(G)$ denote the set of isomorphism classes of irreducible smooth representations of $G$ on a complex vector space, hereafter called \textbf{irreps} of $G$.  Let $\Irr^\circ(G)$ denote the subset consisting of isomorphism classes of supercuspidal irreps.  Let $\Irr_g(G)$ be the subset consisting of isomorphism classes of generic irreps; the adjective ``generic'' is well-defined, since we assume that $\alg{G}$ is adjoint and split over $k$.  Finally, define $\Irr_g^\circ(G) = \Irr_g(G) \cap \Irr^\circ(G)$ to be the set of isomorphism classes of generic supercuspidal irreps.

Let $\hat G$ denote the complex dual group of $\alg{G}$; thus $\hat G$ is a semisimple, simply-connected complex Lie group.  A \textbf{parameter} for $G$ is a continuous homomorphism $\eta \colon W_k \rightarrow \hat G$ such that $\eta(w)$ is semisimple for all $w \in W_k$.  We do not require the extra structure provided by the Weil-Deligne group here.  A parameter $\eta$ is called \textbf{cuspidal} if $Im(\eta)$ is not contained in any proper parabolic subgroup of $\hat G$.  Let $\Par(G)$ denote the set of parameters, and $\Par^\circ(G)$ the set of cuspidal parameters for $G$.  Note that $\hat G$ acts on the sets $\Par(G)$ and $\Par^\circ(G)$ by conjugation, denoted $\Ad$.

An expectation of the local Langlands conjectures is that there is a ``natural'' bijective parameterization:
$$\Phi(G) \colon \Irr_g^\circ(G) \rightarrow \frac{\Par^\circ(G)}{\Ad(\hat G)},$$
whereby the generic supercuspidal irreps of $G$ are parameterized precisely by the $\hat G$-conjugacy classes of cuspidal parameters.

\subsection{The Dichotomy}

When $\alg{G}_2$ is a simple split algebraic group of type $\Type{G}_2$ over $k$, $\hat G_2 = G_2(\CC)$ is the simple complex Lie group of type $\Type{G}_2$.  In this case, the Langlands conjectures predict that the generic supercuspidal irreps of $G_2$ are parameterized by $\hat G_2$-conjugacy classes of cuspidal parameters.  However, the latter can be related to classical groups as follows.

Let $\octs$ denote an octonion  algebra (also called a Cayley algebra) over $\CC$.  Let $\octs_\circ$ denote the subset of trace zero octonions, and realize $G_2(\CC)$ as the group of $\CC$-algebra automorphisms of $\octs$.  Thus, we find an embedding $G_2(\CC) \hookrightarrow SO_7(\CC) = SO(\octs_\circ, N)$, where $N$ denotes the quadratic norm form on $\octs_\circ$.  As $G_2(\CC)$ is simply connected, this embedding extends to an embedding $G_2(\CC) \hookrightarrow Spin_7(\CC)$.  As $Spin_7(\CC)$ is the complex dual group to $PGSp_6$, we find a natural map
$$\Par^\circ(G_2) \rightarrow \Par(PGSp_6).$$

To determine when the image of a cuspidal parameter for $G_2$ is a {\em cuspidal} parameter for $PGSp_6$, we discuss the maximal parabolic subgroups of $G_2(\CC)$ and $Spin_7(\CC)$.  A \textbf{nil-space} in $\octs_\circ$ is a linear subspace $V \subset \octs_\circ$ such that for all $\alpha, \beta \in V$, $\alpha \cdot \beta = 0$.  An \textbf{isotropic subspace} in $\octs_\circ$ is a linear subspace $V \subset \octs_\circ$ such that $N(\alpha) = 0$ for all $\alpha \in V$.  While for one-dimensional subspaces of $\octs_\circ$, nil-spaces coincide with isotropic spaces, this does not hold in higher dimension.

It is known that every maximal parabolic subgroup of $G_2(\CC)$ is the stabilizer of a one-dimensional or two-dimensional nil-space in $\octs_\circ$ (Theorem 3 of Aschbacher \cite{Asc}).  It is also known that every maximal parabolic subgroup of $Spin_7(\CC)$ is the stabilizer of a one-, two-, or three-dimensional isotropic subspace in $\octs_\circ$.

\begin{prop}
\label{GroupTheory}
Suppose that $P$ is a maximal parabolic subgroup of $Spin_7(\CC)$.  Then either $P \cap G_2(\CC)$ is contained in a maximal parabolic subgroup of $G_2(\CC)$ or $P \cap G_2(\CC)$ is contained in a subgroup of $G_2(\CC)$ isomorphic to $SL_3(\CC)$.
\end{prop}
\proof
There are three cases to consider, depending on whether $P$ stabilizes a one-, two-, or three-dimensional isotropic subspace $V \subset \octs_\circ$:
\begin{description}
\item[$\dim(V) = 1$] 
If $\dim(V) = 1$, then any vector in $V$ has norm zero and trace zero, from which it follows that any vector $\alpha \in V$ satisfies $\alpha^2 = 0$.  It follows that $V$ is a nil-space in $\octs_\circ$.  Thus $P \cap G_2(\CC)$ is the maximal parabolic subgroup of $G_2(\CC)$ stabilizing this nil-space.
\item[$\dim(V) = 2$]
If $\dim(V) = 2$, then  every vector $\alpha \in V$ satisfies $\alpha^2 = 0$.  If $V$ is a nil-space, then $P \cap G_2(\CC)$ is the maximal parabolic subgroup of $G_2(\CC)$ stabilizing this nil-space.  If $V$ is not a nil-space, then there exists a basis $\{ \alpha, \beta \} \subset V$ such that $\alpha \cdot \beta = \gamma \neq 0$.  It follows that $V \cdot V \subset \CC \gamma$.  Therefore, if $g \in P \cap G_2(\CC)$, then $g$ stabilizes not only $V$, but also the line spanned by $\gamma$.

Observe that $\gamma^2 = (\alpha \beta) \cdot (\alpha \beta) = \alpha (\beta \alpha) \beta$ by Moufang identities, and $\beta \alpha = - \alpha \beta$ since $(\alpha + \beta)^2 = 0$.  Hence $\gamma^2 = 0$.  Therefore $P \cap G_2(\CC)$ is contained in the maximal parabolic subgroup stabilizing the nil-line $\CC \gamma$.
\item[$\dim(V) = 3$]
If $\dim(V) = 3$, then we begin by choosing a basis $\{ \alpha, \beta, \gamma \}$ of $V$.  There are two possibilities to consider.  First,  if $\gamma \in \CC(\alpha \cdot \beta)$, then $V \cdot V \subset \CC \gamma$, and $\gamma^2 = 0$.  In this case, $P \cap G_2(\CC)$ stabilizes the nil-line $\CC \gamma$, and hence is contained in a maximal parabolic subgroup of $G_2(\CC)$.

If $\gamma \not \in \CC(\alpha \cdot \beta)$, then $[\alpha, \beta, \gamma] \neq 0$, where the bracket denotes the associator:
$$[\alpha, \beta, \gamma] = (\alpha \beta) \gamma - \alpha (\beta \gamma).$$
In this case, we find that $[V,V,V] \subset \CC \cdot [\alpha, \beta, \gamma]$.  Therefore $P \cap G_2(\CC)$ stabilizes the line $\CC \cdot  [\alpha, \beta, \gamma]$.  The stabilizer of a line in $G_2(\CC)$ is either a maximal parabolic subgroup (if the line is a nil-line), or else a subgroup isomorphic to $SL_3(\CC)$.  Thus $P \cap G_2(\CC)$ is contained in a maximal parabolic subgroup of $G_2(\CC)$ or else is contained in a subgroup isomorphic to $SL_3(\CC)$.
\end{description}
\qed

\begin{prop}
\label{BDSApp}
Suppose that $Q$ is a proper parabolic subgroup of $SL_3(\CC)$.  Then, for any embedding of $SL_3(\CC)$ in $G_2(\CC)$, the image of $Q$ is contained in a maximal parabolic subgroup of $G_2(\CC)$.
\end{prop}
\proof
By the theory of Borel and De Siebenthal \cite{BdS}, every embedding of the full rank subgroup $SL_3(\CC)$ in $G_2(\CC)$ arises from a pair of long roots in the root system of type $\Type{G}_2$.  It follows that a parabolic subgroup $Q \subset SL_3(\CC)$ arises from a single long root in the root system of type $\Type{G}_2$; it follows that $Q$ will be contained in the maximal parabolic subgroup of $G_2(\CC)$ corresponding to this long root.  
\qed

The previous propositions now yield the following dichotomy for parameters:
\begin{thm}
Suppose that $\eta \in \Par^\circ(G_2)$ is a cuspidal parameter for $G_2$.  Let $\eta'$ be the associated parameter for $PGSp_6$ obtained by composing $\eta$ with the inclusion $G_2(\CC) \hookrightarrow Spin_7(\CC)$.  Then either $\eta' \in \Par^\circ(PGSp_6)$, i.e., $\eta'$ is a cuspidal parameter, or else there exists a cuspidal parameter $\eta'' \in \Par^\circ(PGL_3)$ such that $\eta$ is obtained from $\eta''$ via an inclusion $SL_3(\CC) \hookrightarrow G_2(\CC)$.
\end{thm}
\proof
If $\eta'$ is not a cuspidal parameter, then there exists a maximal parabolic subgroup $P \subset Spin_7(\CC)$ such that $Im(\eta') \subset P$.  It follows that $Im(\eta) \subset P \cap G_2$.  Since $\eta$ was assumed cuspidal, we find that $P \cap G_2$ is not contained in any maximal parabolic subgroups of $G_2$.  It follows from Proposition \ref{GroupTheory} that $P$ is the stabilizer of a three-dimensional isotropic subspace of $\octs_\circ$, and $P \cap G_2(\CC)$ is contained in a subgroup isomorphic to $SL_3(\CC)$.

Hence if $\eta' \not \in \Par^\circ(PGSp_6)$, then we find that there exists an embedding $\iota \colon SL_3(\CC) \hookrightarrow G_2(\CC)$, and a parameter $\eta'' \in \Par(PGL_3)$ such that $\eta = \iota \circ \eta''$.  If $\eta''$ were not cuspidal, its image would be contained in a maximal parabolic subgroup of $G_2(\CC)$ by Proposition \ref{BDSApp}, contradicting the cuspidality of $\eta$.  Hence $\eta'' \in \Par^\circ(PGL_3)$.
\qed

This theorem demonstrates that to each $\eta \in \Par^\circ(G_2)$, one may associate a cuspidal parameter $\eta' \in \Par^\circ(PGSp_6)$, or else a cuspidal parameter $\eta'' \in \Par^\circ(PGL_3)$.  Since all embeddings of $G_2(\CC)$ in $Spin_7(\CC)$ are $Spin_7(\CC)$-conjugate, we find that $\eta'$ is uniquely determined (up to $Spin_7(\CC)$-conjugacy) by $\eta$ (up to $G_2(\CC)$-conjugacy).  

Similarly, a cuspidal parameter $\eta \in \Par^\circ(G_2)$, which composes to yield a noncuspidal parameter for $PGSp_6$, yields a cuspidal parameter $\eta'' \in \Par^\circ(PGL_3)$ unique up to $G_2(\CC)$-conjugacy.  Note that all embeddings of $SL_3(\CC)$ into $G_2(\CC)$ are $G_2(\CC)$-conjugate; moreover, the $G_2(\CC)$-conjugacy class of a cuspidal parameter $\eta$ determines the cuspidal parameter $\eta''$ uniquely, up to $SL_3(\CC)$-conjugacy {\em and } outer automorphism.  Namely, the outer automorphism of $SL_3(\CC)$ sending $g$ to $(g^\trans)^{-1}$ is realized by conjugating by an element of $G_2(\CC)$.  The normalizer $N(SL_3(\CC))$ in $G_2(\CC)$ is generated by $SL_3(\CC)$ and an element inducing this outer automorphism.  

Putting these observations together, we find:
\begin{thm}[Dichotomy of parameters]
There is a natural injective dichotomy for the set of cuspidal parameters for $G_2$, modulo $G_2(\CC)$-conjugacy:
$$\frac{\Par^\circ(G_2)}{\Ad(G_2(\CC))} \hookrightarrow \frac{\Par^\circ(PGSp_6)}{\Ad(Spin_7(\CC))} \sqcup \frac{\Par^\circ(PGL_3)}{\Ad(N(SL_3(\CC)))}.$$
\end{thm} 

The image of this dichotomy can also be characterized.  First, we observe the following:
\begin{prop}
Suppose that $\eta'' \in \Par^\circ(PGL_3)$.  Then, for any embedding $\iota \colon SL_3(\CC) \hookrightarrow G_2(\CC)$, $\iota \circ \eta'' \in \Par^\circ(G_2)$.
\end{prop}
\proof
It is clear that $\iota \circ \eta'' \in \Par(G_2)$.   If $P$ is a maximal parabolic subgroup of $G_2$, then $P$ stabilizes a nil-line in $\octs_\circ$ or a nil-plane in $\octs_\circ$.  As a representation of $SL_3(\CC)$, the vector space $\octs_\circ$ decomposes into the direct sum of two irreducible three-dimensional representations, and one trivial representation arising from a $SL_3(\CC)$-fixed line in $\octs_\circ$.  Since there is no nil-line nor nil-plane fixed by $SL_3(\CC)$, we find that $P \cap SL_3(\CC)$ fixes a line or plane in one of the irreducible three-dimensional representations of $SL_3(\CC)$.  Hence $P \cap SL_3(\CC)$ is contained in a maximal parabolic subgroup of $SL_3(\CC)$.  The proposition follows immediately.
\qed

We find that the natural dichotomy for cuspidal parameters for $G_2$ includes all cuspidal parameters for $PGL_3$.  However, not all parameters for $PGSp_6$ occur in this dichotomy.  Perhaps the most convenient way of characterizing the parameters for $PGSp_6$ is through the following:
\begin{prop}
Suppose that $\eta' \in \Par^\circ(PGSp_6)$.  Let $L(\eta', Spin, s)$ denote the Artin-Weil L-function associated to $\eta'$ and the 8-dimensional spin representation of $Spin_7(\CC)$.  Then $L(\eta', Spin, s)$ has a pole at $s = 0$ if and only if the image of $\eta'$ is contained in a subgroup of $Spin_7(\CC)$ isomorphic to $G_2(\CC)$.
\end{prop}
\proof
Let $V$ be an 8-dimensional vector space, on which $Spin_7(\CC)$ acts via the spin representation.  The order of the pole of $L(\eta', Spin, s)$ at $s=0$ is precisely the multiplicity of the trivial representation of $W_k$ for its action on $V$.  Thus, we find that $L(\eta', Spin, s)$ has a pole at $s=0$ if and only if $V$ has a nonzero vector fixed by $W_k$.  

Now, the stabilizer of any nonzero vector $v \in V$ in $Spin_7(\CC)$ is either a proper parabolic subgroup of $Spin_7(\CC)$ or else a group isomorphic to $G_2(\CC)$.  Since we assume that $\eta'$ is a cuspidal parameter, its image in not contained in any proper parabolic subgroups of $Spin_7(\CC)$.  Thus $L(\eta', Spin, s)$ has a pole at $s=0$ if and only if $\eta'(W_k)$ lies in an embedded $G_2(\CC)$ in $Spin_7(\CC)$.  
\qed

Define $\Par_{Spin}^\circ(PGSp_6)$ to be the set of cuspidal parameters $\eta'$ for $PGSp_6$, for which $L(\eta', Spin, s)$ has a pole at $s = 0$.  Then, we find the following perfect dichotomy of parameters:
\begin{thm}
\label{DichPar}
There is a bijective dichotomy for the set of cuspidal parameters for $G_2$, modulo $G_2(\CC)$-conjugacy:
$$\frac{\Par^\circ(G_2)}{\Ad(G_2(\CC))} \leftrightarrow \frac{\Par_{Spin}^\circ(PGSp_6)}{\Ad(Spin_7(\CC))} \sqcup \frac{\Par^\circ(PGL_3)}{\Ad(N(SL_3(\CC)))}.$$
\end{thm} 

\subsection{Dichotomy for Irreps of $G_2$}

The dichotomy for parameters in Theorem \ref{DichPar} suggests, via the local Langlands conjectures, a dichotomy for the generic supercuspidal irreps of $G_2$.  Recall that $\Irr_g^\circ(G)$ denotes the set of isomorphism classes of generic supercuspidal irreps of a (semisimple, adjoint, split) group $G$.  

Define $\Irr_{g,Spin}^\circ(PGSp_6)$ to be the subset of $\Irr_g^\circ(PGSp_6)$, consisting of those irreps $\sigma$ for which Shahidi's degree 8 L-function $L(\tau, Spin, s)$ has a pole at $s = 0$.  The main result of this paper is the following:
\begin{thm}
\label{BigThm}
Dual pair correspondences in the simple split adjoint groups $E_6$ and $E_7$ determine a dichotomy function $\Delta$, which is bijective when $p \neq 2$ and injective when $p = 2$:
$$\Delta: \Irr_g^\circ(G_2) \rightarrow \Irr_{g,Spin}^\circ(PGSp_6) \sqcup \frac{\Irr_g^\circ(PGL_3)}{\Contra},$$
where $\Contra$ denotes the equivalence relation given by contragredience.  
\end{thm}

The existence of such a bijection is directly implied by Langlands conjectures and the dichotomy of parameters in Theorem \ref{DichPar}.  The realization of this bijection through theta correspondences is a result of additional interest, and follows many previous realizations of ``Langlands functoriality'' in theta correspondences.  Conversely, this result can be used to parameterize the generic, supercuspidal representations of $G_2$ over a $p$-adic field, using known and perhaps soon-to-be known parameterizations for $PGL_3$ and $PGSp_6$.

Specifically, the local Langlands conjectures have been proven for $PGL_3$ (for $GL_3$ in fact) by Henniart \cite{Hen}, in the sense that
\begin{prop}
There is a natural (compatible with L-functions and $\epsilon$-factors, among other properties) bijection 
$$\Phi(PGL_3) \colon \frac{\Irr_g^\circ(PGL_3)}{\Contra} \rightarrow \frac{\Par^\circ(PGL_3)}{\Ad(N(SL_3(\CC)))}.$$
In particular, the contragredient on irreps corresponds to the change in parameter given by the outer automorphism of $SL_3(\CC)$.
\end{prop}

While parts of the local Langlands conjectures are open for $PGSp_6$, it appears likely that the following will be proven in the not so distant future.
\begin{conj}
There is a bijection 
$$\Phi(PGSp_6) \colon \Irr_g^\circ(PGSp_6) \rightarrow \frac{\Par^\circ(PGSp_6)}{\Ad(Spin_7(\CC))},$$ in which Shahidi's degree 8 Spin L-function on irreps corresponds to the Artin-Weil degree 8 L-function associated to the Spin representation of $Spin_7(\CC)$.  
\end{conj}

The main theorem of this paper implies:
\begin{thm}
Assuming a parameterization $\Phi(PGSp_6)$ satisfying the previous conjecture, and assuming $p \neq 2$, there is a bijective parameterization:
$$\Phi(G_2) \colon \Irr_g^\circ(G_2) \rightarrow \frac{\Par^\circ(G_2)}{\Ad(G_2(\CC))}.$$
\end{thm}

Of course, there are further properties of this parameterization $\Phi(G_2)$ that should be proven; for example, one hopes that $\Phi(G_2)$ is compatible with L-functions and $\epsilon$-factors of various twists.  

\section{Structure Theory}

There are many constructions of exceptional Lie algebras and algebraic groups.  The construction of Allison \cite{All1} using structurable algebras \cite{All2} (with similarities to earlier constructions of Kantor \cite{Kan}), is well-suited to some needs of this paper.  The construction of Koecher \cite{Koe} using Jordan algebras is well-suited to other needs of this paper.  We recall these constructions of Lie algebras, and associated algebraic groups, in this section.  The constructions here are valid whenever $k$ is a field of characteristic zero (and most likely, when $char(k) \neq 2,3$).

\subsection{Composition, Jordan, and structurable algebras}

\subsubsection{Composition algebras}

\begin{defn}
A composition algebra (sometimes called a Hurwitz algebra) over $k$ is a pair $(C,N)$ where $C$ is $k$-algebra, and $N \colon C \rightarrow k$ is a nondegenerate quadratic form which satisfies $N(xy) = N(x) N(y)$ for all $x,y \in C$.
\end{defn}
Given a composition algebra $(C,N)$ over $k$, we write $N$ also for the associated symmetric bilinear form:
$$N(x,y) = N(x+y) - (N(x) + N(y)).$$
The {\em standard involution} on $C$ is given by:
$$\bar x = N(x,1) - x.$$
The {\em norm} and {\em trace} can be recovered from the standard involution:
$$N(x) = x \bar x, \mbox{ and } \Tr(x) = x + \bar x.$$
According to classification results originating with Hurwitz, composition algebras over $k$ have dimension $1$, $2$, $4$, or $8$ as vector spaces over $k$.  A composition algebra of dimension $8$ will be called a {\em Cayley algebra}.  Composition algebras of dimension $1$ and $2$ are commutative and associative.  Composition algebras of dimension $4$ are associative.  Composition algebras of dimension $8$ are alternative:  if $C$ is a Cayley algebra, and $x,y \in C$, then:
$$(xx)y = x(xy) \mbox{ and } (yx)x = y(xx).$$
Although Cayley algebras are nonassociative, the map $(x,y,z) \mapsto \Tr(xyz)$ defines a {\em trilinear form} on a Cayley algebra $C$; the associative law is not required here since
$$\Tr(x(yz)) = \Tr((xy)z), \mbox{ for all } x,y,z \in C.$$

\subsubsection{Algebras with involution}

Suppose that $A$ is a $k$-algebra with involution (denoted $a \mapsto \bar a$).  For $x,y,z \in A$, we define the following:  first, the left- and right-multiplication endomorphisms are defined by $L_x(y) = xy$ and $R_x(y) = yx$.  Thus $L_x, R_x \in \Lie{End}_k(A)$.  Also, $[x,y] = xy - yx$ is the commutator, and $[x,y,z] = (xy)z - x(yz)$ is the associator.  The involution yields a ternary composition 
$$\{x,y,z \} = (x \bar y) z + (z \bar y) x - (z \bar x) y.$$

This ternary composition yields the endomorphism $V_{x,y} \in \Lie{End}_k(A)$, given by $V_{x,y}(z) = \{ x,y,z \}$.  Finally, define the endomorphism $T_x \in \Lie{End}_k(A)$ is given by $T_x = V_{x,1}$.  Then 
$$T_x = L_x + R_{x - \bar x}.$$

Given a $k$-algebra $A$ with involution, one may consider the {\em hermitian} and {\em skew-hermitian} elements of $A$.  The skew-hermitian (or trace zero) elements of $A$ are:
$$A_\circ = \{ a \in A \mbox{ such that } a + \bar a = 0 \}.$$
The hermitian elements of $A$ are denoted:
$$A_+ = \{ a \in A \mbox{ such that } a = \bar a \}.$$
As a $k$-vector space, one may clearly decompose $A$ as a direct sum:  $A = A_\circ \oplus A_+$.

There is a natural alternating $A_\circ$-valued $k$-bilinear form on $A$, defined by:
$$\langle x,y \rangle = x \bar y - y \bar x = (x \bar y) - \overline{x \bar y}.$$
From this form, one may construct the two-step nilpotent Lie algebra:
$$\Lie{h}(A, A_\circ) = A \oplus A_\circ,$$
whose brackets are given by:
$$[(x,r), (y,s)] = (0,\langle x,y \rangle) = (0, x \bar y - y \bar x) \mbox{ for all } x,y \in A, r,s \in A_\circ.$$
One may also directly construct a two-step unipotent algebraic group:
$$\alg{H}(A, A_\circ) = \lset \left(\begin{array}{ccc}1 & x & z \\0 & 1 & \bar x \\0 & 0 & 1\end{array}\right) : x,z \in A \Tr(z) = N(x) \rset,$$
where composition is given by the usual rules for matrix multiplication and the composition in the algebra $A$.

\subsubsection{Jordan algebras}

Let $C$ be a composition algebra over $k$.  Without reviewing the general theory of Jordan algebras, we mention and describe the Jordan algebra $J_C$ of Hermitian-symmetric 3 by 3 matrices with entries in $C$:
$$J_C = \lset  \left(\begin{array}{ccc} a & \gamma & \bar \beta  \\ \bar \gamma & b & \alpha \\ \beta & \bar \alpha & c\end{array}\right) : a,b,c \in k, \alpha, \beta, \gamma \in C \rset.$$ 
On $J_C$, there is the Jordan composition:
$$j_1 \circ j_2 = \frac{1}{2} \cdot (j_1 j_2 + j_2 j_1),$$
where ordinary matrix multiplication is used on the right side above.

But more importantly for our purposes are the quadratic adjoint, cubic determinant, and cross product.  The quadratic adjoint is defined by (following notation of Section 2.4 of \cite{Kru}):
$$\left(\begin{array}{ccc} a & \gamma & \bar \beta  \\ \bar \gamma & b & \alpha \\ \beta & \bar \alpha & c\end{array}\right)^\sharp = \left(\begin{array}{ccc} bc-N(\alpha) & \bar \beta \bar \alpha - c \gamma & \gamma \alpha - b \bar \beta  \\ \alpha \beta - c \bar \gamma & ca - N(\beta) & \bar \gamma \bar \beta - a \alpha \\ \bar \alpha \bar \gamma - b \beta & \beta \gamma - a \alpha & ab - N(\gamma) \end{array}\right).$$
The cross product is the linearization of this quadratic adjoint:
$$j_1 \times j_2 = (j_1 + j_2)^\sharp - (j_1^\sharp + j_2^\sharp).$$
There exists a unique cubic form $\Norm \colon J_C \rightarrow k$, for which
$$j \times j^\sharp = \Norm(j) \cdot j, \mbox{ for all } j \in J_C.$$
There is a natural nondegenerate trace pairing
$$T(j,j') = \Tr(j \circ j').$$

\subsubsection{Structurable algebras}
We define and discuss {\em structurable algebras} here, following the foundational work of Allison \cite{All2} very closely.
\begin{defn}
A $k$-algebra $A$ with involution is called a {\em structurable algebra} if, for all $x,y,z \in A$, the following (quartic polynomial) identity holds:
$$[T_z, V_{x,y} ] = V_{T_z x, y} - V_{x, T_{\bar z} y}.$$
\end{defn}
Such an algebra satisfies:
$$[r,x,y] = [x,y,r] = -[x,r,y], \mbox{ for all } x,y \in A, r \in A_\circ.$$
Let $\Lie{Der}(A)$ denote the Lie algebra over $k$, consisting of derivations of $A$ which commute with the involution.  These are $k$-endomorphisms $D$ of $A$, which satisfy the following identities:
$$D(xy) = (Dx) y + x (Dy), \mbox{ and } D(\bar x) = \overline{Dx} \mbox{ for all } x,y \in A.$$

Important examples of structurable algebras include tensor products of composition algebras.  These have been studied extensively by Allison in \cite{All3}, who proves:
\begin{prop}
Suppose that $B$ and $C$ are composition algebras.  Then $B \otimes_k C$, with the tensor product algebra structure and involution, is a structurable algebra.
\end{prop}

When $A = B \otimes_k C$ is a tensor product of two composition algebras, as above, one may check directly that:
$$A_\circ = (B_\circ \otimes_k k) \oplus (k \otimes_k C_\circ) \isom B_\circ \oplus C_\circ.$$
In this way $\alg{H}(A, A_\circ)$ has central subgroup $B_\circ \oplus C_\circ$, and abelian quotient $B \otimes_k C$.

Another important example of a structurable algebra, from Section 8 of \cite{All2}, is given by a construction of Freudenthal.  From a composition algebra $C$, and the resulting Jordan algebra $J_C$, consider the $k$-vector space
$$F_C = \lset \Matrix{a}{j}{j'}{d} : a,d \in k \mbox{ and } j,j' \in J_C \rset.$$
This space has a natural $k$-algebra structure given by
$$\Matrix{a_1}{j_1}{j_1'}{d_1} \cdot \Matrix{a_2}{j_2}{j_2'}{d_2} = \Matrix{a_1 a_2 + T(j_1, j_2')}{a_1 j_2 + d_2 j_1 + j_1' \times j_2'}{a_2 j_1' + a_2 j_2' + j_1 \times j_2}{T(j_2,j_1') + d_1 d_2}.$$
An involution on $F_C$ is given by
$$\overline{ \Matrix{a}{j}{j'}{d} } = \Matrix{d}{j}{j'}{a}.$$
In \cite{All2}, Allison proves (in fact, he proves much more) that
\begin{prop}
If $C$ is any composition algebra then $F_C$, with product and involution given above, is a structurable algebra.
\end{prop}
Note that the trace zero elements of $F_C$ form a one-dimensional subspace.
$$(F_C)_\circ = \lset \Matrix{a}{0}{0}{-a} : a \in k \rset.$$

\subsection{Lie algebras}
From Jordan algebras and structurable algebras, we may follow constructions of Tits-Koecher and Allison to construct certain Lie algebras over $k$.  We review these constructions here.

\subsubsection{Lie algebras from Jordan algebras}
Suppose that $J$ is a semisimple Jordan algebra.  Then constructions of Tits, Kantor, or Koecher \cite{Koe} (whom we follow here) yield a graded Lie algebra:
$$\Lie{g}_J = \Lie{g}_J^{(-1)} \oplus \Lie{g}_J^{(0)} \oplus \Lie{g}_J^{(1)},$$
where $\Lie{g}_J^{(0)} = \Lie{Str}(J)$ is the subalgebra of $\Lie{End}_k(J)$ generated by derivations of $J$ and left Jordan multiplications $L_j$ (for $j \in J$) and $\Lie{g}_J^{(\pm 1)}$ is identified with $J$ as a $k$-vector space.  The Lie bracket on $\Lie{g}_J$ is given by the following:
\begin{itemize}
\item
For all $j \in J$, let $\alpha_\pm(j)$ denote the element of $\Lie{g}_J^{(\pm 1)}$ associated to $j$.  The Lie algebras $\Lie{g}_J^{(\pm)}$ are abelian, i.e.,
$$[\alpha_+(j), \alpha_+(j')] = [\alpha_-(j), \alpha_-(j')] = 0, \mbox{ for all } j,j' \in J.$$
\item
For all $X \in \Lie{g}_J^{(0)}$ and all $j \in J$, we define Lie brackets by
$$[X, \alpha_+(j)] = \alpha_+ \left( X(j) \right), \mbox{ recalling that } X \in \Lie{Str}(J) \subset \Lie{End}(J).$$
Also, we define
$$[X, \alpha_-(j)] = \alpha_- \left(-X^\ast(j) \right),$$
where $X^\ast$ denotes the adjoint endomorphism of $J$, with respect to the trace pairing on $J$.
\item
For all $j,j' \in J$, we define
$$[\alpha_+(j), \alpha_-(j')] = 2 \left( L_{j \circ j'} + [L_j, L_{j'}] \right) \in \Lie{Str}(J) = \Lie{g}_J^{(0)}.$$
\end{itemize}

In this way, the Lie algebra $\Lie{g}_J$ is naturally endowed with a parabolic subalgebra $\Lie{p}_J = \Lie{g}_J^{(0)} \oplus \Lie{g}_J^{(1)}$ with abelian nilradical $\Lie{u}_J = \Lie{g}_J^{(1)} = J$.

\subsubsection{Lie algebras from structurable algebras}

Suppose that $A$ is a structurable algebra.  Following Allison \cite{All1}, let $\Lie{Strl}(A)$ be the $k$-subspace of $\Lie{End}_k(A)$ spanned by $\Lie{Der}(A)$ and endomorphisms of the form $T_a$ for $a \in A$.  Then $\Lie{Strl}(A)$ is a Lie subalgebra of $\Lie{End}_k(A)$, and contains $\Lie{Der}(A)$ as a Lie subalgebra.  Given $X \in \Lie{Strl}(A)$, $X \in \Lie{Der}(A)$ if and only if $X(1) = 0$.

Many elements of $\Lie{Strl}(A)$ arise from ``inner'' endomorphisms, i.e., endomorphisms arising directly from the composition and involution on $A$.
\begin{itemize}
\item
For all $r \in A$, $T_r \in \Lie{Strl}(A)$ by definition.
\item
For all $x,y \in A$, define a derivation of $A$ by:
$$D_{x,y}(z) = \frac{1}{3} \left[ [x,y] + [\bar x, \bar y], z \right] + [z,y,x] - [z,\bar x, \bar y],$$
for all $z \in A$.   From Section 1 of \cite{All1}, $D_{x,y} \in \Lie{Der}(A) \subset \Lie{Strl}(A) \subset \Lie{End}_k(A)$.
\item
For all $x,y \in A$, one has:
$$V_{x,y} = \frac{1}{3} T_{2x y + \bar yx - \bar x y + y \bar x} + D_{x,\bar y}.$$
Hence $V_{x,y} \in \Lie{Strl}(A)$.
\item
For all $r,s \in A_\circ$,
$$L_r L_s = T_{rs} - V_{r,s}.$$
Hence $L_r L_s \in \Lie{Strl}(A)$.
\end{itemize}
Following \cite{All1}, \cite{All2}, we write $\Lie{Instrl}(A)$ for the subspace of $\Lie{Strl}(A)$ spanned by $V_{x,y}$ for all $x,y \in A$.  We write $\Lie{Inder}(A)$ for the subspace of $\Lie{Der}(A)$ spanned by $D_{x,y}$ for all $x,y \in A$.  Then $\Lie{Instrl}(A)$ is an ideal in $\Lie{Strl}(A)$, and $\Lie{Inder}(A)$ is an ideal in $\Lie{Der}(A)$.  The subspace $\Lie{L}(A)$ spanned by $L_r L_s$ for all $r,s \in A_\circ$ is an ideal in $\Lie{Strl}(A)$, and there is a chain of inclusions:
$$\Lie{L}(A) \subset \Lie{Instrl}(A) \subset \Lie{Strl}(A).$$

For all $X \in \Lie{Strl}(A)$, define $X^\epsilon$ and $X^\delta$ by:
$$X^\epsilon = X - T_{X(1) + \overline{X(1)}}, \mbox{ and } X^\delta = X + R_{\overline{X(1)}}.$$
Then, $X \mapsto X^\epsilon$ is an automorphism of the Lie algebra $\Lie{Strl}(A)$ of order $2$.  The element $X^\delta \in \Lie{End}_k(A)$ preserves the subspace $A_\circ \subset A$, and the resulting map $X \mapsto X^\delta$ is a Lie algebra representation:
$$\Lie{Strl}(A) \rightarrow \Lie{End}_k(A_\circ).$$

From a structurable algebra $A$, Allison (in \cite{All1}) constructs a Lie algebra, with similarities to earlier work of Kantor \cite{Kan}.  This Lie algebra, $\Lie{g}_A$ is constructed with a $\ZZ$-grading, vanishing outside degrees $-2,-1,0,1,2$.  In these degrees, the Lie algebra is constructed as follows:
\begin{itemize}
\item
In degree $\pm 2$, we define $\Lie{g}_A^{(\pm 2)} = A_\circ$.  For all $r \in A_\circ$, we write $\zeta_\pm(r)$ for the corresponding element of $\Lie{g}_A^{(\pm 2)}$.
\item
In degree $\pm 1$, we define $\Lie{g}_A^{(\pm 1)} = A$.  For all $x \in A$, we write $\eta_\pm(x)$ for the corresponding element of $\Lie{g}_A^{(\pm 1)}$.
\item
In degree zero, we define $\Lie{g}_A^{(0)} = \Lie{Instrl}(A)$.
\end{itemize}
The brackets on the Lie algebra $\Lie{g}_A = \bigoplus_{i=-2}^2 \Lie{g}_A^{(i)}$ are defined by the following identities:
\begin{itemize}
\item
The space $\Lie{u}_A = \Lie{g}_A^{(1)} \oplus \Lie{g}_A^{(2)} = A \oplus A_\circ$ is identified as a Lie algebra with $\Lie{h}(A, A_\circ)$.  In other words,
$$[\eta_+(x) + \zeta_+(r), \eta_+(y) + \zeta_+(s)] = \zeta_+(x \bar y - y \bar x),$$
for all $x,y \in A$, and $r,s \in A_\circ$.  The bracket on $\Lie{g}_A^{(-1)} \oplus \Lie{g}_A^{(-2)}$ is defined in the same way:
$$[\eta_-(x) + \zeta_-(r), \eta_-(y) + \zeta_-(s)] = \zeta_-(x \bar y - y \bar x),$$
\item
The elements $X \in \Lie{g}_A^{(0)} = \Lie{Instrl}(A)$ are endomorphisms of the $k$-vector space $A$.  For such elements, $X^\delta$ is an endomorphism of the $k$-vector space $A_\circ$.  Hence, for all $X \in \Lie{g}_A^{(0)}$, it makes sense to define:
$$[X, \eta_+(x) + \zeta_+(r)] = \eta_+(X(x)) + \zeta_+(X^\delta(r)).$$
Recalling that $\epsilon$ is an automorphism of $\Lie{Instrl}(A)$ of order two, it makes sense to define:
$$[X, \eta_-(x) + \zeta_-(r)] = \eta_-(X^\epsilon(x)) + \zeta_-(X^{\epsilon \delta}(r)).$$
\item
For $x,y \in A$, and $r,s \in A_\circ$, define:
$$[\eta_+(x) + \zeta_+(r), \eta_-(y) + \zeta_-(s)] = -\eta_-(sx) + (V_{x,y} + L_r L_s) + \eta_+(ry).$$
\end{itemize}
These identities suffice to determine the Lie algebra structure on all of $\Lie{g}_A$.  Note that $\Lie{g}_A$ is naturally endowed with a parabolic subalgebra 
$$\Lie{p}_A = \Lie{g}_A^{(0)} \oplus \Lie{g}_A^{(1)} \oplus \Lie{g}_A^{(2)},$$
with unipotent radical $\Lie{u}_A$ with center $\Lie{z}_A$.  Furthermore, $\Lie{z}_A$ is identified with $A_\circ$, and $\Lie{u}_A / \Lie{z}_A$ is identified with $A$.

\subsection{Algebraic Groups}

Consider a Jordan algebra $J$, and the Koecher Lie algebra $\Lie{g}_J$ constructed earlier.  Define an algebraic group $\alg{G}_J$ over $k$ as the algebraic subgroup of $\alg{GL}(\Lie{g}_J)$ preserving the Lie bracket and a Killing form.  The three-term grading on $\Lie{g}_J$ yields a parabolic subgroup $\alg{P}_J$ with abelian unipotent radical $\alg{U}_J$, whose $k$-points are identified with $J$ itself.

If $J \subset K$ is an embedding of Jordan algebras (i.e., $J$ and $K$ are Jordan algebras, and $J$ is embedded as a sub-$k$-algebra of $K$), then $\Lie{g}_J$ is naturally a \textbf{graded} Lie subalgebra of $\Lie{g}_K$.  This follows quickly from the fact, proven by Jacobson \cite{Jac} that all derivations of the semisimple Jordan algebras considered are inner derivations -- hence these derivation algebras extend to derivations of larger semisimple Jordan algebras. 

Since $\alg{G}_K$ is an algebraic group with Lie algebra $\Lie{g}_K$, and $\Lie{g}_J$ is a semisimple Lie subalgebra of $\Lie{g}_K$, there is an algebraic subgroup $\alg{G}_J' \subset \alg{G}_K$ and an isogeny $\iota : \alg{G}_J' \rightarrow \alg{G}_J$ (where $\alg{G}_J$ is the adjoint algebraic group associated to $\Lie{g}_J$).  Let $\alg{P}_J' = \iota^{-1}(\alg{P}_J)$ and let $\alg{U}_J'$ be the neutral component of $\iota^{-1}(\alg{U}_J)$.

The embedding of algebraic groups $\alg{G}_J' \subset \alg{G}_K$, is \textbf{compatible with parabolics}:
$$\alg{P}_K \cap \alg{G}_J' = \alg{P}_J', \mbox{ and } \alg{U}_K \cap \alg{G}_J' = \alg{U}_J'.$$

Similarly, consider a structurable algebra $A$, and Allison's Lie algebra $\Lie{g}_A$ constructed previously.  Define an algebraic group $\alg{G}_A$ over $k$ as the algebraic subgroup of $\alg{GL}(\Lie{g}_A)$ preserving the Lie bracket and a Killing form.  The five-term grading on $\Lie{g}_A$ yields a parabolic subgroup $\alg{P}_A$ with two-step unipotent radical $\alg{U}_A \supset \alg{Z}_A$.  The $k$-points of the center $\alg{Z}_A$ can be identified with $A_\circ$, and the $k$-points of the quotient $\alg{U}_A / \alg{Z}_A$ can be identified with $A$ itself.

If $A \subset B$ is an embedding of structurable algebras (i.e., $A$ and $B$ are structurable algebras, and $A$ is embedded as a sub-$k$-algebra with involution into $B$), then $\Lie{g}_A$ is naturally a \textbf{graded} Lie subalgebra of $\Lie{g}_B$ (since elements of $\Lie{Instrl}(A) \subset \Lie{End}_k(A)$ extend naturally to elements of $\Lie{Instrl}(B) \subset \Lie{End}_k(B)$).  As before, one obtains an algebraic subgroup $\alg{G}_A' \subset \alg{G}_B$ together with an isogeny $\iota : \alg{G}_A' \rightarrow \alg{G}_A$.  This embedding is \textbf{compatible with parabolics}:
$$\alg{P}_B \cap \alg{G}_A' = \alg{P}_A', \quad \alg{U}_B \cap \alg{G}_A' = \alg{U}_A', \quad \alg{Z}_B \cap \alg{G}_A' = \alg{Z}_A'.$$

\subsubsection{Automorphisms of composition algebras}

Fix a ``complete chain'' of composition algebras $k \subset K \subset B \subset C$, where $K$,$B$,$C$ are composition algebras of $k$-dimension $2$,$4$,$8$, respectively.   Some interesting algebraic groups arise as automorphism groups of extensions of composition algebras.  Namely, if $H \subset E$ is an embedding of composition algebras over $k$, then let $\alg{Aut}_{E/H}$ denote the algebraic subgroup of $\alg{GL}(E)$ preserving the algebra structure and fixing the subalgebra $H$ element-wise.  For example, $\alg{Aut}_{C/k}$ is a absolutely simple group of type $\Type{G}_2$, and $\alg{Aut}_{C/K}$ is a simply-connected absolutely simple group of type $\Type{A}_2$.  $\alg{Aut}_{B/k}$ is an adjoint absolutely simple group of type $\Type{A}_1$, and $\alg{Aut}_{C/B}$ is a simply-connected absolutely simple group of type $\Type{A}_1$.

\subsubsection{Groups from Jordan algebras}
\label{GJ}
The chain of composition algebras $k \subset K \subset B \subset C$ yields a chain of Jordan algebras $J_k \subset J_K \subset J_B \subset J_C$.  The associated algebraic groups $\alg{G}_J$ with parabolic subgroup $\alg{P}_J = \alg{L}_J \alg{U}_J$ are tabulated below:
\begin{center}
\begin{tabular}{c|cccc}
Composition Algebra & $k$ & $K$ & $B$ & $C$ \\ \hline
Dimension of $J$ & $6$ & $9$ & $15$ & $27$ \\
Type of $\alg{G}_J$ & $\Type{C}_3$ & $\Type{A}_5$ & $\Type{D}_6$ & $\Type{E}_7$ \\
Type of Levi $\alg{L}_J$ &  $\Type{A}_2$ & $\Type{A}_2 \times \Type{A}_2$ & $\Type{A}_5$ & $\Type{E}_6$
\end{tabular}
\end{center}

Given an embedding $H \subset E$ of composition algebras, we find an embedding of Jordan algebras $J_H \subset J_E$, and a subgroup $\alg{G}_{J_H}'$ of $\alg{G}_{J_E}$ together with an isogeny $\alg{G}_{J_H}' \rightarrow \alg{G}_{J_H}$ .  Moreover, the subgroup $\alg{G}_{J_H}'$ commutes with $\alg{Aut}_{E/H}$, naturally embedded in $\alg{G}_{J_E}$.  In this way we find many commuting pairs of subgroups.  We label them only by their type, leaving the precise determination of isogeny type up to the reader.
\begin{center}
\begin{tabular}{c|c|cc}
$H$ & $E$ & $\alg{Aut}_{E/H} \times \alg{G}_{J_H}'$ & $\alg{G}_{J_{E}}$ \\ \hline
$k$ & $C$ & $\Type{G}_2 \times \Type{C}_3$ & $\Type{E}_7$ \\
$K$ & $C$ & $\Type{A}_2 \times \Type{A}_5$ &  $\Type{E}_7$ \\
$k$ & $B$ & $\Type{A}_1 \times \Type{C}_3$ & $\Type{E}_6$ \\
$B$ & $C$ & $\Type{A}_1 \times \Type{E}_6$ & $\Type{E}_7$ 
\end{tabular}
\end{center}

\subsubsection{Tensor products of composition algebras}
\label{TC}
The chain of Hurwitz algebras yields embeddings of structurable algebras from which we examine:
$$k \otimes B \subset k \otimes C \subset K \otimes C \subset B \otimes C \subset C \otimes C.$$
This yields embeddings (up to isogeny) of algebraic groups $\alg{G}_A$, compatible with two-step parabolic subgroups $\alg{P}_A = \alg{L}_A \alg{U}_A$.  We tabulate some possibilities in the following:
\begin{center}
\begin{tabular}{c|ccccc}
$A$ & $k \otimes B$ & $k \otimes C$ & $K \otimes C$ & $B \otimes C$ & $C \otimes C$ \\ \hline
Type of $\alg{G}_A$ & $\Type{C}_3$ & $\Type{F}_4$ & $\Type{E}_6$ & $\Type{E}_7$ & $\Type{E}_8$ \\
Type of Levi $\alg{L}_A$ & $\Type{A}_1 \times \Type{A}_1$ & $\Type{B}_3$ & $\Type{A}_1 \times \Type{A}_2 \times \Type{A}_2$ & $\Type{D}_5 \times \Type{A}_1$ & $\Type{D}_7$ \\
Dimension of $\alg{U}_A / \alg{Z}_A$ & 4 & 8 & 16 & 32 & 64 \\
Dimension of $\alg{Z}_A $ & 3 & 7 & 8 & 10 & 14
\end{tabular}  
\end{center}

This construction also realizes some well-known dual reductive pairs.  Consider three composition algebras $H, H', E$, such that $H \subset E$.  Then, $\alg{Aut}_{E/H}$ naturally acts on the Lie algebra $\Lie{g}_{E \otimes H'}$ and $\alg{Aut}_{E/H}$ fixes the elements of the subalgebra $\Lie{g}_{H \otimes H'}$.  This yields a homomorphism of algebraic groups:
$$\alg{Aut}_{E/H} \times \alg{G}_{H \otimes H'}' \hookrightarrow \alg{G}_{E \otimes H'}.$$
In particular, we find many commuting pairs of subgroups:
\begin{center}
\begin{tabular}{c|c|c|cc}
$H$ & $E$ & $H'$ & $\alg{Aut}_{E/H} \times \alg{G}_{H \otimes H'}'$ & $\alg{G}_{E \otimes H'}$ \\ \hline
$k$ & $C$ & $C$ & $\Type{G}_2 \times \Type{F}_4$ & $\Type{E}_8$ \\
$K$ & $C$ & $C$ & $\Type{A}_2 \times \Type{E}_6$ &  $\Type{E}_8$ \\
$B$ & $C$ & $C$ & $\Type{A}_1 \times \Type{E}_7$ & $\Type{E}_8$ \\
$k$ & $C$ & $B$ & $\Type{G}_2 \times \Type{C}_3$ & $\Type{E}_7$ \\
$B$ & $C$ & $k$ & $\Type{A}_1 \times \Type{C}_3$ & $\Type{F}_4$
\end{tabular}
\end{center}
While such exceptional dual pairs occur often in the literature, this construction is convenient for at least two reasons:  first, it gives dual pairs of nonsplit subgroups which may be otherwise difficult to contruct.  Second, the embeddings are compatible with a distinguished parabolic subgroup, which is convenient later for computation of Jacquet modules.

\subsubsection{Freudenthal structurable algebras}
\label{FA}
Finally, we recall that associated to the chain of composition algebras $k \subset K \subset B \subset C$, there is a chain of Jordan algebras $J_k \subset J_K \subset J_B \subset J_C$, and thus a chain of structurable algebras of Freudenthal type:
$$F_k \subset F_K \subset F_B \subset F_C.$$
Each one of these structurable algebras has a one-dimensional subspace of trace zero elements.  Allison's construction yields embeddings of algebraic groups (up to some isogeny)
$$\alg{G}_{F_k}' \subset \alg{G}_{F_K}' \subset \alg{G}_{F_B}' \subset \alg{G}_{F_C},$$
compatible with two-step ``Heisenberg'' parabolic subgroups $\alg{P}_F = \alg{L}_F \alg{U}_F$.  We tabulate the possibilities in the following:
\begin{center}
\begin{tabular}{c|cccc}
Jordan Algebra  & $J_k$ & $J_K$ & $J_B$ & $J_C$ \\ \hline
Dimension of $F$ & $14$ & $20$ & $32$ & $56$ \\
Type of $\alg{G}_F$ & $\Type{F}_4$ & $\Type{E}_6$ & $\Type{E}_7$ & $\Type{E}_8$ \\
Type of Levi $\alg{L}_F$ & $\Type{C}_3$ & $\Type{A}_5$ & $\Type{D}_6$ & $\Type{E}_7$ 
\end{tabular}
\end{center}

\section{Theta correspondence}
The main result to be proven in this paper is a bijective dichotomy:
$$\Irr_g^\circ(G_2) \leftrightarrow \Irr_{g,Spin}^\circ(PGSp_6) \sqcup \frac{\Irr_g^\circ(PGL_3)}{\Contra}.$$
In this section, we begin the proof of this main result.  We use theta correspondences in $E_6$ and $E_7$ to describe maps for the above dichotomy.  Beginning with a generic supercuspidal irrep $\tau$ of $G_2$,
\begin{itemize}
\item
We will define $\rvec{\Theta}_6(\tau)$, a representation of $PGL_3$, and $\rvec{\Theta}_7(\tau)$, a representation of $PGSp_6$.
\item
If $\rvec{\Theta}_6(\tau) = 0$, then $\rvec{\Theta}_7(\tau)$ has a unique generic supercuspidal irreducible subrepresentation.
\item
Otherwise, and \textbf{if $p \neq 2$}, then $\rvec{\Theta}_6(\tau)$ has a unique, up to contragredience, generic supercuspidal irreducible subrepresentation.  Even if $p = 2$, $\rvec{\Theta}_6(\tau)$ is a multiplicity-free supercuspidal representation of $PGL_3$.
\end{itemize}
By establishing these facts, we establish a map in this section, when $p \neq 2$:
$$\Delta : \Irr_g^\circ(G_2) \rightarrow \Irr_{g}^\circ(PGSp_6) \sqcup \frac{\Irr_g^\circ(PGL_3)}{\Contra}.$$
where $\Delta(\tau)$ is either the unique (up to isomorphism) generic supercuspidal subrepresentation of $\rvec{\Theta}_7(\tau)$ or the unique (up to isomorphism and contragredience) generic supercuspidal subrepresentation of $\rvec{\Theta}_6(\tau)$.  

\subsection{Minimal representations}

Let $\Pi_6$ and $\Pi_7$ denote the minimal representations of the adjoint simple split groups $E_6$ and $E_7$, respectively (we refer to \cite{GS4} for definitions and properties of minimal representations).   Let $\sigma$ be a supercuspidal irrep of $PGSp_6$, let $\tau$ be a supercuspidal irrep of $G_2$, and let $\rho$ be a supercuspidal irrep of $PGL_3$.  We define the following:
$$\lvec{\Theta}_7(\sigma) = \Hom_{PGSp_6}(\sigma, \Pi_7),  \mbox{ and } \rvec{\Theta}_7(\tau) = \Hom_{G_2}(\tau, \Pi_7).$$
Of course, we view $\lvec{\Theta}_7(\sigma)$ as a representation of $G_2$, and $\rvec{\Theta}_7(\tau)$ as a representation of $PGSp_6$, via the dual pair (see Section \ref{TC}):
$$\alg{PGSp}_6 \times \alg{G}_2 \rightarrow \alg{E}_7.$$
Observe here that we consider embeddings of $\sigma$ and $\tau$ as subrepresentations rather than the more commonly used quotients; however, the injectivity and projectivity of supercuspidals in the category of smooth representations implies that nothing is lost.  Note that $\sigma \boxtimes \lvec{\Theta}_7(\sigma)$ is naturally a $(PGSp_6, \sigma)$-isotypic subspace of $\Pi_7$, and $\rvec{\Theta}_7(\tau) \boxtimes \tau$ is naturally a $(G_2, \tau)$-isotypic subspace of $\Pi_7$.

Similarly, we define
$$\lvec{\Theta}_6(\rho) = \Hom_{PGL_3}(\rho, \Pi_6), \mbox{ and } \rvec{\Theta}_6(\tau) = \Hom_{G_2}(\tau, \Pi_6).$$
Here, we view $\lvec{\Theta}_6(\rho)$ as a representation of $G_2$ and $\rvec{\Theta}_6(\tau)$ as a representation of $PGL_3$, via the dual pair
$$\alg{PGL}_3 \times \alg{G}_2 \hookrightarrow \alg{E}_6.$$
Observe that $\rho \boxtimes \lvec{\Theta}_6(\rho)$ is naturally a $(PGL_3, \rho)$-isotypic subspace of $\Pi_6$, and $\rvec{\Theta}_6(\tau) \boxtimes \tau$ is naturally a $(G_2, \tau)$-isotypic subspace of $\Pi_6$.

\subsection{Whittaker functionals}

Let $N_{2}$ and $N_{3}$ be the unipotent radicals of Borel subgroups of $G_{2}$ and $PGSp_{6}$, respectively. Let $\psi_{2} \colon N_{2}\rightarrow \mathbb C^{\times}$ and $\psi_{3} \colon  N_{3}\rightarrow \mathbb C^{\times}$ be generic (principal) characters.  Since $G_{2}$ and $PGSp_{6}$ are of adjoint type these characters are unique up to conjugation by the tori of the respective Borel subgroups.   For this reason, $\tau$ and $\sigma$ are unambiguously called {\em generic} (rather than $\psi_2$-generic and $\psi_3$-generic) if $\tau_{N_{2},\psi_{2}}\neq 0$ and $\sigma_{N_{3},\psi_{3}} \neq 0$ respectively.

More generally, when $\alg{G}$ is a split adjoint semisimple group over $k$, and $\pi$ is a smooth representation of $G$, we write $\Wh_G(\pi)$ for the space of Whittaker functionals on $\pi$, with respect to some maximal unipotent subgroup $\alg{N}$ of $\alg{G}$ and principal character $\psi$ of $N$:
$$\Wh_G(\pi) = \Hom_N(\pi, \psi).$$
 Thus, $\tau$ is called generic if $\Wh_{G_2}(\tau) \neq 0$ and $\sigma$ is called generic if $\Wh_{PGSp_6}(\sigma) \neq 0$.

It is important to recall a few equivalent formulations of Whittaker functionals and genericity.  While well known, a good treatment can be found in the work of Casselman and Shalika \cite{CS}.  First, since $\pi_{N, \psi}$ is the maximal quotient on which $N$ acts via $\psi$, we find canonical isomorphisms
$$\Wh_G(\pi) = \Hom_N(\pi, \psi) \isom \Hom_N(\pi_{N, \psi}, \psi) \isom \Hom_\CC(\pi_{N, \psi}, \CC).$$
In particular, $\dim(\Wh_G(\pi)) = \dim(\pi_{N, \psi})$ if one of these vector spaces is finite-dimensional.

Next, by Frobenius reciprocity, observe that
$$\Wh_G(\pi) = \Hom_N(\pi, \psi) \isom \Hom_G(\pi, \Ind_N^G \psi).$$
If $\pi$ is a generic irrep of $G$, so $\Wh_G(\pi)$ is nonzero, then $\pi$ embeds as a subrepresentation of $\Ind_N^G \psi$.  The image of $\pi$ via such an embedding is uniquely determined by $\pi$; it is called the Whittaker model of $\pi$.

On the other hand, we often consider the \textbf{Gelfand-Graev} representation $\cInd_N^G \psi$; since this is a submodule of $\Ind_N^G \psi$, we find an injective linear map
$$\Hom_G(\pi, \cInd_N^G \psi) \hookrightarrow \Hom_G(\pi, \Ind_N^G \psi) \isom \Wh_G(\pi).$$
In particular, the only irreps of $G$ which occur as \textit{subrepresentations} of a Gelfand-Graev representation are generic irreps, and moreover the uniqueness of Whittaker models implies that
$$\dim \Hom_G(\pi, \cInd_N^G \psi) \leq 1$$
for any irrep $\pi$ of $G$.

While perhaps not all generic irreps occur as subrepresentations of the Gelfand-Graev representation, we can say more about generic supercuspidal irreps.  Corollary 6.5 of \cite{CS} directly implies
\begin{prop}
Suppose that $\pi$ is a generic supercuspidal irrep of $G$.  Then $\pi$ occurs as a subrepresentation of $\cInd_N^G \psi$.
\end{prop}
Namely, the Whittaker model of a generic supercuspidal irrep of $G$ -- a priori a $G$-submodule of $\Ind_N^G \psi$ -- is in fact a $G$-submodule of $\cInd_N^G \psi$.

\subsection{Useful facts}

We will be proving that certain smooth representations of $G_2$ have no generic supercuspidal subrepresentations.  To this end, it is useful to have a few criteria that exclude such representations of $G_2$.  
\begin{prop}
Let $\pi$ be a smooth irrep of $G_2$.  Let $\alg{H}$ be a subgroup of $\alg{G}_2$, such that $\alg{H}$ is isomorphic to $\alg{SL}_3$ over an algebraic closure $\bar k$ of $k$.  If $\pi_H \neq 0$ (there exists a nonzero $H$-invariant linear functional), then $\pi$ is not generic.
\label{G2A2}
\end{prop}
\proof
Every such $\Type{A}_2$ subgroup $\alg{H}$ of $\alg{G}_2$ is conjugate over $\bar k$ (by the theory of Borel and De Siebenthal \cite{BdS}).  All such subgroups arise as stabilizers of quadratic subalgebras of $\octs$.  Lemma 4.10 of \cite{GrS1} now implies the result.
\qed

For $n \geq 4$, consider the commuting pair of split groups over $k$:
$$\alg{B}_3 \times \alg{B}_{n-4} \hookrightarrow \alg{D}_n,$$
where $\alg{B}_3 = \alg{SO}_7$, $\alg{B}_{n-4} = \alg{SO}_{2n-7}$, and $\alg{D}_n = \alg{SO}_{2n}$ are split classical groups labelled by their type.  We regard $\alg{B}_0$ as the trivial group.  Embed $\alg{G}_2$ into $\alg{B}_3$ via the action of $\alg{G}_2$ on $\octs_\circ$.
\begin{prop}
\label{NoGS}
Let $\Pi_{n}$ denote the minimal representation of $\alg{D}_n$ for $n \geq 4$.  Then, as a smooth representation of $G_2$, $\Pi_{n}$ does not have any generic supercuspidal subrepresentations.
\end{prop}
\proof
We prove this by induction on $n$.  For the base step, when $n = 4$, the proposition follows directly from Corollary 5.2 of \cite{HMS}.  

When $n > 4$, consider a maximal parabolic subgroup $\alg{P} = \alg{M} \alg{N}$ of $\alg{D}_n$ whose Levi component $\alg{M}$ satisfies 
$$\alg{G}_2 \subset \alg{B}_3 \subset \alg{D}_{n-1} \subset \alg{M} \isom \alg{GO}_{2n-2}.$$
The adjoint representation of $\alg{M}$ on $\alg{N}$ is the standard representation of $\alg{GO}_{2n-2}$; $N$ is a $(2n-2)$-dimensional vector space over $k$ with nondegenerate symmetric bilinear form.  Let $\Omega \subset N$ be the set of isotropic vectors in $N$.  By Theorem 1.1 of \cite{MS1}, there is a filtration of the minimal representation $\Pi_n$, as a representation of $P$:
$$0 \rightarrow C_c^\infty(\Omega) \rightarrow \Pi_n \rightarrow \left( \Pi_{n-1} \otimes \left| \det \right|^{\frac{1}{2n-2}} \right) \oplus \left| \det \right|^{\frac{n-2}{2n-2}} \rightarrow 0.$$
By induction, the minimal representation $\Pi_{n-1}$ of $D_{n-1}$ does not support any generic supercuspidal representations of $G_2$.  The character $\left| \det \right|^{\frac{n-2}{2n-2}}$ supports nothing but the trivial representation of $G_2$.

Finally, the representation $C_c^\infty(\Omega)$ of $G_2$ arises from the action of $G_2$ on the set of isotropic vectors in $N$.  The stabilizer of such a vector in $G_2$ is a subgroup of type $\Type{A}_2$ as discussed in the previous proposition, a subgroup isomorphic to $[Q,Q]$ for a maximal parabolic $\alg{Q} \subset \alg{G}_2$, or else all of $G_2$.  By the previous proposition, no generic supercuspidal irreps of $G_2$ have vectors fixed by an $\Type{A}_2$ subgroup.  No supercuspidal irreps have vectors fixed by $[Q,Q]$.  No nontrivial irreps have vectors fixed by all of $G_2$.  Hence no generic supercuspidal irreps of $G_2$ occur (as subrepresentations) in the restriction of $\Pi_n$ to $G_2$.
\qed

\subsection{Analysis of the correspondence}

Here we begin the analysis of the theta correspondences in $E_6$ and $E_7$, focusing on generic supercuspidal representations.  We start with the following proposition, which is primarily a consequence of results in the literature.
\begin{prop}
\label{MF7}
Let $\sigma$ be a generic supercuspidal irrep of $PGSp_{6}$.  Then $\lvec{\Theta}_{7}(\sigma)$ is a supercuspidal and multiplicity-free representation of $G_{2}$.  Every irreducible subrepresentation of $\lvec{\Theta}_7(\sigma)$ is generic.
\end{prop}
\proof
First, we prove that $\lvec{\Theta}_{7}(\sigma)$ is supercuspidal.  There are two maximal parabolic subgroups (up to conjugacy) of $\alg{G}_2$ which must be considered.
\begin{center}
\begin{tikzpicture}[scale=2]
\filldraw[fill=gray] (1,0) circle (.05cm);
\draw (1,.2cm) node[above] {Heisenberg};
\draw (1,0) node[below] {$\alpha_1$};
\draw (1.04cm,.035cm) -- (1.96cm,0.035cm);
\draw (1.05cm,0cm) -- (1.95cm, 0 cm);
\draw (1.04cm,-.035cm) -- (1.96cm, -0.035cm);
\draw (1.45cm,-.1cm) -- (1.55cm,0) -- (1.45cm,.1cm);
\draw (2,0) node[below] {$\alpha_2$};
\filldraw[fill=gray] (2,0) circle (.05cm);
\draw (2,.2cm) node[above] {Three-step};
\end{tikzpicture}
\end{center}

\begin{description}
\item[Heisenberg]
Suppose first that $\alg{Q}_2 = \alg{L}_2 \alg{U}_2$ is the Heisenberg parabolic subgroup of $\alg{G}_2$.  If $\lvec{\Theta}_{7}(\sigma)_{U_{2}} \neq 0$ then $\sigma$ occurs $(\Pi_{7})_{U_{2}}$. The structure of $(\Pi_{7})_{U_{2}}$ as an $PGSp_{6}\times L_{2}$-module has been described in \cite{MS1}, Theorem 7.6.  More precisely, one can pick a maximal parabolic subgroup $\alg{Q}_7= \alg{L}_7 \alg{U}_7$ in $\alg{E}_7$ such that $\alg{Q}_7 \cap \alg{G}_2 = \alg{Q}_2$ and $\alg{PGSp}_6 \times \alg{L}_2$ is contained in the Levi factor $\alg{L}_7$ (using the construction of Section \ref{FA}). Then we have a natural map
$$(\Pi_{7})_{U_{2}}\rightarrow (\Pi_{7})_{U_{7}}.$$
By Theorem 7.6 of \cite{MS1}, the kernel of this map does not support any supercuspidal representations of $PGSp_{6}$.  In particular, $\sigma$ must occur in $(\Pi_{7})_{U_{7}}$.  By the same result of \cite{MS1}, the representation $(\Pi_7)_{U_7}$, as a representation of $\alg{L}_7$, has constituents with wave front set supported in the closure of the minimal nilpotent orbit; the constituents are essentially a minimal representation and a trivial representation of $\alg{L}_7$.  Note that $\alg{L}_7$ is a split reductive group $\alg{CSpin}_{12}$ of type $\Type{D}_6$.

The dual pair $PGL_2 \times PGSp_6$ in a group of type $\Type{D}_6$ is addressed in Section 8 of \cite{Sav}, and no generic supercuspidal representations of $PGSp_6$ can occur.  Thus no generic supercuspidal irreps of $PGSp_6$ occur in $(\Pi_{7})_{U_{7}}$.  Therefore $\lvec{\Theta}_7(\sigma)_{U_2} = 0$.
\item[Three-step]
Now, suppose that $\alg{Q}_2 = \alg{L}_2 \alg{U}_2$ is the three-step parabolic subgroup of $\alg{G}_2$.  The structure of $(\Pi_{7})_{U_{2}}$ as an $PGSp_{6}\times L_{2}$-module has been described in \cite{GS3}, Proposition 6.8.  If $\lvec{\Theta}_7(\sigma)_{U_2} \neq 0$, then $\sigma$ occurs in $(\Pi_7)_{U_2}$.

One can pick a maximal parabolic subgroup $\alg{Q}_7= \alg{L}_7 \alg{U}_7$ in $\alg{E}_7$ such that $\alg{Q}_7 \cap \alg{G}_2 = \alg{Q}_2$ and $\alg{PGSp}_6 \times \alg{L}_2$ is contained in the Levi factor $\alg{L}_7$.  Such a parabolic subgroup is discussed and called $P_1$ in Section 4 of \cite{GS3}.  Then we have a natural map
$$(\Pi_{7})_{U_{2}}\rightarrow (\Pi_{7})_{U_{7}}.$$
The results of Proposition 6.8 of \cite{GS3} imply that the kernel does not support any supercuspidal representations of $PGSp_{6}$.  In particular, if $\sigma$ occurs in $(\Pi_7)_{U_2}$, then $\sigma$ occurs in $(\Pi_{7})_{U_{7}}$.  $\alg{L}_7$ is isogenous to $GL_2 \times PGL_6$.

By considering the Iwahori-fixed vectors, any $L_7$ constituent of the representation $(\Pi_7)_{U_7}$ is an Iwahori-spherical representation of $GL_2 \times PGL_6$ associated to the reflection or trivial representation of the Iwahori Hecke algebra of $PGL_6$.  Thus $(\Pi_7)_{U_7}$, as a representation of $PGL_6$ has all constituents appearing in degenerate principal series representations.  Such degenerate principal series restrict to degenerate principal series representations of $PGSp_6$, which are not generic.

 It follows that $\sigma$ cannot occur $(\Pi_7)_{U_7}$.  Therefore $\lvec{\Theta}_7(\sigma)_{U_2} = 0$.
\end{description}
Thus $\lvec{\Theta}_7(\sigma)$ is a supercuspidal representation of $G_2$.  It follows that $\lvec{\Theta}_7(\sigma)$ is semisimple -- a direct sum of supercuspidal irreps.

Next, we recall that $\Wh_{PGSp_6}(\Pi_7) = (\Pi_{7})_{N_{3},\psi_{3}}$ is the Gelfand-Graev module for $G_{2}$
(\cite{GS1}, Proposition 17):
$$\Wh_{PGSp_6}(\Pi_7) \isom \cInd_{N_2}^{G_2}(\psi_{2}).$$
Since $\sigma$ is a generic irreducible supercuspidal representation of $PGSp_6$, $Wh_{PGSp_6}(\sigma)$ is one-dimensional, and the embedding $\sigma \boxtimes \lvec{\Theta}_{7}(\sigma)$ into $\Pi_{7}$ gives an embedding of $\lvec{\Theta}_{7}(\sigma)$ into the Gelfand-Graev module for $G_2$.    

Since generic (and only generic) supercuspidal irreps appear as subrepresentations of the Gelfand-Graev module, and each appears with multiplicity one, we have shown that $\lvec{\Theta}_7(\sigma)$ is a multiplicity-free (though at this point, possibly empty) direct sum of generic supercuspidal irreps of $G_2$.
\qed

To summarize the previous proposition, we have found that if $\sigma$ is a generic supercuspidal irrep of $PGSp_6$, then
$$\lvec{\Theta}_7(\sigma) = \bigoplus_{i \in I} \tau_i,$$
where the right hand side denotes a (possibly empty and possibly infinite) direct sum of distinct (pairwise non-isomorphic) generic supercuspidal irreps of $G_2$.

Next, we consider $\rvec{\Theta}_6(\tau)$, when $\tau$ is a generic supercuspidal irrep of $G_2$, using the same methods as the previous proposition.
\begin{prop}
Let $\tau$ be a generic supercuspidal irrep of $G_2$.  Then $\rvec{\Theta}_{6}(\tau)$ is a supercuspidal and multiplicity-free representation of $PGL_3$.  
\end{prop}
\proof
First, we demonstrate that $\rvec{\Theta}_6(\tau)$ is supercuspidal.  There are two maximal parabolic subgroups (up to conjugacy) of $\alg{PGL}_3$ which must be considered.
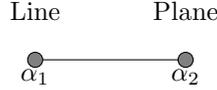
\begin{center}
\begin{tikzpicture}[scale=2]
\filldraw[fill=gray] (1,0) circle (.05cm);
\draw (1,.2cm) node[above] {Line};
\draw (1,0) node[below] {$\alpha_1$};
\draw (1.05cm,0cm) -- (1.95cm, 0 cm);
\draw (2,0) node[below] {$\alpha_2$};
\filldraw[fill=gray] (2,0) circle (.05cm);
\draw (2,.2cm) node[above] {Plane};
\end{tikzpicture}
\end{center}
\begin{description}
\item[Plane-stabilizer]
Let $\alg{Q}_2 = \alg{L}_2 \alg{U}_2$ be the maximal parabolic subgroup of $\alg{PGL}_3$ stabilizing a plane in the standard (projective) representation on $k^3$.  There exists a parabolic subgroup $\alg{Q}_6 = \alg{L}_6 \alg{U}_6$ of $\alg{E}_6$ for which $\alg{Q}_6 \cap \alg{PGL}_3 = \alg{Q}_2$ and $\alg{U}_6 \cap \alg{PGL}_3 = \alg{U}_2$.  

Theorem 4.3 of \cite{MS1} describes $(\Pi_6)_{U_2}$ as a $GL_2 \times G_2$-module; in particular, the kernel of $(\Pi_6)_{U_2} \rightarrow (\Pi_6)_{U_6}$ does not support any supercuspidal representations of $G_2$.  It follows that $\rvec{\Theta}_6(\tau)_{U_2} \boxtimes \tau$ is a $(GL_2 \times G_2)$-submodule of
$$(\Pi_6)_{U_6} \isom (\Pi_5 \otimes \left| \det \right| ) \oplus (1 \otimes \left| \det \right|^2 ),$$
where $\Pi_5$ is the minimal representation of the Levi $\alg{L}_6$ of type $\Type{D}_5$.  But no generic supercuspidal representations of $G_2$ occur in the restriction of the minimal (or trivial) representation of $Spin_{10}$ by Proposition \ref{NoGS}.  Thus $\rvec{\Theta}_6(\tau)_{U_2} = 0$.
\item[Line-stabilizer]
Let now $\alg{Q}_2' = \alg{L}_2' \alg{U}_2'$ be the maximal parabolic subgroup of $\alg{PGL}_3$ stabilizing a line in the standard representation.  Although $\alg{Q}_2'$ in not conjugate to a plane-stabilizing parabolic $\alg{Q}_2$, there exists an outer automorphism of $PGL_3$ which exchanges these two types of maximal parabolic subgroups.  Furthermore, this outer automorphism extends to an outer automorphism of $E_6$.  The uniqueness of the minimal representation of $E_6$ now demonstrates that $\rvec{\Theta}_6(\tau)_{U_2'} = 0$ as well.
\end{description}
Hence we find that $\rvec{\Theta}_6(\tau)$ is supercuspidal.  
Let $\alg{N}_2'$ denote the unipotent radical of a Borel subgroup of $\alg{PGL}_3$, and let $\psi_2'$ be a generic character of $N_2'$.  Let $\alg{N}_2$ be the unipotent radical of a Borel subgroup of $\alg{G}_2$.  By Proposition 17 of \cite{GS1}, it is known that the $G_2$-Whittaker functionals of $\Pi_6$ yield the Gelfand-Graev representation of $PGL_3$:
$$\Wh_{G_2}(\Pi_6) = (\Pi_6)_{N_2, \psi_2} \isom \cInd_{N_2'}^{PGL_3} \CC_{\psi_2'}.$$
Thus since $\tau$ is a generic supercuspidal irrep of $G_2$,  the same arguments as in Proposition \ref{MF7} imply  that $\rvec{\Theta}_6(\tau)$ is a multiplicity-free semisimple representation of $PGL_3$:  $\rvec{\Theta}_6(\tau)$ is a direct sum of pairwise non-isomorphic (automatically generic) supercuspidal irreps.
\qed

It is more complicated to analyze $\rvec{\Theta}_7(\tau)$ when $\tau$ is a generic supercuspidal irrep of $G_2$, since $\rvec{\Theta}_7(\tau)$ may or may not be supercuspidal as a representation of $PGSp_6$.  But we may consider the maximal supercuspidal (as a representation of $PGSp_6$) submodule $\rvec{\Theta}_7^\circ(\tau)$, which fits into a split short exact sequence:
$$0 \rightarrow \rvec{\Theta}_7^\circ(\tau) \rightarrow \rvec{\Theta}_7(\tau) \rightarrow \rvec{\Theta}_7^{ns}(\tau) \rightarrow 0.$$
\begin{prop}
Let $\alg{Q}_3$ denote the Siegel parabolic subgroup of $\alg{PGSp}_6$ (a maximal parabolic subgroup with abelian unipotent radical).  Then the $PGSp_6$-module $\rvec{\Theta}_7^{ns}(\tau)$ is a submodule of $\Ind_{Q_3}^{PGSp_6} \rvec{\Theta}_6(\tau)  \otimes \left| \det \right|$.  In particular, $\rvec{\Theta}_7^{ns}(\tau)$ is a (possibly empty and possibly infinite) direct sum of finite-length representations of $PGSp_6$.  If $\rvec{\Theta}_6(\tau) = 0$, then $\rvec{\Theta}_7(\tau)$ is supercuspidal.
\end{prop}
\proof
We consider the Jacquet modules of $\rvec{\Theta}_7(\tau)$, for the three (conjugacy classes) of maximal parabolic subgroups in $\alg{PGSp}_6$:
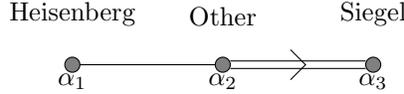
\begin{center}
\begin{tikzpicture}[scale=2]
\filldraw[fill=gray] (0,0) circle (.05cm);
\draw (0,.2cm) node[above] {Heisenberg};
\draw (0,0) node[below] {$\alpha_1$};
\draw (.05cm,0) -- (.95cm,0);
\filldraw[fill=gray] (1,0) circle (.05cm);
\draw (1,.2cm) node[above] {Other};
\draw (1,0) node[below] {$\alpha_2$};
\draw (1.05cm,.025cm) -- (1.95cm,0.025cm);
\draw (1.05cm,-.025cm) -- (1.95cm, -0.025cm);
\draw (1.45cm,-.1cm) -- (1.55cm,0) -- (1.45cm,.1cm);
\draw (2,0) node[below] {$\alpha_3$};
\filldraw[fill=gray] (2,0) circle (.05cm);
\draw (2,.2cm) node[above] {Siegel};
\end{tikzpicture}
\end{center}
The global analogues of the following computations are carried out in Case (4), of the proof of Theorem 3.1 of \cite{GRS}.
\begin{description}
\item[Heisenberg]
First, let $\alg{Q}_3 = \alg{L}_3 \alg{U}_3$ be the ``Heisenberg parabolic'', whose Levi component $\alg{L}_3$ is a split group $CSpin_5 \isom GSp_4$.   We find that $\sigma_{U_3} \boxtimes \tau$ is a quotient of $(\Pi_7)_{U_3}$ as representations of $L_3 \times G_2$.  The unipotent group $\alg{U}_3$ is $5$-dimensional, with $1$-dimensional center $\alg{Z}_3$; there exists a parabolic subgroup $\alg{Q}_7 = \alg{L}_7 \alg{U}_7$ of $\alg{E}_7$ such that $\alg{L}_7$ is isomorphic to $\alg{CSpin}_{12}$, and $\alg{U}_7$ is a Heisenberg group of dimension $33$ (with one-dimensional center $\alg{Z}_7$).  Furthermore, one may choose this parabolic subgroup in such a way that $\alg{Q}_7 \cap \alg{PGSp}_6 = \alg{Q}_3$, $\alg{U}_7 \cap \alg{PGSp}_6 = \alg{U}_3$, and $\alg{Z}_7 = \alg{Z}_3$.  Furthermore, this gives an embedding, $\alg{L}_3 \times \alg{G}_2 \hookrightarrow \alg{L}_7 = \alg{CSpin}_{12}$. 

Now, $\rvec{\Theta}_7(\tau)_{U_3} \boxtimes \tau$ is a subrepresentation of $(\Pi_7)_{U_3}$.  To study $(\Pi_7)_{U_3}$, we examine a commutative diagram with exact rows and columns:
$$\xymatrix{
& 0 \ar[d] & 0 \ar[d] & & \\
0 \ar[r] & W \ar[r] \ar[d] & (\Pi_7)_{Z_3} \ar[r] \ar[d] & (\Pi_7)_{U_3} \ar[d] \ar[r] & 0 \\
0 \ar[r] & C_c^\infty(\Omega) \ar[r] & (\Pi_7)_{Z_7} \ar[r] \ar[d] & (\Pi_7)_{U_7} \ar[r] \ar[d] & 0 \\
& & 0 & 0 &.
}$$
Here $\Omega$ denotes the smallest nontrivial $L_7$-orbit in the 32-dimensional vector space $U_7 / Z_7$; this can be identified with the 15-dimensional quotient $CSpin_{12} / Q_6$, where $Q_6$ is a minuscule maximal parabolic subgroup (with Levi subgroup of type $A_5$) of $CSpin_{12}$.  Geometrically, $\Omega$ can be viewed as a Grassmannian of isotropic $6$-spaces in the 12-dimensional standard representation $V$ of $Spin_{12}$.

From Theorem 6.1 of \cite{MS1}, the kernel of $(\Pi_7)_{Z_7} \rightarrow (\Pi_7)_{U_7}$ can be identified, as a $Q_7$-module, with $C_c^\infty(\Omega)$.  We are led to consider the action and orbits of $G_2 \times Spin_5$ on $\Omega$.  It helps to study the action of $Spin_7 \times Spin_5$ on $\Omega$.  Here the embedding of $Spin_7 \times Spin_5$ in $Spin_{12}$ corresponds to a decomposition $V = V_7 \oplus V_5$ of the standard representation of $Spin_{12}$. 

Such actions are studied by Kudla, in Proposition 3.4 of \cite{Kud}.  If $\omega \in \Omega$ corresponds to an isotropic $6$-space $\Lambda_\omega$, then the projection of $\Lambda_\omega$ onto $V_7$ is at least one-dimensional.  It follows that $\omega$ is stabilized by some maximal parabolic subgroup $Q$ of $Spin_7$.  

It follows that the stabilizer $S_\omega$ of $\omega$ in $G_2$ contains a maximal parabolic subgroup of $G_2$, or $S_\omega$ contains a subgroup of type $\Type{A}_2$ (by the arguments of Proposition \ref{GroupTheory}).  If $S_\omega$ contains a maximal parabolic subgroup of $G_2$, then $C_c^\infty(G_2 / S_\omega)$ does not support any supercuspidal representations of $G_2$.  If $S_\omega$ contains a subgroup of type $\Type{A}_2$, then $C_c^\infty(G_2 / S_\omega)$ does not support any generic supercuspidal representations of $G_2$ by Proposition \ref{G2A2}.  Thus $C_c^\infty(\Omega)$ does not support any generic supercuspidal representations of $G_2$.

Since $Ker \left( (\Pi_7)_{U_3} \rightarrow (\Pi_7)_{U_7} \right)$ is a quotient of $C_c^\infty(\Omega)$ (by the snake lemma), we find that $\rvec{\Theta}_7(\tau)_{U_3} \boxtimes \tau$ is a subrepresentation of $(\Pi_7)_{U_7}$.  But this implies that $\tau$ occurs in the restriction of the minimal representation of $Spin_{12}$ or else $\rvec{\Theta}_7(\tau)_{U_3} = 0$.  By Proposition \ref{NoGS}, no generic supercuspidal representations of $G_2$ occur in this restriction.  It follows that $\rvec{\Theta}_7(\tau)_{U_3} = 0$.
\item[Other]
Next, let $\alg{Q}_3 = \alg{L}_3 \alg{U}_3$ denote the ``Other parabolic'', with $\alg{L}_3$ isogenous to $GL_2 \times SL_2$.  The unipotent radical $\alg{U}_3$ has three-dimensional center $\alg{Z}_3$, and four-dimensional quotient $\alg{U}_3 / \alg{Z}_3$.  $Z_3$ can be identified with the space $M_\circ$ of two-by-two matrices with trace zero, and $U_3 / Z_3$ can be identified with the space $M$ of all two-by-two matrices.

We find that $\rvec{\Theta}_7(\tau)_{U_3} \boxtimes \tau$ is a subrepresentation of $(\Pi_7)_{U_3}$, as representations of $L_3 \times G_2$.

There exists a parabolic subgroup $\alg{Q}_7 = \alg{L}_7 \alg{U}_7$ such that $\alg{L}_7$ is isogenous to $CSpin_{10} \times SL_2$, $\alg{U}_7$ is a two-step unipotent group with $10$-dimensional center $\alg{Z}_7$, and $\alg{Q}_7 \cap \alg{PGSp}_6 = \alg{Q}_3$ and $\alg{L}_7 \cap \alg{PGSp}_6 = \alg{L}_3$.  $Z_7$ can be identified with the space $M_\circ \oplus \octs_\circ$ of pairs $(m,\omega)$, and $U_7 / Z_7$ can be identified with the (32-dimensional) space $M \otimes \octs$.  This arises from the construction of Section \ref{TC}.  This parabolic arises in a similar computation in \cite{GS3}, and our $\alg{Q}_7$ corresponds to the parabolic called $P$ and associated to the vertex $\alpha_4$ in \cite{GS3}.

There are natural short exact sequences which we describe and analyze below:
$$0 \rightarrow C_c^\infty(\Omega, \mc{S}) \rightarrow (\Pi_7) \rightarrow (\Pi_7)_{Z_7} \rightarrow 0,$$
$$0 \rightarrow C_c^\infty(\Omega') \rightarrow (\Pi_7)_{Z_7} \rightarrow (\Pi_7)_{U_7} \rightarrow 0.$$
Here $\Omega$ is the set of nontrivial characters $\omega$ of $Z_7$ for which $(\Pi_7)_{Z_7, \omega} \neq 0$. 
On $\Omega$, $\mc{S}$ is a $\alg{Q}_7$-equivariant sheaf whose fibre over $\omega \in \Omega$ is an irreducible representation of $U_7$ with central character corresponding to $\omega$.  One can compare this to Section 6 of \cite{GS3}.

Similarly, $\Omega'$ is the set of nontrivial characters of $U_7$ for which $(\Pi_7)_{U_7, \omega} \neq 0$.  Identifying characters of $U_7$ with $\octs \otimes M$, a minuscule representation of $L_7$, $\Omega'$ can be identified with the quotient $L_7 / P_6$ where $\alg{P}_6$ is a minuscule parabolic subgroup of $\alg{L}_7$.

Taking $U_3$ co-invariants in each of the short exact sequences, we are led to consider $C_c^\infty(\Omega, \mc{S})_{U_3}$ and $C_c^\infty(\Omega')_{U_3}$.  In the first case, we find that $C_c^\infty(\Omega, \mc{S})_{U_3}$ is a quotient of $C_c^\infty(\Omega, \mc{S})_{Z_3}$.  We compute
$$C_c^\infty(\Omega, \mc{S})_{Z_3} \isom C_c^\infty(\Omega^{\perp Z_3}, \mc{S}),$$
where $\Omega^{\perp Z_3}$ can be identified:
\begin{eqnarray*}
\Omega^{\perp Z_3} & = &\{ (m, \omega) \in M_\circ \oplus \octs_\circ : N(\omega) - N(m) = 0 \mbox{ and } m = 0 \} \\
& = & \{ \omega \in \octs : \omega^2 = 0 \}.
\end{eqnarray*}
It follows that
$$C_c^\infty(\Omega, \mc{S})_{Z_3} \isom \Ind_{P_\omega}^{G_2} \mc{S}_\omega,$$
where $\mc{S}_\omega$ is the fibre of $\mc{S}$ over $\omega \in \octs$, which satisfies $\omega^2 = 0$, and $P_\omega$ is the maximal parabolic subgroup of $G_2$ stabilizing $\omega$.  The representation $\mc{S}_\omega$ of $P_\omega$ factors through the Levi quotient $L_\omega \isom GL_2$ of $P_\omega$.  It follows that $C_c^\infty(\Omega, \mc{S})_{Z_3}$ and hence $C_c^\infty(\Omega, \mc{S})_{U_3}$ does not support any supercuspidal representations of $G_2$.

Next we are led to consider $C_c^\infty(\Omega')$.  Every point of $\Omega'$ corresponds to an isotropic 5-plane $\Lambda$ in $\octs_\circ \oplus M_\circ$ (the standard representation of $CSpin_{10} \subset L_7$), since $CSpin_{10}$ acts via the Spin representation on $U_7 / Z_7$.  The projection of $\Lambda$ onto $\octs_\circ$ is at least 2-dimensional; hence $\Lambda$ is stabilized by a maximal parabolic subgroup of $Spin_7 \subset Spin_{10}$.  Hence $\Lambda$ is stabilized by a maximal parabolic subgroup of $G_2$, or by a subgroup of type $\Type{A}_2$ in $G_2$.  It follows that $C_c^\infty(\Omega')$ does not support any generic supercuspidal representations of $G_2$ using Propositions \ref{G2A2} and \ref{NoGS}.  

By the snake lemma argument as before, we find that  $\rvec{\Theta}_7(\tau)_{U_3} \boxtimes \tau$ occurs as a subrepresentation of $(\Pi_7)_{U_7}$.  The representation $(\Pi_7)_{U_7}$ of $L_7$ has wave front set supported in the minimal orbit.  If $\rvec{\Theta}_7(\tau)_{U_3}$ were nontrivial, then $\tau$ would occur in a theta correspondence $G_2 \times (Spin(3) \times SL_2) \subset CSpin_{10} \times SL_2$.  But no generic supercuspidal representations of $G_2$ occur in such a corresponence, by Proposition \ref{NoGS}.  Hence $\rvec{\Theta}_7(\tau)_{U_3} = 0$.
\item[Siegel]
Finally, let $\alg{Q}_3 = \alg{L}_3 \alg{U}_3$ denote the ``Siegel parabolic'', with $\alg{L}_3 \isom \alg{GL}_3$.   We find that $\rvec{\Theta}_7(\tau)_{U_3} \boxtimes \tau$ is a subrepresentation of $(\Pi_7)_{U_3}$, as representations of $GL_3 \times G_2$.  Let $\alg{Q}_7$ denote a maximal parabolic subgroup of $\alg{E}_7$ whose Levi component has derived subgroup $\alg{E}_6$, such that $\alg{Q}_7 \cap \alg{PGSp}_6 = \alg{Q}_3$.  These embeddings and parabolics arise from the construction of Section \ref{GJ}.  By Theorem 5.3 of \cite{MS1}, the kernel of $(\Pi_7)_{U_3} \twoheadrightarrow (\Pi_7)_{U_7}$ does not support any supercuspidal representations of $G_2$.  It follows that $\rvec{\Theta}_7(\tau)_{U_3} \boxtimes \tau$ is a subrepresentation of $(\Pi_7)_{U_7}$.  

By Theorem 5.3 of \cite{MS1} again, there is a $G_2 \times GL_3$-module isomorphism
$$(\Pi_7)_{U_7} \isom \left( \Pi_6 \otimes \left| \det \right| \right)  \oplus \left( 1 \otimes \left| \det \right|^2 \right).$$
Taking $(G_2,\tau)$-isotypic components, we find an isomorphism of $GL_3$-modules:
$$\rvec{\Theta}_7(\tau)_{U_3} \isom \rvec{\Theta}_6(\tau) \otimes \left| \det \right|.$$
\end{description}
For the rest of the proof, let $\alg{Q}_3 = \alg{L}_3 \alg{U}_3$ denote the Siegel parabolic subgroup of $\alg{PGSp}_6$.  The previous computations and Frobenius reciprocity yield a morphism of $PGSp_6$-modules:
$$\rvec{\Theta}_7(\tau) \rightarrow \Ind_{Q_3}^{PGSp_6} \rvec{\Theta}_6(\tau) \otimes \left| \det \right|.$$
Moreover, the kernel of this morphism is a submodule of $\rvec{\Theta}_7(\tau)$ whose $U_3$-coinvariants vanish.  But since all other (with respect to the Heisenberg parabolic and ``Other'' parabolic) Jacquet modules of $\rvec{\Theta}_7(\tau)$ vanish, the kernel of this morphism is a supercuspidal $PGSp_6$-submodule of $\rvec{\Theta}_7(\tau)$.  Conversely, every supercuspidal $PGSp_6$-submodule of $\rvec{\Theta}_7(\tau)$ is contained in the kernel of the morphism, since supercuspidals do not occur as subrepresentations of parabolically induced representations.

It follows that there is an injective morphism of $PGSp_6$-modules:
$$\rvec{\Theta}_7^{ns}(\tau) \hookrightarrow \Ind_{Q_3}^{PGSp_6} \rvec{\Theta}_6(\tau) \otimes \left| \det \right|.$$
The previous proposition implies that there exists a set of pairwise nonisomorphic supercuspidal irreps $\{ \rho_i \}_{i \in I}$ of $PGL_3$, such that 
$$\rvec{\Theta}_6(\tau) \isom \bigoplus_{i \in I} \rho_i.$$
It follows that there is an injective morphism of $PGSp_6$-modules:
$$\rvec{\Theta}_7^{ns}(\tau) \hookrightarrow \bigoplus_{i \in I} \Ind_{Q_3}^{PGSp_6} \rho \otimes \left| \det \right|.$$
Although there may be an infinite number of summands on the right side above, only finitely many lie in any given Bernstein component for $PGSp_6$.  We find that $\rvec{\Theta}_7^{ns}(\tau)$ is a (possibly infinite and possibly empty) direct sum of finite-length representations of $PGSp_6$.  Moreover, if $\rvec{\Theta}_6(\tau) = 0$, then $\rvec{\Theta}_7^{ns}(\tau)$ vanishes, and so $\rvec{\Theta}_7(\tau)$ is supercuspidal.
\qed

To synthesize the previous propositions, we find that for any generic supercuspidal irrep $\tau$ of $G_2$, there is a set $\{ \sigma_j \}_{j \in J}$ of supercuspidal irreps of $PGSp_6$, a set $\{ \rho_i \}_{i \in I}$ of supercuspidal irreps of $PGL_3$, and a set of finite-length $PGSp_6$ modules $\{ \pi_i \}_{i \in I}$ satisfying:
$$\rvec{\Theta}_6(\tau) \isom \bigoplus_{i \in I} \rho_i,$$
$$\rvec{\Theta}_7(\tau) \isom \bigoplus_{j \in J} \sigma_j \oplus \bigoplus_{i \in I} \pi_i \mbox{ and } \pi_i \subset \Ind_{Q_3}^{PGSp_6} \rho_i \mbox{ for all } i \in I.$$
The above decomposition refines the decomposition of $\rvec{\Theta}_7(\tau)$ into supercuspidal and non-supercuspidal parts:
$$\rvec{\Theta}_7^\circ(\tau) \isom \bigoplus_{j \in J} \sigma_j \mbox{ and } \rvec{\Theta}_7^{ns}(\tau) \isom \bigoplus_{i \in I} \pi_i.$$

\begin{prop}
Let $\tau$ be a generic supercuspidal irrep of $G_{2}$.  Let $\{ \sigma_j \}_{j \in J}$, $\{ \rho_i \}_{i \in I}$ and $\{\pi_i \}_{i \in I}$ be the representations of $PGSp_6$, $PGL_3$, and $PGSp_6$ in the above decomposition.  Then $\rvec{\Theta}_{7}(\tau)$ is nontrivial (so $I \sqcup J \neq \emptyset$).  Moreover, exactly one of the following statements holds: 
\begin{enumerate}
\item
There exists exactly one $j \in J$ such that $\sigma_j$ is generic.  There does not exist $i \in I$ such that $\pi_i$ is generic.
\item
There exists exactly one $i \in I$ such that $\pi_i$ is generic.  There does not exist $j \in j$ such that $\pi_j$ is generic.
\end{enumerate}
\label{G2toC3}
\end{prop}
\proof
For $\tau$ a generic supercuspidal irrep of $G_2$, $\tau$ occurs with multiplicity one in the Gelfand-Graev module:
$$\dim \left( \Hom_{G_2} (\tau, \cInd_{N_2}^{G_2} \psi_2) \right) = 1.$$
But using Proposition 17 of \cite{GS1} again,  
$$\Wh_{PGSp_6}(\Pi_7) = (\Pi_{7})_{N_{3},\psi_{3}} \cong \cInd_{N_2}^{G_2}(\psi_{2}).$$
Thus we find that
$$\dim \left( \Hom_{G_2}(\tau, (\Pi_7)_{N_3, \psi_3}) \right) = \dim \left( \Hom_{G_2}(\tau, \Pi_7) \right)_{N_3, \psi_3} = 1.$$
Thus $\Wh_{PGSp_6}(\rvec{\Theta}_7(\tau))$ is one-dimensional.  In particular, $\rvec{\Theta}_7(\tau)$ is nontrivial.

Now, we apply the decomposition:
$$\rvec{\Theta}_7(\tau) \isom \bigoplus_{j \in J} \sigma_j \oplus \bigoplus_{i \in I} \pi_i \mbox{ and } \pi_i \subset \Ind_{Q_3}^{PGSp_6} \rho_i \mbox{ for all } i \in I.$$
Taking Whittaker functionals, we find
$$\Wh_{PGSp_6}(\rvec{\Theta}_7(\tau)) \isom \bigoplus_{j \in J} \Wh_{PGSp_6}(\sigma_j) \oplus \bigoplus_{i \in I} \Wh_{PGSp_6}(\pi_i).$$
Since the left side is one-dimensional, precisely one summand on the right side is one-dimensional and all other summands on the right side vanish.  The result follows immediately.
\qed

When the residue characteristic $p$ is odd, the representations $\Ind_{Q_3}^{PGSp_6} (\rho \otimes \left| \det \right|)$ are irreducible and generic, whenever $\rho$ is a supercuspidal irrep of $PGL_3$.  This significantly simplifies the analysis of the theta correspondence, in the following way:
\begin{prop}
Suppose that $p \neq 2$.  Let $\tau$ be an generic supercuspidal irrep of $G_2$.  Then if $\rvec{\Theta}_6(\tau) \neq 0$ then $\rvec{\Theta}_6(\tau)$ has a unique irreducible subrepresentation, up to contragredience:  $\rvec{\Theta}_6(\tau) = \rho \oplus \tilde \rho$ for some supercuspidal irrep $\rho$ of $PGL_3$.
\end{prop}
\proof
If $\rho$ is a irreducible subrepresentation of $\rvec{\Theta}_6(\tau)$ (and hence $\rho$ is generic and supercuspidal), then $\rho \boxtimes \tau$ occurs as a quotient (by the injectivity and projectivity of supercuspidals) of the minimal representation $\Pi_6$ of the adjoint group $E_6$.  But we have seen that if $\alg{Q}_3 = \alg{L}_3 \alg{U}_3$ is the Siegel parabolic subgroup of $PGSp_6$, then there is a surjective map of $GL_3$-modules:
$$\rvec{\Theta}_7(\tau)_{U_3} \twoheadrightarrow (\rho \otimes \left| \det \right|).$$
By Frobenius reciprocity, we find a nontrivial map of $PGSp_6$-modules:
$$\rvec{\Theta}_7(\tau) \rightarrow  \Ind_{Q_3}^{PGSp_6} \left( \rho \otimes \left| \det \right| \right).$$
Let $\pi$ denote this induced representation, $\pi =   \Ind_{Q_3}^{PGSp_6}\left( \rho \otimes \left| \det \right| \right)$.  Since $p \neq 2$, the representation $\rho$ is {\em not} self-contragredient, and so $\pi$ is an irreducible generic representation of $PGSp_6$.

Thus $\pi$ must be the {\em unique} generic summand of $\rvec{\Theta}_7(\tau)$ in the decomposition
$$\rvec{\Theta}_7(\tau) \isom \bigoplus_{j \in J} \sigma_j \oplus \bigoplus_{i \in I} \pi_i.$$
By the geometric lemma and Frobenius reciprocity (using the fact that $\rho$ is supercuspidal), the only representations of $GL_3$ which parabolically induce to give this representation $\pi$ of $PGSp_6$ are $\rho$ and its contragredient $\tilde \rho$.  Hence $\rvec{\Theta}_6(\tau)$ contains a unique irreducible subrepresentation up to contragredience, and this irrep and its contragredient are supercuspidal.  This demonstrates that 
$$\rho \subset \rvec{\Theta}_6(\tau) \subset \rho \oplus \tilde \rho.$$  

Lastly note that the map $\rvec{\Theta}_7(\tau) \rightarrow \pi$ is surjective, from which it follows that the map
$$\rvec{\Theta}_7(\tau)_{U_3} \rightarrow \pi_{U_3}$$
is also surjective.  But both $\rho \otimes \left| \det \right|$ and $\tilde \rho \otimes \left| \det \right|$ occur in $\pi_{U_3}$.  Since $\rvec{\Theta}_7(\tau)_{U_3} \isom \rvec{\Theta}_6(\tau) \otimes \left| \det \right|$, we find that both $\rho$ and $\tilde \rho$ occur in $\rvec{\Theta}_6(\tau)$.
\qed

Using the previous propositions, we find that (regardless of residue characteristic):
\begin{thm}
Suppose that $\tau$ is a generic supercuspidal irrep of $G_2$.  Then either there exists a unique generic supercuspidal irreducible subrepresentation $\sigma$ of $\rvec{\Theta}_7(\tau)$ or else there exists a unique -- up to contragredience -- generic supercuspidal irreducible subrepresentation $\rho$ of $\rvec{\Theta}_6(\tau)$ for which the generic summand of $\Ind_{Q_3}^{PGSp_6}(\rho \otimes \left| \det \right|)$ occurs in $\rvec{\Theta}_7(\tau)$.

In this way, the theta correspondences yield a map:
$$\Delta \colon \Irr_g^\circ(G_2) \rightarrow \Irr_{g}^\circ(PGSp_6) \sqcup \frac{\Irr_g^\circ(PGL_3)}{\Contra},$$
$$\tau \mapsto \sigma \mbox{ or } \{ \rho, \tilde \rho \}.$$
\end{thm}

When $p \neq 2$, we find that the dichotomy map is given somewhat simply by
$$\Delta(\tau) = \begin{cases}
          \hfill \sigma &\text{if} \quad \rvec{\Theta}_6(\tau) = 0 \\
           \{\rho, \tilde \rho \} &\text{if} \quad \rvec{\Theta}_6(\tau) \neq 0 \\
         \end{cases} $$
However when $p = 2$, it is possible a priori that a self-contragredient supercuspidal irrep $\rho$ occurs as a summand of $\Theta_6(\tau)$, the {\em non-generic}  summand $\pi$ of $\Ind_{Q_3}^{PGSp_6} \rho \otimes \left| \det \right|$ occurs as a summand of $\rvec{\Theta}_7(\tau)$, and still a generic supercuspidal representation of $PGSp_6$ occurs as a summand of $\rvec{\Theta}_7(\tau)$.  We cannot yet exclude such a strange possibility.
 
\section{Shalika Functionals}

\subsection{The Shalika subgroup}
It is convenient hereafter to view $\alg{GSp}_6$ in the traditional way, as a group of symplectic similitudes.  We let $\alg{M_2}$ denote the abelian unipotent algebraic group of two by two matrices (under addition); if $g$ is a matrix, we write $g^\trans$ for its transpose.  

Let $I = \Matrix{1}{0}{0}{1}$, $J=\Matrix{0}{-1}{1}{0}$, and
$$J_3 = \TMatrix{0}{0}{J}{0}{J}{0}{J}{0}{0}.$$
Let $\alg{GSp}_6$ be the algebraic group of symplectic similitudes:
$$\alg{GSp}_6 = \{ g \in \alg{GL}_6 : g J_3  g^\trans = \Sim(g) \cdot J_3 \mbox{ for some } \Sim(g) \in \alg{GL}_1 \}. $$
The resulting character $\Sim \colon \alg{GSp}_6 \rightarrow \alg{GL}_1$ is called the similitude character.

Let $\alg{Q}_3 = \alg{L}_3 \alg{U}_3$ be the maximal  parabolic subgroup of $\alg{GSp}_{6}$, with Levi component
$$\alg{L}_3 = \lset \TMatrix{g}{0}{0}{0}{h}{0}{0}{0}{\det(g^{-1}h) \cdot g} : g,h \in \alg{GL}_2 \rset,$$
and unipotent radical
$$\alg{U}_3 = \lset \TMatrix{I}{X}{Z}{0}{I}{Y}{0}{0}{I} : X,Y,Z \in \alg{M}_2, XJ + J Y^\trans = 0, ZJ + J Z^\trans = -X J X^\trans \rset.$$
The center of $\alg{U}_3$ is three-dimensional,
$$\alg{Z}_3 = \lset \TMatrix{I}{0}{Z}{0}{I}{0}{0}{0}{I} : ZJ + J Z^\trans = 0 \rset.$$
There is an isomorphism of unipotent groups $\alg{U}_3 / \alg{Z}_3 \rightarrow \alg{M}_2$, given by
$$\TMatrix{I}{X}{Z}{0}{I}{Y}{0}{0}{I} \mapsto X.$$
There is also an isomorphism of reductive groups $\alg{L}_3 \rightarrow \alg{GL}_2 \times \alg{GL}_2$ given by
$$\TMatrix{g}{0}{0}{0}{h}{0}{0}{0}{\det(g^{-1} h) \cdot g} \mapsto (g,h).$$
With these identifications, the conjugation action of $\alg{L}_3$ on $\alg{U}_3 / \alg{Z}_3$ is given by
$$(g,h) \cdot X = g X h^{-1}.$$
Let $\Delta \colon \alg{GL}_2 \rightarrow \alg{GL}_2 \times \alg{GL}_2 \isom \alg{L}_3$ denote the diagonal embedding (there should be no risk of confusing this $\Delta$ with the dichotomy map in other sections).  If $g \in \alg{GL}_2$, then $\Delta(g)$ is identified with an element of $\alg{L}_3 \subset \alg{GSp}_6$:
$$\Delta(g) = \TMatrix{g}{0}{0}{0}{g}{0}{0}{0}{g}.$$  
Then, we write $\alg{S}$ for the ``Shalika subgroup'':
$$\alg{S} = \Delta(\alg{GL}_2) \ltimes \alg{U}_3 \subset \alg{Q}_3.$$

Observe also that the Shalika subgroup has another interpretation:  If $A$ is any $k$-algebra, consider the degenerate cubic $A$-algebra $A[\epsilon] / \langle \epsilon^3 \rangle$.  Then there is a natural inclusion (of codimension 1):
$$\alg{S}(A) \subset \alg{GL}_2\left( A[\epsilon] / \langle \epsilon^3 \rangle \right).$$

Define a character $\psi_3$ of $U_3$ by $\psi_{3}(u)=\psi_k(-\Tr(X))$ (for a matrix $u \in U_3$ projecting to $X \in M_2 \isom U_3 / Z_3$).  $\Delta(GL_2)$ is precisely the centralizer of the character $\psi_3$ in $L$; hence the character $\psi_{3}$ can be extended uniquely to a character $\psi_S$ of $S$ such that $\psi_S(\Delta(g))=1$ for all $g \in GL_2$. 

When $\sigma$ is a smooth representation of $GSp_6$, we define the space of \textbf{Shalika functionals}  by
$$\Sh(\sigma) = \Hom_S(\sigma, \psi_S).$$
Note that, if $\sigma$ has a nonzero Shalika functional, then the central character of $\sigma$ is trivial.  The main goal of this section is to demonstrate that for supercuspidal irreps $\sigma$ of $PGSp_6$, $\dim(\Sh(\sigma)) \leq 1$ -- the ``uniqueness'' of Shalika functionals.

Our methods are similar to many other papers; we mention the work of Jacquet and Rallis \cite{J-R}, who prove uniqueness of Shalika models for $GL_{2n}$.  The $k$-points of their ``Shalika subgroup'' can be identified with $GL_n\left( k[\epsilon] / \langle \epsilon^2 \rangle \right)$.  While their Shalika functionals are related to a degenerate quadratic algebra, ours are related to a degenerate cubic algebra.

\subsection{Double cosets}

If $g \in GSp_6$, then its transpose $g^\trans$ is also an element of $GSp_6$, and the transpose is an involution (anti-automorphism of order two) of $\alg{GSp}_6$.  If $\alg{H} \subset \alg{G}$ is an algebraic subgroup, we write $\alg{H}^\trans$ for its transpose.

We will require an explicit description of the double cosets $\alg{Q}_3^\trans \backslash \alg{GSp}_6 / \alg{Q}_3$ as well as $\alg{S}^\trans \backslash \alg{GSp}_6 / \alg{S}$.  As $\alg{Q}_3$ is a maximal parabolic subgroup of $\alg{GSp}_6$, the first is a routine computation; it suffices to find representatives for double cosets in the Weyl group of type $\Type{C}_3$, modulo the parabolic subgroup of type $\Type{A}_1 \times \Type{A}_1$.  For this, we define elements of $\alg{GSp}_6$ corresponding to simple root reflections $a,b,c$ (though we refrain from identifying a maximal torus, Borel subgroup, et cetera):
\begin{center}
\begin{tikzpicture}[scale=2]
\filldraw[fill=gray] (0,0) circle (.05cm);
\draw (0,0) node[below] {$a$};
\draw (.05cm,0) -- (.95cm,0);
\filldraw[fill=gray] (1,0) circle (.05cm);
\draw (1,0) node[below] {$b$};
\draw (1.05cm,.025cm) -- (1.95cm,0.025cm);
\draw (1.05cm,-.025cm) -- (1.95cm, -0.025cm);
\draw (1.45cm,-.1cm) -- (1.55cm,0) -- (1.45cm,.1cm);
\draw (2,0) node[below] {$c$};
\filldraw[fill=gray] (2,0) circle (.05cm);
\end{tikzpicture}
\end{center}
$$a = \TMatrix{J}{0}{0}{0}{I}{0}{0}{0}{J}, c = \TMatrix{I}{0}{0}{0}{J}{0}{0}{0}{I},$$
$$b = \left(\begin{array}{ccc|ccc}1 & 0 & 0 & 0 & 0 & 0 \\0 & 0 & 1 & 0 & 0 & 0 \\0 & -1 & 0 & 0 & 0 & 0 \\ \hline 0 & 0 & 0 & 0 & 1 & 0 \\0 & 0 & 0 & -1 & 0 & 0 \\0 & 0 & 0 & 0 & 0 & 1 \end{array}\right) .$$
\begin{prop}
The algebraic variety $\alg{GSp}_6$ can be decomposed as a finite disjoint union
$$\alg{GSp}_6 = \bigsqcup_{\sigma \in \Sigma}  \alg{Q}_3^\trans \sigma \alg{Q}_3,$$
where 
$$\Sigma = \{ 1, b, bcb, bacb, bcabacb \}.$$
\end{prop}
\proof
The nontrivial shortest representatives for double cosets in the Weyl group are given by words in $a$, $b$, $c$, which begin and end with $b$.  These can be found by direct computation, using the relations in the Coxeter group.
\qed

Define an embedding $\eta$ of $\alg{GL}_2$ into $\alg{L}_3$ by:
$$\eta(g) = \TMatrix{I}{0}{0}{0}{g}{0}{0}{0}{\det(g) \cdot I},$$
Then it can be easily verified that:
$$ \alg{L}_3  = \Delta(\alg{GL}_2) \eta(\alg{GL}_2) =  \eta(\alg{GL}_2)  \Delta(\alg{GL}_2).$$
The previous proposition now implies
\begin{cor}
The algebraic variety $\alg{GSp}_6$ can be decomposed as a disjoint union
$$\alg{GSp}_6 = \bigsqcup_{\sigma \in \Sigma} \alg{S}^\trans \eta(\alg{GL}_2) \sigma \eta(\alg{GL}_2) \alg{S}.$$
\end{cor}

Let $\alg{R} =  \alg{Q}_3^\trans \sigma \alg{Q}_3$ be a double coset in $\alg{GSp}_6$.  Then we find that
$$\alg{R}^\trans = \alg{Q}_3^\trans  (\sigma^\trans) \alg{Q_3} = \alg{R}.$$
If $s \in S$, then we define a character $\psi_S^\trans$ of $S^\trans$ by
$$\psi_S^\trans(s) = \psi_S(s^\trans).$$

\subsection{Distributions}

If $X$ is a subset of $GSp_6$ and $X = S^\trans X S$, then there is a natural action $\ell \times \rho$ of $S^\trans \times S$ on $C_c^\infty(X)$, given by:
$$[\ell(s) \rho(t) f](x) = [\rho(t) \ell(s) f](x) = f(s^{-1} x t),$$
for all  $s \in  S^\trans$, $t \in S$, $x \in X$, $f \in C_c^\infty(X)$.   If $T$ is a distribution on $X$, i.e. $T$ is a linear functional on $C_c^\infty(X)$, then we say that $T$ is $(S, \psi, \trans)$-invariant if for all $s \in S^\trans$, $t \in S$, $f \in C_c^\infty(X)$,
$$T((\ell(s) \rho(t) f) = \psi_S^\trans(s) \psi_S(t^{-1}) T(f).$$

We frequently apply the following restrictions on the support of such distributions:  
\begin{description}
\item[R1]
If $s \in S$, $g \in G$, $g s g^{-1} \in S^\trans$, and $\psi_S(s) \neq \psi_S^\trans(g s g^{-1})$, then the coset $S^\trans g S$ does not support any $(S, \psi, \trans)$-invariant distributions.
\item[R2]
If $s \in S$, $g \in G$, $g^{-1} s^\trans g \in S$, and $\psi_S^\trans(s^\trans) \neq \psi_S(g^{-1} s^\trans g)$, then the coset $S^\trans g S$ does not support any $(S, \psi, \trans)$-invariant distributions.
\end{description}
These restrictions follow directly from Bernstein's localization principle; this method is used often in the study of Shalika and Whittaker models, and we point to the recent work of Jiang, Nien, and Qin \cite{JNQ} for an example similar in spirit.

We will also apply the following criterion to prove transpose-invarinace of distributions:
\begin{description}
\item[TI]
If $g \in G$, and there exist $s_1, s_2 \in S \cap S^\trans = \Delta(GL_2)$ such that $s_1 g s_2 = g^\trans$, then any $(S, \psi, \trans)$-invariant distribution on $S^\trans g S$ is also transpose-invariant.
\end{description}

Following the methods of Gelfand-Kazhdan \cite{GeK}, we prove the following
\begin{thm}
\label{DistThm}
Let $R = Q_3^\trans \sigma Q_3$ be a double coset in $GSp_6$.  Suppose that $T$ is a $(S, \psi, \trans)$-invariant distribution on $R$.  Then $T$ is transpose-invariant.  
\end{thm}  
\proof
We prove this theorem, by analyzing the five cosets $Q_3^\trans \sigma Q_3$ individually.  We whittle down the support of such a distribution $T$ using the restrictions (R1) and (R2), and prove transpose-invariance using criterion (TI).
\begin{description}
\item[$\sigma = 1$]
For $\sigma = 1$, we are led to consider distributions $T$ on $R =  Q_3^\trans \cdot Q_3 = U_3^\trans \cdot L_3 \cdot U_3$.  Since $L_3$ normalizes both $U_3$ and $U_3^\trans$ the $(S, \psi, \trans)$-invariant distributions $T$ on $R$ are in natural correspondence with distributions on $L_3$ which are $\Delta(GL_2)$ bi-invariant.  

Thus we are led to consider the orbits for the action $\alpha$ of $\Delta(GL_2) \times \Delta(GL_2) \isom GL_2 \times GL_2$ on $L_3 \isom GL_2 \times GL_2$, given by:
$$\alpha(g,h)(x,y) = (gxh, gyh).$$
Clearly every element $(x,y) \in L_3 \isom GL_2 \times GL_2$ is in the same orbit as $(1, x^{-1} y)$.  Furthermore we find that for all $g \in GL_2$, $(1, x^{-1} y)$ is in the same orbit as $(1, g x^{-1} y g^{-1})$.  Finishing this analysis, we find that the orbits of $GL_2 \times GL_2$ on $L_3$ are in natural bijection with the orbits of $GL_2$ on $GL_2$ by conjugation.  Furthermore, this bijection is compatible with the transpose (on $L_3$ and on $GL_2$).  

It follows that the $(GL_2 \times GL_2)$-invariant distributions on $L_3$ are in bijection with the $GL_2$-invariant distributions on $GL_2$ (for the conjugation action).  Since every element of $GL_2$ is conjugate to its transpose, we find that conjugation-invariant distributions on $GL_2$ are also transpose-invariant.  It follows that $(GL_2 \times GL_2)$-invariant distributions on $L_3$ are also transpose-invariant, finishing this case.
\item[$\sigma = b$]
For $\sigma = b$, we first whittle down the support of $(S, \psi, \trans)$-invariant distributions $T$ on $R = Q_3^\trans \sigma Q_3$.  Consider a general $(S^\trans,S)$ coset representative in $R$:  $g = \eta(u) \sigma \eta(v)$.  We require explicit forms for the entries of $u$ and $v$:
$$u = \Matrix{u_1}{u_2}{u_3}{u_4}, v = \Matrix{v_1}{v_2}{v_3}{v_4}.$$
We must consider two cases.  
\begin{description}
\item[$u_3 = 0$]
If $u_3 = 0$, then choose $\lambda_1, \lambda_2$ so that 
$$\psi_k(\lambda_1 v_1 + \lambda_2 v_2) \neq 1,$$ 
using the fact that $v$ is nonsingular.  Define
$$x = \Matrix{x_1}{x_2}{x_3}{x_4} = \Matrix{\lambda_1 v_1}{\lambda_1 v_2}{\lambda_2 v_1}{\lambda_2 v_2}, \mbox{ and }$$ 
$$s = \TMatrix{I}{x}{0}{0}{I}{J x^\trans J^{-1}}{0}{0}{I}.$$
We compute
$$g s g^{-1} =  \left(
\begin{array}{cccccc}
 1 & \lambda_1 & 0 & 0 &0 & 0 \\
 0 & 1 & 0 &0  & 0 & 0\\
 0 & u_1 \lambda_2 & 1 & 0 &0 & 0 \\
 0 & 0 & 0 & 1 & 0 & 0 \\
 0 & 0 &0 & -u_1 \lambda_2 & 1 & \lambda_1 \\
 0 & 0 & 0 & 0 & 0 & 1
\end{array}
\right)
.$$
We find that 
$$\psi_S(s) = \psi_k(x_1 + x_4) = \psi_k(\lambda_1 v_1 + \lambda_2 v_2) \neq 1,$$
$$\psi_S^\trans(g x g^{-1}) = \psi_k(0) = 1.$$ 
By criterion (R1), $R = S^\trans g S$ does not support any $(S,\psi,\trans)$-invariant distributions.
\item[$u_3 \neq 0$]
Suppose that $u_3 \neq 0$.  If $v_1 \neq 0$, then define
$$x = \Matrix{x_1}{x_2}{x_3}{x_4} = \Matrix{\lambda_1 v_1}{\lambda_1 v_2}{\lambda_2 v_1}{\lambda_2 v_2},$$
where $\lambda_1$ and $\lambda_2$ are chosen in such a way that $\lambda_2 \neq 0$ and 
$$x_1 + x_4 = \lambda_1 v_1 + \lambda_2 v_2 = 0.$$
Then since $v$ is nonsingular, we find that
$$x_4 v_3 - x_3 v_4 = \lambda_2 v_2 v_3 - \lambda_2 v_1 v_4 = -\lambda_2 \det(v) \neq 0.$$
Simplifying, 
$$g s g^{-1} = 
\left(
\begin{array}{cccccc}
 1 & \lambda_1 & 0 & 0 &0 & 0 \\
 0 & 1 & 0 &0  & 0 & 0\\
 0 & \lambda_2 u_1 & 1 & 0 &0 & 0 \\
 0 &  - \lambda_2 u_3 & 0 & 1 & 0 & 0 \\
 0 & 0 &\lambda_2 u_3  & - \lambda_2 u_1 & 1 & \lambda_1 \\
 0 & 0 & 0 & 0 & 0 & 1
\end{array}
\right)
.$$
By scaling the vector $(\lambda_1, \lambda_2)$ if necessary, we find that 
$$\psi_S(s) = \psi_k(x_1 + x_4) = \psi_k(0) = 1,$$
$$\psi_S^\trans(g s g^{-1}) = \psi_k(-\lambda_2 u_3) \neq 1.$$  
By criterion (R1), $R = S^\trans g S$ does not support any $(S,\psi,\trans)$-invariant distributions.

If $v_1 = 0$, and $v_2 \neq -u_3$ then we may choose $\lambda_2$ such that $\psi(-\lambda_2 u_3) \neq \psi(\lambda_2 v_2)$.  From this it follows that $\psi_S^\trans(g s g^{-1}) \neq \psi_S(s)$.  It follows that $R = S^\trans g S$ does not support any $(S,\psi,\trans)$-invariant distributions.
\end{description}

We find that all $(S,\psi,\trans)$-invariant distributions $T$ must be supported on cosets $S^\trans g S$ for which $g = \eta(u) \sigma \eta(v)$ with $v_1 = 0$, $u_3 \neq 0$, and $v_2 = - u_3$.  Applying (R2) instead of (R1), we can whittle down the support further (in a symmetric way), and we find that all $(S,\psi,\trans)$-invariant distributions $T$ must be supported on 
$$X = \{ S^\trans \eta(u) \sigma \eta(v) S : u_1 = v_1 = 0, v_2 = -u_3 \}.$$

Now, if $g = \eta(u) \sigma \eta(v)$ and $u_1 = v_1 = 0$, and $v_2 = -u_3$, consider the elements $z,y \in S$ given by:
$$z = \left(
\begin{array}{cccccc}
1 & 0 & 0 & 0 & 0 & 0 \\
0 & 1 & 0 & 0 & u_4 v_3 & 0 \\
0 & 0 & 1 & 0 & 0 & 0 \\
0 & 0 & 0 & 1 & 0 & 0 \\
0 & 0 & 0 & 0 & 1 & 0 \\
0 & 0 & 0 & 0 & 0 & 1
\end{array}
\right),
y = \left(
\begin{array}{cccccc}
1 & 0 & 0 & 0 & 0 & 0 \\
0 & 1 & 0 & 0 & v_4 u_2  & 0 \\
0 & 0 & 1 & 0 & 0 & 0 \\
0 & 0 & 0 & 1 & 0 & 0 \\
0 & 0 & 0 & 0 & 1 & 0 \\
0 & 0 & 0 & 0 & 0 & 1
\end{array}
\right).$$
Then we find that
$$y^\trans gz = \left(
\begin{array}{cccccc}
 1 & 0 & 0 & 0 & 0 & 0 \\
 0 & 0 & 0 & v_2 & 0 & 0 \\
 0 & 0 & 0 & 0 & -u_2 v_2 v_3 & 0 \\
 0 & v_2 & 0 & 0 & 0 & 0 \\
 0 & 0 & -u_2 v_2 v_3 & 0 & 0 & 0 \\
 0 & 0 & 0 & 0 & 0 & \det(uv)
\end{array}
\right).$$
Observing that $y^\trans g z$ is equal to its transpose, and $\psi_S(y) = \psi_S(z) = 1$, we find that $(S, \psi, \trans)$-invariant distributions on $X$ are also transpose-invariant.
\item[$\sigma = bcb$]
For $\sigma = bcb$, consider a general coset representative $g = \eta(u) \sigma \eta(v)$.  With $u,v$ as before, define $w = uv$, so that
$$w = \Matrix{w_1}{w_2}{w_3}{w_4} = \Matrix{u_1 v_1 + u_2 v_3}{u_1 v_2 + u_2 v_4}{u_3 v_1 + u_4 v_3}{u_3 v_2 + u_4 v_4}.$$
If $(w_3, w_4) \neq (-\det(v), 0)$, then there exist $x_3, x_4$ such that
$$\psi_k \left( \frac{w_4 x_3 - w_3 x_4}{\det(v)} \right) \neq 1, \mbox{ and } \psi_k(x_4) = 1.$$
In this case, we set $x_1 = x_2 = 0$ and $x_3, x_4$ satisfying the above conditions.  Define as in the previous case
 $$x = \Matrix{x_1}{x_2}{x_3}{x_4}, \mbox{ and } s = \TMatrix{I}{x}{0}{0}{I}{J x^\trans J^{-1}}{0}{0}{I}.$$
Then we compute
$$
g s g^{-1} = \left(
\begin{array}{cccccc}
 1 & 0 &  0 & 0 & 0 & 0 \\
 0 & 1 & 0 & 0 & 0 & 0 \\
 0 & \frac{w_2 x_3 - w_1 x_4}{\det(v)} & 1 & 0 & 0 & 0 \\
 0 & \frac{w_4 x_3 - w_3 x_4}{\det(v)} & 0 & 1 & 0 &  0 \\
 0 & 0 &  -\frac{w_4 x_3 - w_3 x_4}{\det(v)} &  \frac{w_2 x_3 - w_1 x_4}{\det(v)}  & 1 & 0 \\
 0 & 0 & 0 & 0 & 0 & 1
\end{array}
\right).$$
We find that 
$$\psi_S(s) = \psi_k(x_4) = 1,$$
$$\psi_S^\trans(g s g^{-1}) = \psi_k \left(  \frac{w_4 x_3 - w_3 x_4}{\det(v)} \right) \neq 1.$$
By (R1) it follows that if $T$ is a $(S,\psi,T)$-invariant distribution on $R$ then $T$ is supported on $(S^\trans,S)$-cosets of the form $S^\trans \eta(u) \sigma \eta(v) S$ for $u,v \in GL_2$ satisfying:
\begin{equation}
\label{uvEq}
uv = \Matrix{\ast}{\det(u)}{-\det(v)}{0}.
\end{equation}
For such $u,v$ we compute
$$g = 
\left(
\begin{array}{cccccc}
 1 & 0 & 0 & 0 & 0 & 0 \\
 0 & 0 & 0 & 0 & -\det(v) & 0 \\
 0 & 0 & \ast & \det(u) & 0 &0 \\
 0 & 0 & -\det(v) & 0 & 0 &  0 \\
 0 & \det(u) &  0 & 0  & 0 & 0 \\
 0 & 0 & 0 & 0 & 0 & \det(uv)
\end{array}
\right).
$$
A direct computation yields
$$\Delta \Matrix{\det(v)}{0}{0}{-\det(u)} g  \Delta \Matrix{\det(v)}{0}{0}{-\det(u)}^{-1} =  g^\trans.$$
By criterion (TI), we find that all $(S, \psi, T)$-invariant distributions on these cosets are also transpose-invariant.

\item[$\sigma = bacb$]
For $\sigma = bacb$, consider a general coset representative $g = \eta(u) \sigma \eta(v)$ with $u,v$ as before.

First, if $u_1 = 0$, then we may choose $x_1, x_2, x_3, x_4$ such that:
$$\psi_k(x_1 + x_4) \neq 1 \mbox{ and } v_2 x_1 - v_1 x_2 = 0.$$
For this choice, there exists $\lambda$ such that $x_1 = \lambda v_1$ and $x_2 = \lambda v_2$.  Define
$$x = \Matrix{x_1}{x_2}{0}{x_4}, s = \TMatrix{I}{x}{0}{0}{I}{J x^\trans J^{-1}}{0}{0}{I}.$$
Then we compute 
$$g s g^{-1} = \left(
\begin{array}{cccccc}
 1 & \frac{-v_1 x_4}{\det(v)} &0 & 0& 0 & 0 \\
 0 & 1 & 0 & 0 & 0 & 0 \\
 0& 0 & 1 & 0 & 0 & 0 \\
 u_3 \lambda & 0 & 0 & 1 & 0 & 0 \\
 \frac{x_4 v_3}{\det(u^{-1} v)} & 0 & 0 & 0 & 1 & \frac{-v_1 x_4}{\det(v)} \\
 0 & -\frac{x_4 v_3}{\det(u^{-1} v)} & u_3 \lambda &0& 0 & 1
\end{array}
\right).
$$
We find that $g s g^{-1} \in S^\trans$, 
$$\psi_S(s) = \psi_k(x_1 + x_4) \neq 1$$
$$\psi_S^\trans(g s g^{-1}) = \psi_k(0) = 1.$$
By (R1) the double coset $S^\trans \eta(u) \sigma \eta(v) S$ does not support any $(S, \psi, \trans)$-invariant distributions.

Next suppose that $u_1 \neq 0$.  Choose $\lambda$ such that $\psi_k(u_1 \lambda) \neq 1$.  Define
$$x = \Matrix{\lambda v_1}{\lambda v_2}{0}{- \lambda v_1}, \mbox{ and } $$
$$s = \TMatrix{I}{x}{0}{0}{I}{J x^\trans J^{-1}}{0}{0}{I}.$$

We compute
$$g s g^{-1} = 
\left(
\begin{array}{cccccc}
 1 & \frac{\lambda v_1^2}{\det(v)} & 0 & 0 & 0 & 0 \\
 0 & 1 & 0 & 0 & 0 & 0 \\
 u_1 \lambda  & 0 & 1 & 0 & 0 & 0 \\
 u_3 \lambda & 0 & 0 & 1 & 0 & 0 \\
 \frac{- \lambda v_1 v_3}{\det(u^{-1} v)} & 0 & 0 & 0 & 1 & \frac{\lambda v_1^2}{\det(v)} \\
 0 & \frac{\lambda v_1 v_3}{\det(u^{-1}v)} & u_3 \lambda & -u_1 \lambda & 0 & 1
\end{array}
\right).
$$
We find that $g s g^{-1} \in S^\trans$,
$$\psi_S(s) = \psi_k(x_1 + x_4) = \psi_k(0) = 1,$$
$$\psi_S^\trans(g s g^{-1}) = \psi_k(u_1 \lambda) \neq 1.$$
By (R1) the double coset $S^\trans \eta(u) \sigma \eta(v) S$ does not support any $(S, \psi, \trans)$-invariant distributions.

\item[$\sigma = bcabacb$]
Suppose that $g = \eta(u) \sigma \eta'(v)$ for $u,v \in GL_2$, where $\eta'(v) = \Delta(v) \eta(v)^{-1}$; we find it convenient to use slightly different coset representatives here, using $\eta'$ instead of $\eta$.  There are two cases to consider.

First, suppose that $u \det(v)^{-1} =-v \in GL_2$.  Then we find that
$$g = \TMatrix{0}{0}{-u}{0}{u}{0}{-u}{0}{0}.$$
Note that  there exists $\gamma \in GL_2$ such that $\gamma u \gamma^{-1} = u^\trans$.  It follows that
$$\Delta(\gamma) g \Delta(\gamma)^{-1} = g^\trans.$$
By (TI), any $(S, \psi, \trans)$-invariant distribution on $S^\trans g S$ will be transpose-invariant.

Next, suppose that $u \det(v)^{-1} \neq -v$.  Then we may choose $X \in M_2$ such that
$$\psi_k(-\Tr(X)) \neq \psi_k(\Tr(\det(v)^{-1} u X v^{-1})).$$ 
Define an element $s \in S$ by
$$s = \TMatrix{I}{X}{0}{0}{I}{J X^\trans X}{0}{0}{I}.$$
Then we find that
$$g s g^{-1} = \TMatrix{I}{0}{0}{\det(v)^{-1} u J X^\trans J  v^{-1}}{I}{0}{0}{\det(u)v X u^{-1}}{I}.$$
We find that $g s g^{-1} \in S^\trans$ and
$$\psi_S(s) = \psi_k(-\Tr(X)),$$
$$\psi_S(g s g^{-1}) = \psi_k(\Tr(\det(v)^{-1} u X v^{-1})).$$
By (R1), the coset $S^\trans g S$ does not support any $(S, \psi, \trans)$-invariant distributions.

\end{description}
\qed

With this technical work done, we now find
\begin{thm}
\label{USF}
Suppose that $\sigma$ is a supercuspidal irrep of $GSp_6$.  Then the space of Shalika functionals for $\sigma$ is at most one-dimensional:
$$\dim \left( \Sh(\sigma) \right) \leq 1.$$
\end{thm}
\proof
Our previous results on distributions, with the methods of Gelfand, Kazhdan, and Bernstein imply that the pair $(GSp_6, S)$ is a Gelfand pair, in the sense of Condition 4.1 of \cite{Gro} (though we work with the character $\psi_S$ of $S$ rather than the trivial representation of $S$).  To be precise, for an irrep $\sigma$ of $GSp_6$, with contragredient $\tilde \sigma$, we find (cf. Proposition 4.2 of \cite{Gro}) that
$$\dim \left( \Sh(\sigma) \right) \cdot \dim \left( \Sh(\tilde \sigma) \right) \leq 1.$$

So it remains to check that $\sigma$ has a nonvanishing Shalika functional if and only if $\tilde \sigma$ has a nonvanishing Shalika functional.  

Since $S$ is a unimodular subgroup of $GSp_6$, there is a nondegenerate $GSp_6$-invariant pairing:
$$\cInd_S^{GSp_6} \psi_S \times \cInd_S^{GSp_6} \tilde \psi_S \rightarrow \CC,$$
given by integration of functions on $S \backslash G$:
$$\langle f_1, f_2 \rangle = \int_{S \backslash G} f_1(g) f_2(g) dg.$$
Now, if $\sigma$ is a supercuspidal irrep of $GSp_6$ with nonvanishing Shalika functional, then $\sigma$ occurs as a subrepresentation of $\cInd_S^{GSp_6} \psi_S$.  The nondegeneracy of the pairing above (and the injectivity of supercuspidals) implies that $\tilde \sigma$ occurs as a subrepresentation of $\cInd_S^{GSp_6} \tilde \psi_S$.  It follows that $\tilde \sigma$ has a nonvanishing Shalika functional, with respect to the character $\tilde \psi_S$.  But since $\tilde \psi_S$ and $\psi_S$ are conjugate (via an element of $GSp_6$) characters of $S$, we find that $\Sh(\tilde \sigma) \neq 0$.
\qed

\subsection{Theta Correspondence}

The importance of Shalika functionals in the theta correspondence is the following:
\begin{lem}
Suppose that $\sigma$ is a generic supercuspidal irrep of $PGSp_6$.  Then there is a linear isomorphism:
$$\Wh_{G_2}(\lvec{\Theta}_7(\sigma)) = \lvec{\Theta}_7(\sigma)_{N_2, \psi_2} \isom \Sh(\sigma).$$
\end{lem}
\proof
Here  $\alg{N}_{2}$ be the unipotent radical of a Borel subgroup of $\alg{G}_{2}$ and $\psi_{2}$ is a principal character of $N_{2}$.  Let $\alg{Q}_2= \alg{L}_2 \alg{U}_2$ be a maximal parabolic subgroup of $\alg{G}_2$ such that $\alg{U}_2$ is contained in $\alg{N}_{2}$ and $\alg{N}_{2}/\alg{U}_{2}$ corresponds to a short simple root. 

Then $\Wh_{G_2}(\Pi_7) = (\Pi_{7})_{N_2,\psi_{2}}$ can be computed in two stages:
$$(\Pi_{7})_{N_2,\psi_{2}} = \left( (\Pi_{7})_{U_{2},\psi_{2}} \right)_{N_2, \psi_2}.$$
Lemma 2.9 on page 213 in \cite{GrS1} shows how to compute the co-invariants of $\Pi_{7}$ with respect to any character of $U_{2}$. The characters of $U_{2}$ are parameterized by cubic $k$-algebras, and the restriction of $\psi_{2}$ to $U_{2}$ corresponds to the degenerate cubic algebra $k[\epsilon] / \langle \epsilon^3 \rangle$.

Let $S^{\circ}\subseteq S$ be the semidirect product of $GL_{2}$ with $U^{\circ}_{3} \subseteq U_{3}$, where $U^{\circ}_{3}$ contains the center $Z_{3}$ and $U^{\circ}_{3}/Z_{3}$ corresponds to trace zero matrices in $U_{3}/Z_{3}\cong M_{2}(k)$. Then
$$(\Pi_{7})_{U_{2},\psi_{2}}\cong \cInd^{GSp_{6}}_{S^{\circ}}( \complex ).$$
Under this identification, one observes that the action of $N_{2}$ on $(\Pi_{7})_{U_{2},\psi_{2}}$ (which restricts to the character $\psi_2$ on $U_2$) is identified with the action of $S / S^\circ$ by left translation on $\cInd_{S^{\circ}}^{GSp_{6}}( \complex )$. This implies that 
$$\Wh_{G_2}(\Pi_7) = (\Pi_{7})_{N_{2},\psi_{2}} \isom \cInd_{S}^{GSp_{6}}(\psi_{S}),$$
as representations of $GSp_6$.  

Applying $\Hom_{GSp_6}(\sigma, \cdot)$ to both sides above, we find that 
$$\Wh_{G_2}(\lvec{\Theta}_7(\sigma)) = \lvec{\Theta}_7(\sigma)_{N_2, \psi_2} \isom \Sh(\sigma).$$
\qed

Since $\lvec{\Theta}_7(\sigma)$ is multiplicity-free, supercuspidal, and every subrepresentation is generic, we immediately find that:
\begin{prop}
\label{PerProp}
Suppose that $\sigma$ is a generic supercuspidal irrep of $PGSp_6$, with trivial central character.  Then $\lvec{\Theta}_7(\sigma)$ is nonzero if and only if $\Sh(\sigma) \neq 0$.  Moreover, if $\lvec{\Theta}_7(\sigma) \neq 0$ then $\lvec{\Theta}_7(\sigma)$ is a generic supercuspidal irrep of $G_2$.
\end{prop}
\proof
This proposition directly follows from the previous lemma, and the ``uniqueness of Shalika functionals'' of Theorem \ref{USF}
\qed

\subsection{Injectivity of Dichotomy}

We can now demonstrate the following
\begin{thm}
The dichotomy map is injective:
$$\Delta : \Irr_g^\circ(G_2) \hookrightarrow \Irr_{g}^\circ(PGSp_6) \sqcup \frac{\Irr_g^\circ(PGL_3)}{\Contra}.$$
\end{thm}
If two generic supercuspidal irreps $\tau, \tau'$ of $G_2$ have the property that $\rvec{\Theta}_6(\tau)$ and $\rvec{\Theta}_6(\tau')$ have a common supercuspidal subrepresentation, then $\tau$ is isomorphic to $\tau'$ by Theorem 19 of \cite{GS1}.

If two generic supercuspidal irreps $\tau, \tau'$ of $G_2$ have the property that $\rvec{\Theta}_7(\tau)$ and $\rvec{\Theta}_7(\tau')$ have a common generic supercuspidal subrepresentation $\sigma$, then $\tau$ is isomorphic to $\tau'$ by Proposition \ref{PerProp} since both $\tau$ and $\tau'$ must be subrepresentations of the irrep $\lvec{\Theta}_7(\sigma)$.  
\qed

\section{L-functions and periods}

Now that we have proven that the dichotomy map is injective, it remains to characterize its image.  In fact, all supercuspidal irreps of $PGL_3$ occur in the theta correspondence with a generic supercuspidal irrep of $G_2$, by Theorem 19 of Gan and Savin \cite{GS1}
\begin{prop}
Suppose that $\rho$ is a supercuspidal irrep of $PGL_3$.  Then there exists a unique generic supercuspidal irrep $\tau$ of $G_2$ occurring in $\lvec{\Theta}_6(\rho)$.
\end{prop}
This immediately implies that
\begin{cor}
The image of dichotomy in $\Irr_g^\circ(PGL_3)$ includes all non-self-contragredient supercuspidal irreps.  In particular, $\Delta$ surjects onto $\Irr_g^\circ(PGL_3)$ when $p \neq 2$.
\end{cor}

On the other hand, the image of dichotomy in $\Irr_g^\circ(PGSp_6)$ is so far only characterized as
$$\Delta \left( \Irr_g^\circ(G_2) \right) \cap \Irr_g^\circ(PGSp_6) = \{ \sigma \in \Irr_g^\circ(PGSp_6) : \Sh(\sigma) \neq 0 \}.$$
In this section, we demonstrate that the image of dichotomy can be described not only by the Shalika functional, but also by the degree 8 spin L-function.  The goal of this section is to prove the following
\begin{thm}
\label{SpinShalika}
Suppose that $\sigma$ is a generic supercuspidal irrep of $PGSp_6$.  Let $L(\sigma, Spin, s)$ denote Shahidi's L-function, associated to the 8-dimensional spin representation of $Spin_7(\CC)$.  Then $\Sh(\sigma) \neq 0$ if and only if $L(\sigma, Spin, s)$ has a pole at $s = 0$.
\end{thm}

One direction in this theorem -- that a nonvanishing Shalika functional implies that $L(\sigma, Spin, s)$ has a pole at $s=0$ -- follows from Shahidi's work, examination of a reducibility point, and properties of the minimal representation of $E_8$.  The other direction relies on an integral representation for the spin L-function due to Bump-Ginzburg \cite{BG} and studied by Vo \cite{Vo}.  We prove that these two incarnations of the spin L-function have the same poles, using global methods.

\subsection{A reducibility point}

Let $\alg{P}_4 = \alg{M}_4 \alg{N}_4$ be the Heisenberg parabolic subgroup of $\alg{F}_4$, with Levi component $\alg{M}_4 \isom \alg{GSp}_6$ and $\Sim$ is the similitude character of $\alg{GSp}_6$.  The modular character, for the adjoint action of $M_4$ on $N_4$ can then be expressed as:
$$\delta_{P_4}(m) = \vert \Sim(m) \vert^8.$$

For $\sigma$ a generic supercuspidal irrep of $PGSp_6$, consider the family of representations of $F_4$:
$$I(\sigma, s) = \Ind_{P_4}^{F_4} (\sigma \otimes \vert \Sim \vert^{s+4} ),$$
where $\Sim$ is the similitude character of $GSp_6$.  The normalization factor $\vert \Sim \vert^{4}$ is chosen so that $I(\sigma, 0)$ is unitary when $\sigma$ is unitary.

Let $L(\sigma, Spin, s)$ be Shahidi's $L$-function, where $Spin$ is the $8$-dimensional representation of the dual Levi $\alg{\hat M}_4 \isom \alg{CSpin}_7(\complex)$ on the abelian quotient of the unipotent radical of the parabolic $\alg{\hat P}_4$ dual to $\alg{P}_4$.  The following result is essentially due to Shahidi \cite{Sha}:
\begin{lem}
The L-function $L(\sigma, Spin, s)$ has a pole at $s = 0$ if and only if $I(\sigma, -1)$ is reducible, in which case it has a composition series of length two.   In this case, the unique irreducible submodule $J(\sigma)$ of $I(\sigma, -1)$ is not generic.
\end{lem}
\proof
To compute this reducibility point, we compute some constants discussed in \cite{Sha}.  Let $\alpha_1, \ldots, \alpha_4$ denote the simple roots in a root system of type $\Type{F}_4$, numbered as below.
\begin{center}
\begin{tikzpicture}[scale=2]
\filldraw[fill=gray] (0,0) circle (.05cm);
\draw (0,.2cm) node[above] {$P_4=GSp_6 \ltimes N_4$};
\draw (0,0) node[below] {$\alpha_1$};
\draw (.05cm,0) -- (.95cm,0);
\draw (1,0) circle (.05cm);
\draw (1,0) node[below] {$\alpha_2$};
\draw (1.05cm,.025cm) -- (1.95cm,0.025cm);
\draw (1.05cm,-.025cm) -- (1.95cm, -0.025cm);
\draw (1.45cm,-.1cm) -- (1.55cm,0) -- (1.45cm,.1cm);
\draw (2,0) circle (.05cm);
\draw (2,0) node[below] {$\alpha_3$};
\draw (2.05cm,0) -- (2.95cm,0);
\draw (3,0) circle (.05cm);
\draw (3,0) node[below] {$\alpha_4$};
\end{tikzpicture}
\end{center}
Let $\beta$ denote the highest root, so that:
$$\beta = 2 \alpha_1 + 3 \alpha_2 + 4 \alpha_3 + 2 \alpha_4.$$
Observe that the maximal parabolic subgroup $\alg{P}_4$ is associated to the root $\alpha_1$, which is adjacent to $- \beta$ in the extended (affine) Dynkin diagram.

Let $\rho_P$ denote the half-sum of the roots occurring in $\alg{N}_4$.  Then $\rho_P = 4 \beta$.  It follows that:
$$\tilde \alpha_1 = \langle \alpha_1^\vee, \rho_P \rangle^{-1} \cdot \rho_P = \langle \alpha_1^\vee, \beta \rangle^{-1} \beta = \beta.$$
Since $\beta$ corresponds precisely to the similitude character of $\alg{M}_4 = \alg{GSp}_6$, it follows that $I(\sigma,s)$ is normalized as in Shahidi \cite{Sha}.  The result now follows directly from \cite{Sha}; a helpful exposition of the results from Shahidi can be found in Section 2 of \cite{Zha}.
\qed

To demonstrate a connection between nonvanishing of a theta correspondence and L-functions, we use a method of Mui{\'c}-Savin \cite{MuS} and consider a theta correspondence in a larger group.  The following lemma plays a similar role in this section to Proposition 4.1 in \cite{MuS}.  
\begin{lem}
Let $\Pi_8$ denote the minimal representation of $E_8$.  Let $\phi_4$ be a generic character of a maximal unipotent subgroup $U_4$ of $F_4$ .  Then 
$$\Wh_{F_4}(\Pi_8) = (\Pi_8)_{U_4, \phi_4}  = 0.$$
\end{lem}
\proof
We study the Whittaker functionals $\Wh_{F_4}(\Pi_8) = (\Pi_8)_{U_4, \phi_4}$ in stages:
$$(\Pi_8)_{U_4, \phi_4} = \left( \left( (\Pi_8)_{N_4, \psi_4} \right)_{N_3, \psi_3} \right)_{U_2, \psi_2}.$$
where $N_4$ is a 15-dimensional Heisenberg group in $F_4$, $N_3$ is a 6-dimensional abelian unipotent subgroup of $GSp_6$, and $U_2$ is a maximal unipotent subgroup of $SL_3$.  

\textit{Stage 1:  The $N_4,\psi_4$ coinvariants.}
We view $\alg{F}_4$ here as the algebraic group associated to the 14-dimensional structurable algebra of Freudenthal type
$$F_k \isom k \oplus J_k \oplus J_k \oplus k.$$
Similarly, we view $\alg{E}_8$ as the algebraic group associated to the 56-dimensional structurable algebra of Freudenthal type
$$F_\octs \isom k \oplus J_\octs \oplus J_\octs \oplus k.$$
The construction of these algebras and groups follows Section \ref{FA}.  As a result, $\alg{F}_4$ is endowed with a parabolic subgroup $\alg{P}_4 = \alg{M}_4 \alg{N}_4$, and $\alg{E}_8$ contains a parabolic subgroup $\alg{P}_8 = \alg{M}_8 \alg{N}_8$, such that:
\begin{enumerate}
\item
$\alg{N}_4$ and $\alg{N}_8$ are two-step unipotent groups with one-dimensional centers $\alg{Z}_4$ and $\alg{Z}_8$.
\item
$\alg{N}_4 / \alg{Z}_4$ is naturally identified with $F_k$, and $\alg{N}_8 / \alg{Z}_8$ is naturally identified with $F_\octs$.
\item
The parabolics are aligned, in the sense that $\alg{P}_8 \cap \alg{F}_4 = \alg{P}_4$, $\alg{N}_8 \cap \alg{F}_4 = \alg{N}_4$, and $\alg{Z}_8 = \alg{Z}_4$.
\end{enumerate}

Let $\psi_4$ denote the restriction of $\phi_4$ to $N_4$.  Then $\psi_4$ is in the minimal $GSp_6$-orbit in the space of characters of $N_4$.  By conjugation, we may assume that $\psi_4$ corresponds to the element $\Matrix{1}{0}{0}{0} \in F_k \isom k \oplus J_k \oplus J_k \oplus k$ (identified with $N_4^- / Z_4^-$). 

The space $(\Pi_8)_{N_4, \psi_4}$ is a quotient of the kernel
$$Ker \left( (\Pi_8)_{Z_8} \rightarrow (\Pi_8)_{N_8} \right).$$
From Corollary 11.12 of \cite{GS4}, this kernel can be identified with $C_c^\infty(\Omega)$, where $\Omega$ is the minimal nonzero $M_8$-orbit in the 56-dimensional minuscule representation $N_8^- / Z_8^-$.  

Then the characters of $N_8$ which restrict to $\psi_4$ on $N_4$, and also are in the $N_8 / Z_8$-support of $\Pi_8$ correspond to elements
$$\Matrix{1}{j}{j^\sharp}{0} \in F_\octs \isom k \oplus J_\octs \oplus J_\octs \oplus k,$$
where $j \in J_\octs$, the entries of $j$ and $j^\sharp$ are in $\octs_\circ$, $j^\sharp$ is the quadratic adjoint of $j$, and $\Norm(j) = 0$.  Here we refer to Proposition 11.2 and Section 10 of \cite{GS4} for a description of the orbit $\Omega$ and Jordan algebras.
  
Thus the representation $(\Pi_8)_{N_4, \psi_4}$ is identified with $C_c^\infty(\Omega^\perp)$, where
$$\Omega^\perp = \lset
\left(\begin{array}{ccc}0 & \alpha & - \beta \\ - \alpha & 0 & \gamma \\ \beta & - \gamma & 0\end{array}\right) : \begin{array}{c}
\alpha^2 = \beta^2 = \gamma^2 = 0  \\
\Tr(\alpha \beta) = \Tr(\beta \gamma) = \Tr(\gamma \alpha) = 0 \end{array},  \Tr(\alpha \beta \gamma) = 0 \rset.$$
Equivalently, we may view
$$\Omega^\perp = \{ (\alpha, \beta, \gamma) \in \octs_\circ^3 : \Span_k(\alpha, \beta, \gamma) \mbox{ is isotropic}, \Tr(\alpha \beta \gamma) = 0 \}.$$

\textit{Stage 2:  The $N_3,\psi_3$ coinvariants.}  Now we are led to consider
$$\left( (\Pi_8)_{N_4, \psi_4} \right)_{N_3, \psi_3} \isom C_c^\infty(\Omega^\perp)_{N_3, \psi_3}.$$

First, we describe the subgroup $N_3$ of $GSp_6$, and the character $\psi_3$.  Here, $\alg{Q}_3 = \alg{M}_3 \alg{N}_3$ denotes the ``Siegel parabolic'' in $GSp_6$, whose derived subgroup is $\alg{M}_3' \isom \alg{SL}_3 \ltimes \alg{N}_3$.  This derived subgroup $\alg{M}_3'$ stabilizes the character $\psi_4$ of $N_4$.  The unipotent radical $\alg{N}_3$ of $\alg{Q}_3$ is abelian, and its $k$-points $N_3$ are identified naturally with the space $J_k$ of symmetric 3 by 3 matrices with entries in $k$.  Then $N_3$ acts on $\Omega^\perp$ in the following way:
$$\kappa \star j = j + (j^\sharp \times \kappa) \mbox{ for all } \kappa \in N_3 = J_k, j \in \Omega^\perp \subset J_\octs.$$
In particular, we can compute
$$
\left(\begin{array}{ccc} 1 & 0 & 0 \\ 0 & 0 & 0 \\ 0 & 0 & 0 \end{array}\right) \star \left(\begin{array}{ccc}0 & \alpha & - \beta \\ - \alpha & 0 & \gamma \\ \beta & -\gamma & 0\end{array}\right)  = 
\left(\begin{array}{ccc} 0 & \alpha & - \beta \\ -\alpha & 0 & \gamma -\alpha \beta \\ \beta & \gamma + \alpha \beta & 0\end{array}\right),
$$
or in shorthand,
$$\left(\begin{array}{ccc} 1 & 0 & 0 \\ 0 & 0 & 0 \\ 0 & 0 & 0 \end{array}\right) \star (\alpha, \beta, \gamma) = (\alpha, \beta, \gamma - \alpha \beta).$$
Let $\psi_3$ be the character of $N_3$ given by
$$\psi_3 \left(\begin{array}{ccc} a & r & s \\ r & b & t \\ s & t & c \end{array}\right)  = \psi_k(a).$$

We now decompose $\Omega^\perp$ into two subsets:
$$\Omega_1^\perp = \{ (\alpha, \beta, \gamma) \in \Omega^\perp : \alpha \beta = 0 \},$$
$$\Omega_2^\perp = \{ (\alpha, \beta, \gamma) \in \Omega^\perp : \alpha \beta \neq 0 \}.$$

We find almost immediately that $C_c^\infty(\Omega_1^\perp)_{N_3, \psi_3} = 0$:  indeed, if $j \in \Omega_1^\perp$, $t \in k$ and $\psi_k(t) \neq 0$, then 
$$\psi_3 \left(\begin{array}{ccc} t & 0 & 0 \\ 0 & 0 & 0 \\ 0 & 0 & 0 \end{array}\right) \neq 1 \mbox{ and } \left(\begin{array}{ccc} t & 0 & 0 \\ 0 & 0 & 0 \\ 0 & 0 & 0 \end{array}\right)  \star j = j.$$

In other words $\Omega_1^\perp$ does not support any $(N_3,\psi_3)$-invariant distributions.  It follows that
$$\left( (\Pi_8)_{N_4, \psi_4} \right)_{N_3, \psi_3} = C_c^\infty(\Omega^\perp)_{N_3, \psi_3} = C_c^\infty(\Omega_2^\perp)_{N_3, \psi_3}.$$
Note that if $j = (\alpha, \beta, \gamma) \in \Omega_2^\perp$, then $\alpha \beta$ is nonzero and orthogonal to $\alpha$, $\beta$, $\gamma$, and itself (with respect to the trace pairing on $\octs_\circ$.  But the maximal dimension of an isotropic subspace in $\octs_\circ$ is three, so there must exist $a,b,c \in k$, not all zero, such that
$$\alpha \beta = a \alpha + b \beta + c \gamma.$$
Multiplying through by $\alpha$ or by $\beta$, we find that
$$b \alpha \beta = -c \alpha \gamma,$$
$$a \alpha \beta = -c \gamma \beta.$$
Hence $\Span _k(\alpha \beta, \beta \gamma, \gamma \alpha)$ is one-dimensional.  Note also that $c \neq 0$ in the above relations, since otherwise $\alpha \beta = 0$.  

\textit{Stage 3:  The $U_2,\psi_2$ coinvariants.}  Now we are led to consider the coinvariants
$$\Wh_{F_4}(\Pi_8) = (\Pi_8)_{U_4, \phi_4} = \left( \left( (\Pi_8)_{N_4, \psi_4} \right)_{N_3, \psi_3} \right)_{U_2, \psi_2} = \left( C_c^\infty(\Omega_2^\perp)_{N_3, \psi_3} \right)_{U_2, \psi_2}.$$
We describe the subgroup $U_2$ and character $\psi_2$ here.  There is a chain of embeddings:
$$U_2 \subset SL_3 \subset Q_3' = SL_3 \ltimes J_k  \subset Q_3 \subset GSp_6 \subset Q_4 = GSp_6 \ltimes F_k \subset F_4.$$  
The resulting action of $SL_3$ on $F_k \isom k \oplus J_k \oplus J_k \oplus k$ is given by the formulas of Section 3.1 of Krutelevich \cite{Kru}.
$$\gamma \cdot \Matrix{a}{A}{B}{b} = \Matrix{a}{\gamma A \gamma^\trans}{(\gamma^{-1})^{\trans} B \gamma^{-1}}{b}, \mbox{ for all } \gamma \in SL_3.$$
In particular, let $U_2$ denote the standard maximal unipotent subgroup of this $SL_3$; the action of $U_2$ on $\Omega^\perp$ is given by
$$\left(\begin{array}{ccc}1  & x & z \\ 0 & 1 & y \\ 0 & 0 & 1 \end{array}\right) \cdot (\alpha, \beta, \gamma) = (\alpha + (xy-z) \gamma + y  \beta, \beta - x \gamma, \gamma).$$
We define $\psi_2$ to be the principal character of $U_2$ given by
$$\psi_2 \left(\begin{array}{ccc}1  & x & z \\ 0 & 1 & y \\ 0 & 0 & 1 \end{array}\right) = \psi_k(x - y).$$

Together with the action of $N_3$, we find an action of $SL_3 \ltimes N_3$ on $\Omega^\perp$:  for all $\gamma \in SL_3$ and all $\kappa \in N_3$,
$$\gamma \cdot (\kappa \star j) = (\gamma \kappa \gamma^{\trans}) \star (\gamma \cdot j).$$
The subgroup $U_2$ of $SL_3$ stabilizes the character $\psi_3$ of $N_3$:
$$\psi_3(u \kappa u^\trans) = \psi_3(\kappa), \mbox{ for all } u \in U_2, \kappa \in N_3.$$

Now for $j = (\alpha, \beta, \gamma) \in \Omega_2^\perp$, so that $\alpha \beta \neq 0$, we find three possibilities:
\begin{description}
\item[$\alpha \gamma \neq 0$]
If $\alpha \gamma \neq 0$, then there exists $x \in k$ such that $\alpha \beta - x \alpha \gamma = 0$.  We find that  $j = (\alpha, \beta, \gamma)$ is in the same $U_2$-orbit as $j' =(\alpha, \beta - x \gamma, \gamma)$, and $(\alpha)(\beta - x \gamma) = 0$.  Thus $j' \in \Omega_1^\perp$, and cannot be contained in the support a $(N_3, \psi_3)$-invariant distribution by the result of Stage 2.  
\item[$\beta \gamma \neq 0$]
If $\beta \gamma \neq 0$, then there exists $z \in k$ such that $\alpha \beta - z \gamma \beta = 0$.  We find that $j=(\alpha, \beta, \gamma)$ is in the same $U_2$-orbit as $j'=(\alpha - z \gamma, \beta, \gamma)$, and $(\alpha - z \gamma)(\beta) = 0$.  Such elements $j' \in \Omega_1^\perp$ cannot be in the support a $(N_3, \psi_3)$-invariant distribution, again by the result of Stage 2.
\item[$\alpha \gamma = \beta \gamma = 0$]
If $\alpha \gamma = \beta \gamma = 0$, then we find that $a = b = 0$ in the linear dependence
$$\alpha \beta = a \alpha + b \beta + c \gamma.$$
Thus $\alpha \beta = c \gamma$.  Define an element $j'$ in the $N_3$-orbit of $j$ by
$$j' = \left(\begin{array}{ccc} c^{-1} & 0 & 0 \\ 0 & 0 & 0 \\ 0 & 0 & 0 \end{array}\right)  \star j = (\alpha, \beta, 0).$$
Then $j'$ cannot be in the support of a $(U_2, \psi_2)$-invariant distribution, since for any $x \in k$ such that $\psi_k(x) \neq 1$, we have
$$\psi_2 \left(\begin{array}{ccc}1  & x & 0 \\ 0 & 1 & 0 \\ 0 & 0 & 1 \end{array}\right)  \neq 1 \mbox{ and } \left(\begin{array}{ccc}1  & x & 0 \\ 0 & 1 & 0 \\ 0 & 0 & 1 \end{array}\right) \cdot j' = j'.$$
\end{description}
It follows that $\left( C_c^\infty(\Omega_2^\perp)_{N_3, \psi_3} \right)_{U_2, \psi_2} = 0$, and so
$$\Wh_{F_4}(\Pi_8) = 0.$$
\qed

Now we can demonstrate a connection between a nonvanishing Shalika functional and a pole in Shahidi's L-function.
\begin{thm}
\label{RedProp}
Suppose that $\sigma$ is a generic supercuspidal irrep of $PGSp_6$, with nonzero Shalika functional.  Let $\tau = \lvec{\Theta}_{7}(\sigma)$, a generic supercuspidal irrep by Proposition \ref{PerProp}.  Then the following statements are true:
\begin{enumerate}
\item
The $L$-function $L(\sigma,Spin,s)$ has a pole at $s=0$.
\item
If $J(\sigma)$ is the unique irreducible subrepresentation of $I(\sigma, -1)$, then $J(\sigma) \boxtimes \tau$ occurs as a quotient of the minimal representation $\Pi_8$ of $E_8$,  restricted to the dual pair $F_4 \times G_2$.
\end{enumerate}
\end{thm}
\proof
We use the Heisenberg parabolic subgroups $\alg{P}_4$, $\alg{P}_8$ of $\alg{F}_4$, $\alg{E}_8$, discussed in the previous result.

Recall that $\sigma \boxtimes \tau$ occurs as a quotient (or subrepresentation) of the minimal representation $\Pi_7$ of $E_7$.  By Theorem 6.1 of \cite{MS1} (following Proposition 4.1 of \cite{Sav}, and not requiring any condition on residue characteristic), the Jacquet functor (along $N_8$) of the minimal representation $\Pi_8$ of $E_8$ can be identified as a representation of $CE_7$:
$$(\Pi_8)_{N_8} \isom \left( \Pi_7 \otimes \vert \det \vert^{3/28} \right) \oplus \vert \det \vert^{5/28}.$$

Since $N_4 \subset N_8$, $(\Pi_8)_{N_8}$ is a quotient of $(\Pi_8)_{N_4}$, as representations of $GSp_6 \times G_2 \subset CE_7$.  It follows that there is a surjective $GSp_6 \times G_2$ intertwining map:
$$(\Pi_8)_{N_4} \twoheadrightarrow \Pi_7 \otimes \vert \det \vert^{3/28}.$$
Since $\sigma \boxtimes \tau$ occurs as a quotient of $\Pi_7$, restricted to $GSp_6 \times G_2$, we find a surjective $GSp_6 \times G_2$ intertwining map:
$$(\Pi_8)_{N_4} \twoheadrightarrow (\sigma \boxtimes \tau) \otimes \vert \det \vert^{3/28}.$$

It follows by Frobenius reciprocity that there is a nontrivial $F_4 \times G_2$ intertwining map:
$$\Pi_8 \rightarrow \Ind_{P_4}^{F_4} (\sigma \otimes \vert \det \vert^{3/28}) \boxtimes \tau.$$
In order to identify the restriction of $\det$ (the determinant for the action of $CE_7$ on a 56-dimensional space) to $GSp_6$, we consider the coroot $\alpha^\vee$ of $\alg{F}_4$ which satisfies:
$$\alpha^\vee(t) \in Z(GSp_6), \mbox{ for all } t \in k^\times, \mbox{ and } \Sim(\alpha^\vee(t)) = t^2.$$
This is the coroot of the $\alg{SL}_2$ which commutes with $\alg{Sp}_6 = [\alg{M}_4, \alg{M}_4]$ in $\alg{F}_4$.  This $\alg{SL}_2$ is identified with the $\alg{SL}_2$ which commutes with $\alg{E}_7 = [\alg{M}_8, \alg{M}_8]$ in $\alg{E}_8$.  Indeed, both copies of $\alg{SL}_2$ arise as $\alg{Aut}_{\octs/M_2}$, embedded in $\alg{G}_{\octs \otimes \octs} \isom \alg{E}_8$ and in $\alg{G}_{k \otimes \octs} \isom \alg{F}_4$.  We refer to Section \ref{TC} for a construction of these groups from tensor products of composition algebras.

The character $\det$ of $CE_7$, considered above, pairs with $\alpha^\vee$, in such a way that:
$$\det(\alpha^\vee(t)) = t^{56} = \Sim(\alpha^\vee(t))^{28}.$$
Indeed, $\alpha^\vee(t)$ acts on the 56-dimensional space $N_8/ [N_8, N_8]$ by the scalar $t$, and the determinant is computed above.  Comparing with the similitude character, for every element $m$ of the subgroup $GSp_6 \subset CE_7$, one has:
$$\vert \det(m) \vert^{3/28} = \vert \Sim(m) \vert^{3}.$$
Hence we find a nontrivial $F_4 \times G_2$ intertwining map:
$$\Pi_8 \rightarrow \Ind_{P_4}^{F_4} (\sigma \otimes \vert \Sim \vert^{3}) \boxtimes \tau = I(\sigma, -1) \boxtimes \tau.$$

Since $I(\sigma, -1)$ is generic, we find that $\Wh_{F_4}(I(\sigma, -1)) = I(\sigma, -1)_{U_4, \phi_4}$ is nonzero, where $\phi_4$ is a generic character of $U_4$ as before.  But $\Wh_{F_4}(\Pi_8) = (\Pi_8)_{U_4, \phi_4} = 0$ by the previous lemma.  It follows that the image of the above intertwining map must be a proper submodule of $I(\sigma, -1) \boxtimes \tau$.  Thus we get both statements at once:
\begin{enumerate}
\item
$I(\sigma, {-1})$ is reducible and by the work of Shahidi, $L(\sigma, Spin, s)$ has a pole at $s=0$.
\item
$J(\sigma) \boxtimes \tau$ occurs as a quotient of $\Pi_8$ (restricted from $E_8$ to $F_4 \times G_2$).
\end{enumerate}
\qed

\subsection{Eisenstein series}
Here we review Eisenstein series on $GL_{2}$, as they are used in the construction of the spin L-function by Bump and Ginzburg \cite{BG}.  Let $F$ be a global field with adele ring $\adeles$. Following \cite{GS}, page 47, for every place $v$ of $F$ and $s$ in $\mathbb C$ we define $V(s)$ to be the local unramified principal series representation of $\alg{GL}_{2}(F_{v})$, 
unnormalized, so that the trivial representation is a submodule of $V(0)$ and a quotient of $V(1)$.  Here $F_v$ is the completion of $F$ at $v$; $q_v$ will denote the cardinality of the residue field at $v$, if $v$ is a finite place. 

We have an intertwining operator $M_{v}(s): V(s)\rightarrow V(1-s)$ defined by 
$$M_{v}(s)(f_{v,s})(g)= \int_{N_{v}} f_{v,s}(wng)~dn$$
where $f_{v,s}$ is in $V(s)$. 
Let $f^0_{v,s}$ be the spherical vector in $V(s)$ normalized so that $f_{v,s}^0(1)=1$. Then, (see \cite{GS} page 51)  
$$M_{v}(s) f^{0}_{v,s}= \frac{L_{v}(2s-1)}{L_{v}(2s)} f_{v,s}^{0}$$
where $L_{v}(s)=(1-q_{v}^{-s})^{-1}$. We normalize the operator $M_{v}(s)$ by defining
$$M_{v}^{\ast}(s)= \gamma_{v}(2s-1) \cdot M_{v}(s)$$
where $\gamma_{v}(s)$ is the $\gamma$-factor attached to the trivial representation of $GL_{1}$.  In particular, $\gamma_{v}(s)=L_{v}(1-s)/L_{v}(s)$ for finite places $v$. (Note that, since $\prod_{v}\gamma_{v}(s)=1$, this normalization has no effect globally.) An advantage of this normalization is that 
$$M_{v}^{*}(1-s)\circ M_{v}^{*}(s)= \Id.$$
Moreover, we normalize the spherical vector by defining $f_{v,s}^\ast = L_v(2s) \cdot f^0_{v,s}$.  The advantage of this normalization is that  
$$M^{*}_{v}(s)(f_{v,s}^\ast)=f_{v,1-s}^\ast.$$

In order to define Eisenstein series, as in \cite{GS} page 52, we define admissible sections $f_{s}=\otimes f_{v,s}$ as follows. Let $S$ be a finite set of places containing all archimedean places. Then define $f_{v,s}= f_{v,s}^\ast$ for all $v \not \in S$ and, for $v \in S$, we take $f_{v,s}$ to be one of the two: 
\begin{enumerate}
\item $f_{v,s}$ is a constant section, i.e., its restriction to a maximal compact $K_{v}$ does not depend on $s$. 
\item $f_{v,s}=M^{\ast}_{v}(1-s)(g_{v,s})$ where $g_{v,s}$ is a constant section. 
\end{enumerate}

Note that if $f_{v,s}$ is defined by (2), then $f_{v,s}$ has a pole at $s=0$ with residue contained 
in the trivial submodule of $V(0)$.  For an admissible section $f_s$, define Eisenstein series by 
$$ E(s, g, f_{s})= \sum_{\gamma\in \alg{B}(F)\backslash \alg{GL}_{2}(F)} f_{s}(\gamma g),$$
where $\alg{B}$ is the standard Borel subgroup of upper-triangular matrices in $\alg{GL}_2$.

\subsection{Zeta integrals}

Let $\sigma=\bigotimes_{v}\sigma_{v}$ be a generic cuspidal automorphic representation of $\alg{GSp}_{6}(\adeles)$. Then, for an admissible section $f_s$, we have a zeta integral 
$$Z(s, \phi , f_{s})= \int_{\alg{Z}(\adeles)\alg{GL}_{2}(F)\backslash \alg{GL}_{2}(\adeles)}
\int_{\alg{U}(F)\backslash \alg{U}(\adeles)} \phi(\Delta(g)u) \psi_{U}(u)  E(s, g, f_{s}) ~du dg,$$
where $\phi$ is an automorphic form in the space of  $\sigma$, $f_s$ is an admissible section, $\alg{U}$ is the two-step unipotent radical of the Shalika subgroup $\alg{S}$, and $\alg{Z}$ is the center of $\alg{GL}_2$.
   
Let $W_{\phi}=\bigotimes_v W_{v}$ be the Whittaker function associated to $\phi$. The global zeta integral can be rewritten as a product of local zeta integrals $\prod_{v} Z(s, W_{v}, f_{v,s})$, where the local factor $Z(s,W_{v}, f_{v,s})$ is
\begin{eqnarray*}
\int_{\alg{B}(F_{v})\backslash \alg{GL}_{2}(F_{v})}\int_{F_{v}^{\times}}
\int_{F_{v}^{2}} 
W_{v}\left(\left( \begin{array}{cccccc} 
y &   &  &  &  &  \\
   &  y &  &  &  &  \\
   &     & y &  &  &  \\
   &   z  &  x & 1 &  &  \\
   &      &  z &    &  1 & \\
   &      &     &    &     & 1 
   \end{array}\right)  w \Delta(g)\right) & \\ 
   \times |y|^{s-3} f_{v,s}(g) ~ dx  dz  d^\times y dg &
\end{eqnarray*}
where 
\[ 
w=
\left( \begin{array}{crcccc} 
1 &  0  & 0 & 0 & 0 & 0 \\
 0 &  0 & 0 & 1 & 0 & 0 \\
 0 & -1  & 0 & 0 & 0 & 0 \\
 0 &  0 &  0 & 0 & 1 & 0 \\
 0 &  0 & 1 & 0 &  0 & 0\\
 0 &  0 &  0 &  0 &   0 & 1 
   \end{array}\right). \]

In contrast to the formula in \cite{BG} we do not have the factor $L_{v}(2s)$  as we have built it into the definition of $f_{v,s}$. If $\sigma_{v}$ is supercuspidal and $f_{v,s}$ is a constant section then the local zeta integral converges for all $s$. The following is of crucial interest to us:  assume that $v$ is finite and take  $f_{v,s}=f^{0}_{v,s}$ in the local zeta integral. Since $f^{0}_{v,0}$ is the constant function on $\alg{GL}_{2}(F_{v})$, the zeta integral at $s=0$ defines a Shalika functional. 

The following is claimed as Theorem 1 in \cite{BG}. 
\begin{prop} \label{unramified} 
Assume that $\sigma_{v}$ is unramified, and let $W_{v}$ be the corresponding (spherical) Whittaker function.  Then 
$$Z(s,W_{v}, f_{v,s}^\ast)= L(\sigma_{v}, Spin, s), $$
where the spin L-function on the right side is given by the appropriate Euler factor from the Satake parameters of $\sigma_v$.
\end{prop}

For every finite place $v$, we can now define a local $\gamma$-factor by  
$$\gamma(\sigma_{v},s) Z(s, W_{v}, f_{v,s})= Z(1-s, W_{v}, M^{*}_{v}(s)(f_{v,s})).$$
The fact that the definition of $\gamma(\sigma_{v},s)$ is independent of $W_{v}$ and $f_{v,s}$, and that it is a rational function in $q_{v}^{s}$, was proved by Vo in \cite{Vo}.  Since $M^{*}_{v}(1-s)\circ M^{*}_{v}(s)=\Id$,  we have a local functional equation 
$$\gamma(\sigma_{v},1-s)\gamma(\sigma_{v},s)=1$$
for every finite place $v$.  In particular, $\gamma(\sigma_{v}, s)$ has a pole at $s=1$ if and only if it has a zero at $s=0$.
 
Assume that $\sigma_{v}$ is unramified. Since $M^{*}_{v}(s)(f_{v,s}^\ast)=f_{v,1-s}^\ast$, Proposition \ref{unramified} implies that 
$$\gamma(\sigma_{v},s)= \frac{L(\sigma_{v},Spin,1-s)}{L(\sigma_{v},Spin,s)}.$$

\begin{prop}
\label{GammaFunctional}
Let $\sigma_{v}$ be a generic supercuspidal irrep of $GSp_6$. If $\gamma(\sigma_{v},s)$ has a zero at $s=0$ then $\sigma_{v}$ has a Shalika functional. 
\end{prop}
\proof 
Let $f_{v,s}=M^{*}_{v}(1-s) (g_{v,s})$ where $g_{v,s}$ is a constant section. Note that $f_{v,s}$ can have a pole at $s=0$ with the residue contained in the trivial subrepresentation of $V(0)$.  Consider the functional equation  
$$\gamma(\sigma_{v},s) Z(s, W_{v}, f_{v,s}) = Z(1-s, W_{v}, M^{\ast}_{v}(s)(f_{v,s})) = Z(1-s, W_{v}, g_{v,s}).$$
Since $\sigma_{v}$ is supercuspidal, the local zeta integral $Z(1-s,W_{v}, g_{v,s})$ converges for all $s$ and can be arranged to be nonzero at $s=0$ by a result of Vo (Proposition 10.4 in \cite{Vo}). 

Thus the functional equation and $\gamma(\sigma_{v},0)=0$ implies that the zeta integral $Z(s, W_{v}, f_{v,s})$ has a pole at $s=0$, for some choice of $g_{v,s}$.  After taking the residue of $f_{v,s}$ at $s=0$, the zeta integral gives a Shalika functional. 
\qed

\begin{prop} 
Let $\sigma$ be a generic supercuspidal representation of $GSp_6 = \alg{GSp}_6(k)$ with trivial central character.  Let $\gamma(\sigma, s)$ be the local factor defined above by means of zeta integrals.  Let $\gamma'(\sigma,s)$ be the analogous local factor constructed by Shahidi \cite{Sha}. Then the poles and zeros of 
$\gamma(\sigma,s)$, counted with multiplicity, coincide with poles and zeros of $\gamma'(\sigma,s)$. 
\end{prop}  
\proof The proof of this is global and uses the idea of \cite{GRS4}. Assume, as we may, that the global field $F$ contains a place $v$ such that $F_v \isom k$. Let $\Sigma$ be a global generic cuspidal automorphic representation such that $\Sigma_w$ is unramified for all finite places $w \neq v$ and $\Sigma_v \cong \sigma$.
   
The functional equation for Eisenstein series $E(s, g, f_{s})$ (\cite{GJ} page 232) implies a functional equation of the global zeta integral: 
$$Z(s,\phi, f_{s})=Z(1-s, \phi, M_{s}(f_{s})).$$
This in turn, implies  that 
$$\gamma(\Sigma_{\infty},s) \gamma(\sigma,s) = \frac{L_S(\Sigma, Spin, 1-s)}{L_S(\Sigma, Spin, s)}$$
where $S=S_{\infty}\cup\{v\}$ is the set of places consisting of  all archimedean places $S_{\infty}$ and $v$, $L_S(\Sigma, Spin, s)$ is the corresponding partial $L$-function, and 
$$\gamma(\Sigma_{\infty},s)=\prod_{w \in S_{\infty}}
\frac{Z(1-s, W_w, M^{\ast}_{w,s}(f_{w,s}))}{Z(s,W_w, f_{w,s})}.$$
We have a similar global equation satisfied by Shahidi's $\gamma$-factors. Combining the two gives 
$$\gamma(\Sigma_{\infty},s)\gamma(\sigma,s)=\gamma'(\Sigma_{\infty},s)\gamma'(\sigma,s).$$
Note that, as a consequence,  $\gamma(\Sigma_{\infty},s)$ does not depend on the choice of $f_{v,s}$. 

We need to understand the location of poles and zeros of $\gamma(\Sigma_{\infty},s)$.  Fortunately, in Proposition 12.1 of \cite{Vo}, Vo shows that for every archimedian place $w$ and every $s_0$, one can pick a constant section $f_{w,s}$ such that $Z(s_0, W_w, f_{w,s_0})\neq 0$. He also shows (see Proposition 11.1 and Lemma 11.5 in  \cite{Vo}) that the poles of the zeta integral for a constant section at archimedean places lie among the poles of $\Gamma(s_{0}+s)$ for finitely many complex numbers $s_0$. Since the poles of $M^{\ast}_{w,s}(f_{w,s})$ are contained on the real axis, it follows that poles of $\gamma(\Sigma_\infty,s)$ are located on finitely many lines parallel to the real axis. The same is true for zeros since 
$$\gamma(\Sigma_{\infty}, s) \gamma(\Sigma_{\infty},1-s)=1.$$

On the other hand, since $\gamma(\sigma,s)$ is a rational function in $q^s$, if $s_0$ is a zero or a pole then so is $s_0 + \frac{2\pi i n}{\log q}$ for every integer $n$. The same is true for Shahidi's factors -- poles and zeros of $\gamma'(\Sigma_{\infty},s)$ lie on finitely many lines parallel to the real axis, while the poles or zeros of $\gamma'(\sigma, s)$ lie on lines parallel to imaginary axis.  In view of the identity 
$$\gamma(\Sigma_{\infty},s) \gamma(\sigma, s)=\gamma'(\Sigma_\infty, s)\gamma'(\sigma,s),$$
it follows that poles and zeros of 
$\gamma(\sigma, s)$ must coincide with poles and zeros of $\gamma'(\sigma, s)$ as desired.
\qed 
 
 We can now demonstrate Theorem \ref{SpinShalika}, which is encompassed by the theorem below.
 \begin{thm}
 Let $\sigma$ be a generic supercuspidal irrep of $PGSp_6$.  Then the following three conditions are equivalent:
 \begin{enumerate}
 \item
 $\sigma$ has a nonvanishing Shalika functional.
 \item
 Shahidi's L-function $L(\sigma, Spin, s)$ has a pole at $s=0$.
 \item
 The Bump-Ginzburg-Vo L-function $L(\sigma, Spin, s)$ has a pole at $s=0$.
 \end{enumerate}
 \end{thm}
 \proof
 We prove the full circle of implications, from the results earlier in the section.
 \begin{description}
 \item[(1) implies (2)]  If $\sigma$ has a nonvanishing Shalika functional, then Shahidi's $L(\sigma, Spin, s)$ has a pole at $s=0$ by Theorem \ref{RedProp}.
 \item[(2) implies (3)]  If Shahidi's L-function $L(\sigma, Spin, s)$ has a pole at $s=0$, then Shahidi's local factor $\gamma'(\sigma, s)$ has a zero at $s=0$.  By the previous proposition, the local factor for the Bump-Ginzburg-Vo L-function $\gamma(\sigma, s)$ must also have a zero at $s=0$.  It follows that the Bump-Ginzburg-Vo L-function $L(\sigma, Spin, s)$ has a pole at $s=0$.
 \item[(3) implies (1)]  If the Bump-Ginzburg-Vo L-function $L(\sigma, Spin, s)$ has a pole at $s=0$, the local factor $\gamma(\sigma, s)$ has a zero at $s=0$.  Then Proposition \ref{GammaFunctional} implies that $\sigma$ has a nonvanishing Shalika functional.
 \end{description}
\qed

\bibliographystyle{plain}
\bibliography{Dichotomy}

\def\cprime{$'$} \def\cprime{$'$}
\begin{thebibliography}{10}

\bibitem{All2}
B.~N. Allison.
\newblock A class of nonassociative algebras with involution containing the
  class of {J}ordan algebras.
\newblock {\em Math. Ann.}, 237(2):133--156, 1978.

\bibitem{All1}
B.~N. Allison.
\newblock Models of isotropic simple {L}ie algebras.
\newblock {\em Comm. Algebra}, 7(17):1835--1875, 1979.

\bibitem{All3}
B.~N. Allison.
\newblock Tensor products of composition algebras, {A}lbert forms and some
  exceptional simple {L}ie algebras.
\newblock {\em Trans. Amer. Math. Soc.}, 306(2):667--695, 1988.

\bibitem{Asc}
Michael Aschbacher.
\newblock Chevalley groups of type {$G\sb 2$} as the group of a trilinear form.
\newblock {\em J. Algebra}, 109(1):193--259, 1987.

\bibitem{BdS}
A.~Borel and J.~De~Siebenthal.
\newblock Les sous-groupes ferm\'es de rang maximum des groupes de {L}ie clos.
\newblock {\em Comment. Math. Helv.}, 23:200--221, 1949.

\bibitem{BG}
Daniel Bump and David Ginzburg.
\newblock Spin {$L$}-functions on symplectic groups.
\newblock {\em Internat. Math. Res. Notices}, (8):153--160, 1992.

\bibitem{CS}
W.~Casselman and J.~Shalika.
\newblock The unramified principal series of {$p$}-adic groups. {II}. {T}he
  {W}hittaker function.
\newblock {\em Compositio Math.}, 41(2):207--231, 1980.

\bibitem{CKPSS}
J.~W. Cogdell, H.~H. Kim, I.~I. Piatetski-Shapiro, and F.~Shahidi.
\newblock Functoriality for the classical groups.
\newblock {\em Publ. Math. Inst. Hautes \'Etudes Sci.}, (99):163--233, 2004.

\bibitem{GS2}
Wee~Teck Gan and Gordan Savin.
\newblock Real and global lifts from {$\rm PGL\sb 3$} to {$G\sb 2$}.
\newblock {\em Int. Math. Res. Not.}, (50):2699--2724, 2003.

\bibitem{GS1}
Wee~Teck Gan and Gordan Savin.
\newblock Endoscopic lifts from {${\rm PGL}\sb 3$} to {$G\sb 2$}.
\newblock {\em Compos. Math.}, 140(3):793--808, 2004.

\bibitem{GS4}
Wee~Teck Gan and Gordan Savin.
\newblock On minimal representations definitions and properties.
\newblock {\em Represent. Theory}, 9:46--93 (electronic), 2005.

\bibitem{GJ}
Stephen Gelbart and Herv{\'e} Jacquet.
\newblock Forms of {${\rm GL}(2)$} from the analytic point of view.
\newblock In {\em Automorphic forms, representations and {$L$}-functions
  ({P}roc. {S}ympos. {P}ure {M}ath., {O}regon {S}tate {U}niv., {C}orvallis,
  {O}re., 1977), {P}art 1}, Proc. Sympos. Pure Math., XXXIII, pages 213--251.
  Amer. Math. Soc., Providence, R.I., 1979.

\bibitem{GS}
Stephen Gelbart and Freydoon Shahidi.
\newblock {\em Analytic properties of automorphic {$L$}-functions}, volume~6 of
  {\em Perspectives in Mathematics}.
\newblock Academic Press Inc., Boston, MA, 1988.

\bibitem{GeK}
I.~M. Gel{\cprime}fand and D.~A. Kajdan.
\newblock Representations of the group {${\rm GL}(n,K)$} where {$K$} is a local
  field.
\newblock In {\em Lie groups and their representations ({P}roc. {S}ummer
  {S}chool, {B}olyai {J}\'anos {M}ath. {S}oc., {B}udapest, 1971)}, pages
  95--118. Halsted, New York, 1975.

\bibitem{GJ1}
David Ginzburg and Dihua Jiang.
\newblock Periods and liftings: from {$G\sb 2$} to {$C\sb 3$}.
\newblock {\em Israel J. Math.}, 123:29--59, 2001.

\bibitem{GRS}
David Ginzburg, Stephen Rallis, and David Soudry.
\newblock A tower of theta correspondences for {$G\sb 2$}.
\newblock {\em Duke Math. J.}, 88(3):537--624, 1997.

\bibitem{GRS4}
David Ginzburg, Stephen Rallis, and David Soudry.
\newblock On a correspondence between cuspidal representations of {${\rm GL}\sb
  {2n}$} and {$\widetilde{\rm Sp}\sb {2n}$}.
\newblock {\em J. Amer. Math. Soc.}, 12(3):849--907, 1999.

\bibitem{GrS2}
David Ginzburg, Stephen Rallis, and David Soudry.
\newblock On explicit lifts of cusp forms from {${\rm GL}\sb m$} to classical
  groups.
\newblock {\em Ann. of Math. (2)}, 150(3):807--866, 1999.

\bibitem{Gro}
Benedict~H. Gross.
\newblock Some applications of {G}el\cprime fand pairs to number theory.
\newblock {\em Bull. Amer. Math. Soc. (N.S.)}, 24(2):277--301, 1991.

\bibitem{GrS1}
Benedict~H. Gross and Gordan Savin.
\newblock Motives with {G}alois group of type {$G\sb 2$}: an exceptional
  theta-correspondence.
\newblock {\em Compositio Math.}, 114(2):153--217, 1998.

\bibitem{HT}
Michael Harris and Richard Taylor.
\newblock {\em The geometry and cohomology of some simple {S}himura varieties},
  volume 151 of {\em Annals of Mathematics Studies}.
\newblock Princeton University Press, Princeton, NJ, 2001.
\newblock With an appendix by Vladimir G. Berkovich.

\bibitem{Hen}
Guy Henniart.
\newblock La conjecture de {L}anglands locale pour {${\rm GL}(3)$}.
\newblock {\em M\'em. Soc. Math. France (N.S.)}, (11-12):186, 1984.

\bibitem{Hen2}
Guy Henniart.
\newblock La conjecture de {L}anglands locale pour {${\rm GL}(p)$}.
\newblock {\em C. R. Acad. Sci. Paris S\'er. I Math.}, 299(3):73--76, 1984.

\bibitem{Hen3}
Guy Henniart.
\newblock Une preuve simple des conjectures de {L}anglands pour {${\rm GL}(n)$}
  sur un corps {$p$}-adique.
\newblock {\em Invent. Math.}, 139(2):439--455, 2000.

\bibitem{HMS}
Jing-Song Huang, Kay Magaard, and Gordan Savin.
\newblock Unipotent representations of {$G\sb 2$} arising from the minimal
  representation of {$D\sb 4\sp E$}.
\newblock {\em J. Reine Angew. Math.}, 500:65--81, 1998.

\bibitem{Jac}
N.~Jacobson.
\newblock Derivation algebras and multiplication algebras of semi-simple
  {J}ordan algebras.
\newblock {\em Ann. of Math. (2)}, 50:866--874, 1949.

\bibitem{J-R}
Herv{\'e} Jacquet and Stephen Rallis.
\newblock Uniqueness of linear periods.
\newblock {\em Compositio Math.}, 102(1):65--123, 1996.

\bibitem{JNQ}
Dihua Jiang, Chufeng Nien, and Yujun Qin.
\newblock Local {S}halika models and functoriality.
\newblock {\em Manuscripta Math.}, 127(2):187--217, 2008.

\bibitem{Kan}
I.~L. Kantor.
\newblock Certain generalizations of {J}ordan algebras.
\newblock {\em Trudy Sem. Vektor. Tenzor. Anal.}, 16:407--499, 1972.

\bibitem{Koe}
Max Koecher.
\newblock Imbedding of {J}ordan algebras into {L}ie algebras. {I}.
\newblock {\em Amer. J. Math.}, 89:787--816, 1967.

\bibitem{Kru}
Sergei Krutelevich.
\newblock Jordan algebras, exceptional groups, and {B}hargava composition.
\newblock {\em J. Algebra}, 314(2):924--977, 2007.

\bibitem{Kud}
Stephen~S. Kudla.
\newblock On the local theta-correspondence.
\newblock {\em Invent. Math.}, 83(2):229--255, 1986.

\bibitem{KM}
Philip Kutzko and Allen Moy.
\newblock On the local {L}anglands conjecture in prime dimension.
\newblock {\em Ann. of Math. (2)}, 121(3):495--517, 1985.

\bibitem{LS1}
Hung~Yean Loke and Gordan Savin.
\newblock On local lifts from {${\rm G}\sb 2(\Bbb R)$} to {${\rm Sp}\sb 6(\Bbb
  R)$} and {${\rm F}\sb 4(\Bbb R)$}.
\newblock {\em Israel J. Math.}, 159:349--371, 2007.

\bibitem{MS1}
K.~Magaard and G.~Savin.
\newblock Exceptional {$\Theta$}-correspondences. {I}.
\newblock {\em Compositio Math.}, 107(1):89--123, 1997.

\bibitem{MuS}
Goran Mui{\'c} and Gordan Savin.
\newblock Symplectic-orthogonal theta lifts of generic discrete series.
\newblock {\em Duke Math. J.}, 101(2):317--333, 2000.

\bibitem{Sav}
Gordan Savin.
\newblock Dual pair {$G_{\mathcal J} \times PGL_2$} where {$G_{\mathcal J}$ }
  is the automorphism group of the {J}ordan algebra {$\mathcal J$}.
\newblock {\em Invent. Math.}, 118(1):141--160, 1994.

\bibitem{GS3}
Gordan Savin and Wee~Teck Gan.
\newblock The dual pair {$G\sb 2\times {\rm PU}\sb 3(D)$} ({$p$}-adic case).
\newblock {\em Canad. J. Math.}, 51(1):130--146, 1999.

\bibitem{Sha}
Freydoon Shahidi.
\newblock A proof of {L}anglands' conjecture on {P}lancherel measures;
  complementary series for {$p$}-adic groups.
\newblock {\em Ann. of Math. (2)}, 132(2):273--330, 1990.

\bibitem{Vo}
San~Cao Vo.
\newblock The spin {$L$}-function on the symplectic group {${\rm GSp}(6)$}.
\newblock {\em Israel J. Math.}, 101:1--71, 1997.

\bibitem{Zha}
Yuanli Zhang.
\newblock {$L$}-packets and reducibilities.
\newblock {\em J. Reine Angew. Math.}, 510:83--102, 1999.

\end{thebibliography}
\end{document}